\documentclass[french, english, 3p]{elsarticle}
\journal{}
\biboptions{longnamesfirst, sort&compress}
\makeatletter
\def\ps@pprintTitle{%
   \def\@oddfoot{\reset@font\phantom{\today}\hfil\thepage\hfil\today}
}
\makeatother

\usepackage[utf8]{inputenc}
\usepackage[T1]{fontenc}
\usepackage{lmodern}

\usepackage{babel}
\usepackage[babel=true]{microtype}
\usepackage{amsmath, amsfonts, bm, enumitem, mathtools, url}

\usepackage{graphicx, tikz, tkz-euclide}
\usetkzobj{all}
\newcommand{\FigPath}{.}
\usepackage{siunitx}
\usepackage{hyperref, cleveref}
\hypersetup{%
   pdfauthor = {Klaus Gärtner, Lennard Kamenski},
   pdftitle = {Why do we need Voronoi cells and Delaunay meshes?
      Essential properties of the Voronoi finite volume method.},
   pdfkeywords = {finite volume method, boundary conforming Delaunay triangulation,
      Voronoi cells},
   pdfsubject = {2010 MSC: 65M08, 65M50},
}
\crefname{equation}{}{}
\crefname{figure}{Fig.}{Figs.}
\Crefname{figure}{Figure}{Figures}

\providecommand{\abs}[1]{\lvert #1 \rvert}
\providecommand{\D}[1]{\mathrm{d}#1}
\providecommand{\R}{\mathbb{R}}
\providecommand{\Th}{\mathcal{T}_h}

\begin{document}

%--- Title frame -----------------------------------------------------
\renewcommand*{\today}{\emph{\small{}July 31, 2019}}

\begin{frontmatter}
\title{%
   Why do we need Voronoi cells and~Delaunay meshes?\\
   Essential properties of the Voronoi finite volume method.%
}

\author{Klaus Gärtner}
\author{Lennard Kamenski}
\address{m4sim GmbH, Berlin, Germany (\url{https://m4sim.de})}

\begin{abstract}%
Unlike other schemes that locally violate the essential stability properties of the analytic parabolic and elliptic problems, Voronoi finite volume methods (FVM) and boundary conforming Delaunay meshes provide good approximation of the geometry of a problem and are able to preserve the essential qualitative properties of the solution for any given resolution in space and time as well as changes in time scales of multiple orders of magnitude.
This work provides a brief description of the essential and useful properties of the Voronoi FVM, application examples, and a motivation why Voronoi FVM deserve to be used more often in practice than they are currently.
\\[0.5\baselineskip]
{\footnotesize{}%
\textcopyright~2019 The Authors.
Preprint of the Work accepted for publication in CMMP
by Pleiades Publishing (\url{http://pleiades.online}).}
\end{abstract}


\end{frontmatter}

\section{Introduction}

Many application problems result in highly nonlinear elliptic or parabolic partial differential equations.
Many of them have been studied for a long time and analytic results are available, which are often based on a variational formulation.
Using entropy or free energy functionals allows the derivation of solution properties and decay properties of the free energy without restrictive smallness assumption on the free energy or the boundary conditions.
Typically, a~priori bounds for the solutions, such as positivity and boundedness, have to be derived and decay of the free energy along trajectories can be shown.

In the discrete case, it is less obvious.
Often the discretization parameters in space and time, $h$ and $\tau$, which are in fact limits, are regarded as finite and universally achievable.
However, this is not always the case in situations with limited smoothness, such as in states further away from the thermodynamic equilibrium or in case of higher free energy densities, and particularly not, if the approximation of the non-smooth limits is required.
Solutions may have plateaus and internal or boundary layers connecting
the plateaus and the boundary values.
Energy densities are often very large if compared to a typical human environment. 
A well-known example is the charge transport in semiconductors, where the large difference in force results from the fact that the  difference of the electrostatic potential of merely \SI{1}{\volt} corresponds to a height difference of
\SI{4000}{\km} in the gravitational field of the earth (mechanical heat
equivalent relation at \SI{300}{\kelvin}); accordingly, \SI{100}{\volt}, which is not much for such problems, correspond to a height difference from the earth to somewhere behind the moon.
Many questions in electro-chemistry, chemical reactions, and other areas result in the same type of problems and require the preservation of qualitative properties in their discrete form.

In such cases, the classic discretization error has only a limited significance with respect to the guarantee of a stable discreet solution.
It is often even physically irrelevant whether a stationary limit value at the smallest time scale of \SI{1}{\second} of the problem is reached at $10^{19}$ or $10^{20}$, since the relevant human response time and life span are \SI{1}{\second} and \SI{100}{years}, respectively ($\SI{100}{years} \approx 100 \times 400 \times 25 \times \SI{3600}{\second} = \SI{3.6e+9}{s}$).
In contrast, of great interest are stationary states and qualitative statements, such as whether design A provides a better solution in the \si{\giga\hertz} range then design B.
Therefore, preservation of the structural properties of the analytical problems in space and time for finite and not too small discretization parameters is an essential.

Moreover, often, the provided data is not precise and internal layers are 
separating
volumes with different properties (generation on one side of the layer,
recombination on the other). 
To estimate these volumes within the order of
percent of the total volume may be sufficient.
What is often needed is the position of the layer but not its resolution.
As long as the discretization scheme inherits
the stability properties of the analytic problem independently of the
spatial and time step sizes, the requirement to resolve everything
may be not adequate. % in order to solve these problems.

It is not a priori clear which combination of space-time discretization will provide unconditionally stable
discrete equations sufficiently well approximating the original problem or whether such a scheme is efficient.
Voronoi-Delaunay meshes and implicit Euler time discretization play a special role here.

Stable discretizations for reduced problems are dating back to
1950s~\cite{AllenSouthwell55,Ilin69,Sch_Gu} and are still relevant today
in a generalized setting.  \emph{Voronoi cells}~\cite{Voronoi07}
and their dual, \emph{Delaunay triangulation}~\cite{Delaunay}, have
been used as early as~\cite{Macneal53} (e.g.,\ the classical book
by Varga~\cite{Varga62}).
Delaunay triangulations are optimal for a variety of criteria, e.g.,\
maximizing the smallest angle~\cite{Lawson77}, minimizing the
\emph{roughness} of the mesh~\cite{Rip90} (the squared
$L^2$ norm of the gradient of the piecewise linear interpolating
surface) in two dimensions, and other (see, e.g.,~\cite[Sect.~4]{She08}
and the references therein).
It is hard to follow the history of the development of Delaunay
meshes in detail but the central points have been related to $C^1$
surface interpolation and to applications such as neutron transport
and charge transport in semiconductors.
\emph{Admissible meshes}~\cite{EymardGallouetHerbin97} summarize
the minimal mesh attributes required to prove convergence.

Vienna may have been the birthplace of the term \emph{boundary conforming
Delaunay} meshes (e.g.,~comp.~\cite{PFleischmann-99}),
which are required in case of interfaces with jumps in material properties (e.g.,
diffusion coefficients), or nonlinear boundary conditions of the
third kind, or some other special cases to ensure the boundary and
internal interfaces of the problem domain to be a part of the
Voronoi / Delaunay mesh, whereas the \emph{constraint Delaunay}
property~\cite{Chew89}, which gained popularity recently~\cite{She02b,Si15},
is not really helpful.
Local orthogonal coordinates would be the best choice %wonderful
and, although they will often
not exist globally, locally they might be the key to combine stability,
good approximation properties, and low complexity of the discrete problem.

The following list summarizes some desired properties of a discrete
problem in relation to its analytic counterpart:
\setlist{nolistsep}
\begin{itemize}[noitemsep]
\item the same qualitative stability properties as the original, analytic problem,
\item no artificial smallness or smoothness assumptions,
\item a weak form of the discrete problem, known analytic test functions and techniques,
\item hence, a weak discrete energy functional,
\item some weak convergence for distributions on the right-hand side,
\item internal layer positions as bounded jumps (unresolved layers) in the solution.
\end{itemize}
Some examples, where the list is partially fulfilled, are: transport of electrons and holes in semiconductors and phase separation (non-Cahn-Hilliard)~\cite{GajGar05a,GajGar05,Gar09,ga15,glgae-09,GajGri06}.

%In well understood analytic cases, one has a given functional, e.g., describing the free energy, and a corresponding parabolic evolution equation with solutions usually in the corresponding Sobolev space $W^{1,2}$ with limited energy for this kind of equations.

%On the other side of the time scale, charge transport in semiconductors such as silicon has a fast phonon scattering (\SI{e-13}{s}), charge carrier thermalization on wavelength scales of about \SI{1}{\mm} is about \SI{e-12}{s}, and the drift-diffusion equations are valid only up to about \SI{e-10}{s}, if the volumes are sufficiently large compared with wavelengths and free path lengths.
%Thus, paths of high-energy particles on the \si{\mm} scale can be considered as instantaneous and are not described by the Nernst-Plack-type parabolic equations.

%Semiconductor equations are a well-known example but there are considerably more problems exhibiting similar structure, however, not in our usual, macroscopic living environment of biological objects.

%--- Definitions ----------------------------------------
\section{Delaunay meshes, Voronoi volumes, and boundary conforming Delaunay}

In the following, we consider a discretization (\emph{mesh}) $\Th(\Omega)$ of a domain $\Omega \subset \mathbb{R}^d$ and assume that $\Omega$ is a union of finite $d$-dimensional subdomains, $\Omega = \cup_i \Omega_i$.
Each subdomain contains only one material and consists of $d$-dimensional simplices $\bm{E}^d_j$, $\Omega_i = \cup_j \bm{E}^d_j$.
The vertices of the mesh simplices are denoted by $\bm{v}_{k}$.
Consequently, all interfaces between materials and boundaries of the problem domain are represented by low-dimensional simplices $\bm{E}^m_{k}$, $1 \le m < d$.
It is assumed that the domain geometry representation contains all critical lines and points.

%--- Delaunay ----------------------------------------
\subsection{Delaunay triangulations}
A discretization (mesh) $\Th(\Omega)$ of a given $d$-dimensional domain $\Omega$ by simplices is called a \emph{Delaunay mesh} or \emph{Delaunay triangulation} if it fulfills the empty circumball property, that is, if no mesh vertex is inside the circumball of any simplex, i.e.,
\[
   \forall \bm{v}_i, \bm{E}^d_j
   \in \Th(\Omega),
   ~
   \bm{v}_i \not\in \mathrm{B^o}(\bm{E}^d_j)
   .
\]
For a given mesh, a mesh simplex is called \emph{Delaunay} if its circumball is empty (\cref{fig:delaunay}).
A mesh edge is called \emph{Delaunay} if there exist an empty sphere containing this edge but no other mesh vertices inside.
Consequently, a simplex is a Delaunay simplex iff all of its edges are Delaunay.
Note that the Delaunay edges are not necessarily the shortest edges possible (see \cref{fig:delaunay}, right)
but a Delaunay triangulation maximises the minimum angle of a triangulation~\cite{Lawson77} and it is a minimal roughness triangulation in two dimensions~\cite{Rip90}.

A Delaunay triangulation is unique if no $d+2$ vertices lie on a common empty $d$-sphere.

\begin{figure}[t]%
   \centering{}%
   \begin{tikzpicture}[scale = 0.75]%
      \useasboundingbox (-0.69, -1.67) rectangle (5.15,3.86);
      % define vertices 1, 2, 3, 4
      \tkzDefPoint(0.0,  0.0){p1}
      \tkzDefPoint(4.5, -0.3){p2}
      \tkzDefPoint(2.8,  3.2){p3}
      \tkzDefPoint(0.1,  2.5){p4}
      % triangles 1-2-3 (gray) & 2-3-4
      \tkzFillPolygon[fill = gray!15](p1,p2,p3)
      \tkzDrawPolygon[thin](p1,p2,p3,p4)
      \tkzDrawSegment[thin](p1,p3)
      % circumcircle through 1-2-3
      \begin{pgfinterruptboundingbox}
      \tkzCircumCenter(p1,p2,p3)\tkzGetPoint{c}
      \end{pgfinterruptboundingbox}
      \tkzDrawCircle[dotted](c,p1)
      % draw vertices
      \tkzDrawPoints(p1,p2,p3,p4)
      \tkzLabelPoint[left](p1){$\bm{v}_1$}
      \tkzLabelPoint[below right](p2){$\bm{v}_2$}
      \tkzLabelPoint[above](p3){$\bm{v}_3$}
      \tkzLabelPoint[above left](p4){$\bm{v}_4$}
   \end{tikzpicture}%
   \hfill{}%
   \begin{tikzpicture}[scale = 0.75]%
      \useasboundingbox (-0.69, -1.67) rectangle (5.15,3.86);
      % define vertices 1, 2, 3, 4
      \tkzDefPoint(0.0,  0.0){p1}
      \tkzDefPoint(4.5, -0.3){p2}
      \tkzDefPoint(2.8,  3.2){p3}
      \tkzDefPoint(0.1,  2.5){p4}
      % triangles 1-2-4 (gray) & 2-3-4
      \tkzFillPolygon[fill = gray!15](p1,p2,p4)
      \tkzDrawPolygon[thin](p1,p2,p3,p4)
      \tkzDrawSegment[thin](p2,p4)
      % circumcircle through 1-2-4
      \begin{pgfinterruptboundingbox}
      \tkzCircumCenter(p1,p2,p4)\tkzGetPoint{c}
      \end{pgfinterruptboundingbox}
      \tkzDrawCircle[dotted](c,p1)%
      % draw vertices
      \tkzDrawPoints(p1,p2,p3,p4)
      \tkzLabelPoint[left](p1){$\bm{v}_1$}
      \tkzLabelPoint[below right](p2){$\bm{v}_2$}
      \tkzLabelPoint[above](p3){$\bm{v}_3$}
      \tkzLabelPoint[above left](p4){$\bm{v}_4$}
   \end{tikzpicture}%
   \hfill{}%
   \begin{tikzpicture}[scale = 0.75]%
      %\useasboundingbox (-0.69, -1.67) rectangle (5.15,3.86);
      \useasboundingbox (-0.69, -1.67) rectangle (5.15,4.46);
      % define vertices 1, 2, 3, 4
      \tkzDefPoint( 0.0,  1.6){p1}
      \tkzDefPoint( 4.5, -0.3){p2}
      \tkzDefPoint( 2.3,  4.0){p3}
      \tkzDefPoint( 1.0,  3.3){p4}
      % triangles 1-2-4 (gray) & 2-3-4
      \tkzFillPolygon[fill = gray!15](p1,p2,p4)
      \tkzDrawPolygon[thin](p1,p2,p3,p4)
      \tkzDrawSegment[thin](p2,p4)
      \tkzDrawSegment[dashed](p1,p3)
      % circumcircle through 1-2-4
      \begin{pgfinterruptboundingbox}
      \tkzCircumCenter(p1,p2,p4)\tkzGetPoint{c}
      \end{pgfinterruptboundingbox}
      \tkzDrawCircle[dotted](c,p1)%
      % draw vertices
      \tkzDrawPoints(p1,p2,p3,p4)
      \tkzLabelPoint[left](p1){$\bm{v}_1$}
      \tkzLabelPoint[below right](p2){$\bm{v}_2$}
      \tkzLabelPoint[above](p3){$\bm{v}_3$}
      \tkzLabelPoint[above left](p4){$\bm{v}_4$}
   \end{tikzpicture}%
   \caption{%
      $\triangle \bm{v}_1 \bm{v}_2 \bm{v}_3$ is Delaunay
      because its circumball is empty (left), whereas
      $\triangle \bm{v}_1 \bm{v}_2 \bm{v}_4$ is not Delaunay
      because its circumball contains $\bm{v}_3$ (center).
      The Delaunay edge may not be the shortest edge:
      $\bm{v}_2 \bm{v}_4$ is Delaunay but not the shortest (right)%
   }\label{fig:delaunay}%
\end{figure}

%--- Voronoi volumes etc. ----------------------------------------
\subsection{Voronoi volumes, surfaces, and cells}

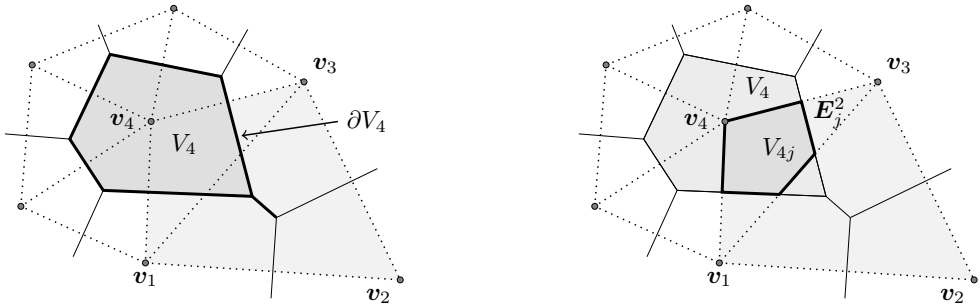
\begin{figure}[t]%
   \centering{}%
   \begin{tikzpicture}[scale = 0.8]%
      % define vertices 1, 2, 3, 4
      \tkzDefPoint(0.0,  0.0){p1}
      \tkzDefPoint(0.1,  2.5){p2}
      \tkzDefPoint(2.8,  3.2){p3}
      \tkzDefPoint(4.5, -0.3){p4}
      \tkzDefPoint( 0.5,  4.5){p5}
      \tkzDefPoint(-2.0,  3.5){p6}
      \tkzDefPoint(-2.2,  1.0){p7}
      \tkzDefPoint( 0.3,  2.1){labelVV}
      \tkzDefPoint( 3.4,  2.5){labelVS}
      % circumcenters
      \begin{pgfinterruptboundingbox}
      \tkzCircumCenter(p1,p2,p3)\tkzGetPoint{c123}
      \tkzCircumCenter(p1,p3,p4)\tkzGetPoint{c134}
      \tkzCircumCenter(p2,p3,p5)\tkzGetPoint{c235}
      \tkzCircumCenter(p2,p5,p6)\tkzGetPoint{c256}
      \tkzCircumCenter(p2,p6,p7)\tkzGetPoint{c267}
      \tkzCircumCenter(p2,p7,p1)\tkzGetPoint{c271}
      \end{pgfinterruptboundingbox}
      % midpoints
      \tkzDefMidPoint(p1,p4)\tkzGetPoint{mid14}
      \tkzDefMidPoint(p3,p4)\tkzGetPoint{mid34}
      \tkzDefMidPoint(p3,p5)\tkzGetPoint{mid35}
      \tkzDefMidPoint(p5,p6)\tkzGetPoint{mid56}
      \tkzDefMidPoint(p6,p7)\tkzGetPoint{mid67}
      \tkzDefMidPoint(p7,p1)\tkzGetPoint{mid71}
      \tkzDefMidPoint(c235,c123)\tkzGetPoint{vs}
      % voronoi volumes
      \tkzFillPolygon[fill = gray!10](p1,p2,p3,p4)
      \tkzFillPolygon[fill = gray!25](c123,c235,c256,c267,c271)
      \tkzDrawPolygon[very thick](c123,c235,c256,c267,c271)
      \tkzDrawLines[very thick, add = 0 and 0](c123,c134)
      \tkzDrawLines[thin, add = 0 and 1/2](c134,mid14)
      \tkzDrawLines[thin, add = 0 and 1/3](c134,mid34)
      \tkzDrawLines[thin, add = 0 and 1/2](c235,mid35)
      \tkzDrawLines[thin, add = 0 and 2/3](c256,mid56)
      \tkzDrawLines[thin, add = 0 and 1/2](c267,mid67)
      \tkzDrawLines[thin, add = 0 and 1/2](c271,mid71)
      % triangles 1-2-4 (dotted)
      \tkzDrawPolygon[dotted](p1,p4,p3,p5,p6,p7)
      \tkzDrawLines[  dotted, add = 0 and 0](p2,p1)
      \tkzDrawLines[  dotted, add = 0 and 0](p2,p3)
      \tkzDrawLines[  dotted, add = 0 and 0](p1,p3)
      \tkzDrawLines[  dotted, add = 0 and 0](p2,p5)
      \tkzDrawLines[  dotted, add = 0 and 0](p2,p6)
      \tkzDrawLines[  dotted, add = 0 and 0](p2,p7)
      % draw vertices
      \tkzDrawPoints(p1,p2,p3,p4,p5,p6,p7)
      \tkzLabelPoint[right](labelVV){$V_4$}
      \tkzLabelPoint[right](labelVS){$\partial V_4$}
      \tkzDrawSegment[add = 0 and -1/20, arrows = ->](labelVS, vs)
      \tkzLabelPoint[below](p1){$\bm{v}_1$}
      \tkzLabelPoint[left = 1mm](p2){$\bm{v}_4$}
      \tkzLabelPoint[above right](p3){$\bm{v}_3$}
      \tkzLabelPoint[below left](p4){$\bm{v}_2$}
   \end{tikzpicture}%
   \hspace{0.15\linewidth}%
   \begin{tikzpicture}[scale = 0.8]%
      % define vertices 1, 2, 3, 4
      \tkzDefPoint(0.0,  0.0){p1}
      \tkzDefPoint(0.1,  2.5){p2}
      \tkzDefPoint(2.8,  3.2){p3}
      \tkzDefPoint(4.5, -0.3){p4}
      \tkzDefPoint( 0.5,  4.5){p5}
      \tkzDefPoint(-2.0,  3.5){p6}
      \tkzDefPoint(-2.2,  1.0){p7}
      \tkzDefPoint( 0.3,  3.1){labelVV}
      \tkzDefPoint( 0.6,  2.0){labelVS}
      \tkzDefPoint( 1.5,  2.65){labelEj}
      % circumcenters
      \begin{pgfinterruptboundingbox}
      \tkzCircumCenter(p1,p2,p3)\tkzGetPoint{c123}
      \tkzCircumCenter(p1,p3,p4)\tkzGetPoint{c134}
      \tkzCircumCenter(p2,p3,p5)\tkzGetPoint{c235}
      \tkzCircumCenter(p2,p5,p6)\tkzGetPoint{c256}
      \tkzCircumCenter(p2,p6,p7)\tkzGetPoint{c267}
      \tkzCircumCenter(p2,p7,p1)\tkzGetPoint{c271}
      \end{pgfinterruptboundingbox}
      % midpoints
      \tkzDefMidPoint(p1,p2)\tkzGetPoint{mid12}
      \tkzDefMidPoint(p1,p4)\tkzGetPoint{mid14}
      \tkzDefMidPoint(p3,p4)\tkzGetPoint{mid34}
      \tkzDefMidPoint(p3,p5)\tkzGetPoint{mid35}
      \tkzDefMidPoint(p5,p6)\tkzGetPoint{mid56}
      \tkzDefMidPoint(p6,p7)\tkzGetPoint{mid67}
      \tkzDefMidPoint(p7,p1)\tkzGetPoint{mid71}
      \tkzDefMidPoint(c235,c123)\tkzGetPoint{vs}
      \tkzDefMidPoint(p2,p3)\tkzGetPoint{mid23}
      \tkzDefMidPoint(p2,p4)\tkzGetPoint{mid24}
      % intersections
      \tkzInterLL(p1,p3)(c271,c123)\tkzGetPoint{i1}
      \tkzInterLL(p1,p3)(c123,c235)\tkzGetPoint{i2}
      % voronoi volume
      \tkzFillPolygon[fill = gray!10](p1,p2,p3,p4)
      \tkzFillPolygon[fill = gray!15](c123,c235,c256,c267,c271)
      \tkzDrawPolygon[thin](c123,c235,c256,c267,c271)
      \tkzDrawPolygon[dotted](p1,p2,p3)
      \tkzDrawPolygon[thin, ](c123,c235,c256,c267,c271)
      \tkzFillPolygon[fill = gray!25](p2,mid12,i1,i2,mid23)
      \tkzDrawPolygon[very thick](p2,mid12,i1,i2,mid23)
      \tkzDrawLines[thin, add = 0 and 0](c123,c134)
      \tkzDrawLines[thin, add = 0 and 1/2](c134,mid14)
      \tkzDrawLines[thin, add = 0 and 1/3](c134,mid34)
      \tkzDrawLines[thin, add = 0 and 1/2](c235,mid35)
      \tkzDrawLines[thin, add = 0 and 2/3](c256,mid56)
      \tkzDrawLines[thin, add = 0 and 1/2](c267,mid67)
      \tkzDrawLines[thin, add = 0 and 1/2](c271,mid71)
      % triangles 1-2-4 (dotted)
      \tkzDrawPolygon[dotted](p1,p4,p3,p5,p6,p7)
      \tkzDrawLines[  dotted, add = 0 and 0](p2,p5)
      \tkzDrawLines[  dotted, add = 0 and 0](p2,p6)
      \tkzDrawLines[  dotted, add = 0 and 0](p2,p7)
      % draw vertices
      \tkzDrawPoints(p1,p2,p3,p4,p5,p6,p7)
      \tkzLabelPoint[right](labelVV){$V_4$}
      \tkzLabelPoint[right](labelVS){$V_{4j}$}
      \tkzLabelPoint[right](labelEj){$\bm{E}^2_j$}
      %\tkzDrawSegment[add = 0 and 1/10, arrows = ->](labelVS, vs)
      %\tkzLabelPoint[left = 1mm](p2){$\bm{v}_i$}
      \tkzLabelPoint[below](p1){$\bm{v}_1$}
      \tkzLabelPoint[left = 1mm](p2){$\bm{v}_4$}
      \tkzLabelPoint[above right](p3){$\bm{v}_3$}
      \tkzLabelPoint[below left](p4){$\bm{v}_2$}
      %\tkzLabelPoint[above](p5){$\bm{v}_5$}
      %\tkzLabelPoint[left](p6){$\bm{v}_6$}
      %\tkzLabelPoint[below left](p7){$\bm{v}_7$}
   \end{tikzpicture}%
   \caption{%
      Voronoi volume $V_4$,
      Voronoi surface $\partial V_4$ (left), and the
      Voronoi cell $V_4j = V_4 \cap \bm{E}^2_j$
      (right)%
   }\label{fig:voronoi}%
\end{figure}

The \emph{Voronoi volume} $V_i$ of a vertex $\bm{v}_i$ is the region of $\mathbb{R}^d$ consisting of all points closer to $\bm{v}_i$ than to any other mesh vertex,
\[
   V_i := \{ \bm{x} \in \mathbb{R}^d :
      \|\bm{x} - \bm{v}_i\| < \|\bm{x} - \bm{v}_j \|
      ~
      \forall \bm{v}_j \in \Th(\Omega) \}
      .
 \]
Its boundary $\partial V_i$ is the corresponding \emph{Voronoi
surface} (\cref{fig:voronoi}, left).

The Voronoi decomposition for a set of vertices is dual to its Delaunay triangulation.
Each Voronoi facet corresponds to an edge of the Delaunay triangulation.

The intersection of a Voronoi volume $V_i$
with a simplex $\bm{E}^d_j$, $V_{ij} = V_i \cap
\bm{E}^d_j$, is the \emph{Voronoi volume element (cell)} of
$\bm{v}_i$ with respect to $\bm{E}^d_j$
(\cref{fig:voronoi}, right).
Being an intersection of half-spaces, Voronoi cells and
facets are always convex  and their measures are continuous functions
of vertex positions.

%
%\end{definition}

%--- BCD ----------------------------------------
\subsection{Boundary conforming Delaunay triangulation}

\begin{figure}[t]%
   \begin{tikzpicture}[scale = 0.85]%
      \useasboundingbox (-1.23, -1.07) rectangle (4.78, 4.27);
      % define vertices
      \tkzDefPoint(0.0,  0.0){p1}
      \tkzDefPoint(0.1,  2.5){p2}
      \tkzDefPoint(2.8,  3.2){p3}
      \tkzDefPoint(4.5, -0.3){p4}
      \begin{pgfinterruptboundingbox}
      \tkzDefPointBy[reflection=over p1--p3](p2)\tkzGetPoint{p22}
      \tkzCircumCenter(p1,p2,p3) \tkzGetPoint{c}
      \end{pgfinterruptboundingbox}
      \tkzInterLL(p1,p3)(p2,p22) \tkzGetPoint{i13}
      % triangles, boundary segments
      \tkzDrawPolygon[thin, fill = gray!15](p1,p2,p3)
      \tkzDrawLines[very thick, add = 0 and 1/7](p3,p1)
      \tkzDrawLines[very thick, add = 0 and 1/7](p3,p4)
      % reflected point
      \tkzDrawLines[thin, dashed,  add = 0 and 0](p2,p22)
      \tkzDrawSegments[thin, dashed](p1,p22 p22,p3)
      \tkzDrawSegments[thin, dashed](p2,i13 i13,p22)
      \tkzMarkSegments[mark = ||](p2,i13 p22,i13)
      \tkzMarkRightAngle(p2,i13,p3)
      % circumcircle
      \tkzDrawCircle[dotted](c,p1)
      % draw vertices
      \tkzDrawPoints(p1,p2,p3,p4,p22)
      \tkzLabelPoint[left = 0.07](p1){$\bm{v}_1$}
      \tkzLabelPoint[above left](p2){$\bm{v}_4$}
      \tkzLabelPoint[above](p3){$\bm{v}_3$}
      \tkzLabelPoint[below left](p4){$\bm{v}_2$}
      \tkzLabelPoint[right](p22){$\bm{v}_4^*$}
      \end{tikzpicture}%
   \hfill{}%
   \begin{tikzpicture}[scale = 0.85]%
      \useasboundingbox (-1.23, -1.07) rectangle (4.78, 4.27);
      % vertices
      \tkzDefPoint(0.0,  0.0){p1}
      \tkzDefPoint(0.1,  2.5){p2}
      \tkzDefPoint(2.8,  3.2){p3}
      \tkzDefPoint(4.5, -0.3){p4}
      % circumceters of 1-2-3 & 1-3-4
      \begin{pgfinterruptboundingbox}
      \tkzCircumCenter(p1,p2,p3)\tkzGetPoint{c123}
      \tkzCircumCenter(p1,p3,p4)\tkzGetPoint{c134}
      \tkzDefPointBy[projection=onto p1--p3](p2)\tkzGetPoint{p22}
      \end{pgfinterruptboundingbox}
      % edge midpoints
      \tkzDefMidPoint(p1,p2)\tkzGetPoint{mid12}
      \tkzDefMidPoint(p1,p3)\tkzGetPoint{mid13}
      \tkzDefMidPoint(p2,p3)\tkzGetPoint{mid23}
      \tkzDefMidPoint(p3,p4)\tkzGetPoint{mid34}
      \tkzDefMidPoint(p1,p4)\tkzGetPoint{mid14}
      % tdraw riangles and boundary segments
      \tkzDrawPolygon[thin, fill = gray!15](p1,p2,p3)
      \tkzDrawLines[very thick, add = 0 and 1/7](p3,p1)
      \tkzDrawLines[very thick, add = 0 and 1/7](p3,p4)
      % circumball of 1-3
      \tkzDrawCircle[dotted](mid12,p2)
      \tkzDrawCircle[dotted](mid23,p2)
      % admissible vertex
      \tkzDrawLines[dashed,  add = 0 and 0](p2,p22)
      \tkzMarkRightAngle(p2,p22,p3)
      % draw Voronoi volumes
      \tkzDrawLines[thin, add = 0 and 0](c123,c134)
      \tkzDrawLines[thin, add = 0 and 1/2](c134,mid14)
      \tkzDrawLines[thin, add = 0 and 1/2](c134,mid34)
      \tkzDrawLines[thin, add = 0 and 1/3](c123,mid12)
      \tkzDrawLines[thin, add = 0 and 1/3](c123,mid23)
      % draw vertices
      \tkzDrawPoints(p1,p2,p3,p4)
      \tkzLabelPoint[left = 0.07](p1){$\bm{v}_1$}
      \tkzLabelPoint[above left](p2){$\bm{v}_4$}
      \tkzLabelPoint[above](p3){$\bm{v}_3$}
      \tkzLabelPoint[below left](p4){$\bm{v}_2$}
   \end{tikzpicture}%
   \hfill{}%
   \begin{tikzpicture}[scale = 0.85]%
      \useasboundingbox (-1.23, -1.07) rectangle (4.78, 4.27);
      % vertices
      \tkzDefPoint(0.0,  0.0){p1}
      \tkzDefPoint(0.1,  2.5){p2}
      \tkzDefPoint(2.8,  3.2){p3}
      \tkzDefPoint(4.5, -0.3){p4}
      % edge midpoints and projection
      \tkzDefMidPoint(p1,p2)\tkzGetPoint{mid12}
      \tkzDefMidPoint(p2,p3)\tkzGetPoint{mid23}
      \tkzDefMidPoint(p3,p4)\tkzGetPoint{mid34}
      % circumceters of 1-2-3 & 1-3-4
      \begin{pgfinterruptboundingbox}
      \tkzDefPointWith[orthogonal](mid12,p2)\tkzGetPoint{h12}
      \tkzDefPointWith[orthogonal](mid23,p3)\tkzGetPoint{h23}
      \tkzDefPointWith[orthogonal](mid34,p4)\tkzGetPoint{h34}
      \tkzDefPointBy[projection=onto p1--p3](p2)\tkzGetPoint{p22}
      \tkzDefPointBy[projection=onto p1--p3](mid12)\tkzGetPoint{h122}
      \tkzDefPointBy[projection=onto p1--p3](mid23)\tkzGetPoint{h232}
      \end{pgfinterruptboundingbox}
      % tdraw riangles and boundary segments
      \tkzDrawPolygon[thin, fill = gray!15](p1,p2,p3)
      \tkzDrawLines[         add = 0 and 0](p2,p22)
      \tkzDrawLines[very thick, add = 0 and 1/2](h122,p1)
      \tkzDrawLines[very thick, add = 0 and 0](h232,p3)
      \tkzDrawLines[very thick, dashed, add = 0 and 0](h122,h232)
      \tkzDrawLines[very thick, add = 0 and 1/7](p3,p4)
      % circumballs
      \tkzDrawCircle[dotted](h122,p1)
      \tkzDrawCircle[dotted](h232,p3)
      % draw Voronoi volumes
      \tkzDrawLines[thin, add = 0 and -0.5](mid12,h12)
      \tkzDrawLines[thin, add = 0 and -0.6](mid23,h23)
      \tkzDrawLines[thin, add = 0 and -0.5](mid34,h34)
      \tkzDrawLines[thin, add = 0 and  0  ](mid12,mid23)
      \tkzDrawLines[thin, add = 0 and  0  ](mid12,h122)
      \tkzDrawLines[thin, add = 0 and  0  ](mid23,h232)
      \tkzMarkRightAngle(p2,p22,p3)
      % draw vertices
      \tkzDrawPoints(p1,p2,p3,p4,p22)
      \tkzLabelPoint[left = 0.07](p1){$\bm{v}_1$}
      \tkzLabelPoint[above left](p2){$\bm{v}_4$}
      \tkzLabelPoint[above](p3){$\bm{v}_3$}
      \tkzLabelPoint[below left](p4){$\bm{v}_2$}
      \tkzLabelPoint[above right = -2mm and 1mm](p22){$\tilde{\bm{v}}_4$}
   \end{tikzpicture}%
   \caption{%
      If $\angle \bm{v}_1 \bm{v}_4 \bm{v}_3$ is obtuse,
      the ghost vertex $\bm{v}_4^*$ to treat the Neumann boundary
      condition will violate the Delaunay property (left).
      Iserting a new vertex as the orthogonal projection of $\bm{v}_4$
      onto $\bm{v}_1 \bm{v}_3$ will not violate the Delaunay
      property of $\bm{v}_1 \bm{v}_4$ and $\bm{v}_4 \bm{v}_3$ (center)
      while restoring the boundary conformity of the mesh (right);
      the dashed bold line corresponds to the boundary condition
      influencing $\tilde{\bm{v}}_4$%
   }\label{fig:bcd}%
\end{figure}
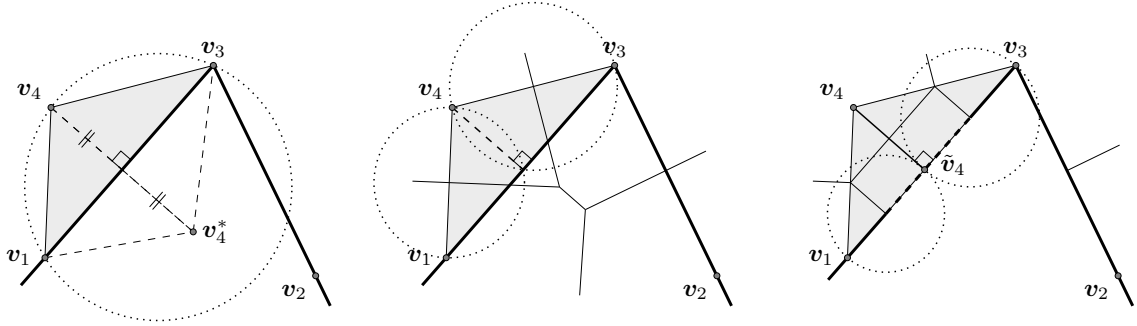

Introducing boundary or interface conditions into the mesh might violate the Delaunay property.
For instance, consider a homogeneous Neumann boundary or interface condition in \cref{fig:voronoi} along the edges $\bm{v}_1 \bm{v}_3$ and $\bm{v}_3 \bm{v}_2$ (\cref{fig:bcd} (left)).
In cell centered schemes, a common technique to treat the Neumann boundary is to introduce a ghost vertex $\bm{v}_4^*$ as a reflection of $\bm{v}_4$ over $\bm{v}_1 \bm{v}_3$.
If the ghost vertex $\bm{v}_4^*$ is inside the circumball of $\triangle \bm{v}_1 \bm{v}_3 \bm{v}_4$ (which will be the case iff $\bm{v}_4$ is inside the circumball of $\bm{v}_1 \bm{v}_3$) then $\triangle \bm{v}_1 \bm{v}_4 \bm{v}_3$ becomes non-Delaunay.
The problem extends to higher dimensions as well and motivates the following definition.

%--- Problem with Delaunay and interfaces ------------------------------

A Delaunay mesh is called \emph{boundary conforming Delaunay} if no mesh vertex is inside the smallest (diametrical) circumball of any low-dimensional boundary or interface simplex, i.e.,
\[
   \forall \bm{v}_i \in \Th(\Omega),
   \forall \bm{E}^m_j \in \partial \Omega_i,
   ~
   \bm{v}_i \not\in \mathrm{B^o}(\bm{E}^m_j)
   \quad(1 \le m < d)
   .
\]
The property of an empty diametrical ball is sometimes referred to as the \emph{Gabriel} property~\cite{GabSok69}.
In brief, a boundary conforming Delaunay mesh is a mesh, where each simplex is Delaunay and each boundary simplex is Gabriel.

In the situation of \cref{fig:bcd} (left), boundary conformity can be 
enforced by splitting $\triangle \bm{v}_1 \bm{v}_3 \bm{v}_4$ with a new
vertex on the boundary edge  $\bm{v}_1 \bm{v}_3$.
Taking the intersection of the diametrical disks of $\bm{v}_1 \bm{v}_4$
and $\bm{v}_4 \bm{v}_3$ with the boundary edge $\bm{v}_1 \bm{v}_3$,
which is exactly the projection of $\bm{v}_4$ onto $\bm{v}_1 \bm{v}_3$
(\cref{fig:bcd}, center), guarantees that all other mesh simplices remain Delaunay
while restoring the boundary conformity of the mesh (\cref{fig:bcd}, right).

Details on the boundary conforming Delaunay mesh generation can be
found in~\cite{SiGarFuh10} and the references therein.

%--- Finite Volumes ----------------------------------------
\section{Voronoi finite volumes and their essential properties}

If applying the Gauss theorem to a Voronoi cell or its intersection
with $\partial \Omega$, each facet is related to an edge and therefore
is a set of simplices.
It is more convenient to assemble the finite volume equations using
simplices of the dual Delaunay mesh, in particular it is more
convenient to define material properties per simplex than per Voronoi
volume.

\subsection{Voronoi Finite Volume Equations}

%--- Differences along edges ----------------------------------------
The \emph{adjacency matrices} $\tilde G_i$ of $i$-dimensional
simplices mapping vertices to edges are defined recursively through
the relation
\[
   \tilde{G}_{i+1} =
   \begin{pmatrix*}[l]
     \phantom{-}\tilde{G}_i & \bm{0}_{\frac{i(i-1)}{2}} \\
      -\mathbb{I}_{i+1}       & \bm{1}_{i+1}
   \end{pmatrix*},
   \quad i = 1, 2, \dotsc,
\]
with $\mathbb{I}_j$ denoting the $j \times j$ identity matrix,
$\bm{0}_j = ( \overbrace{ 0,\dotsc, 0}^{j}{)}^\mathrm{T}$,
$\bm{1}_j = ( \overbrace{1, \dotsc, 1}^{j}{)}^\mathrm{T}$,
and $\tilde{G}_1 = (-1, 1)$ describing the \emph{difference along an interval},
\[
   \tilde{G}_1 \begin{pmatrix} \bm{x}_1\\\bm{x}_2 \end{pmatrix}
   =
      \begin{pmatrix} -1 & 1 \end{pmatrix}
      \begin{pmatrix} \bm{v}_1\\\bm{v}_2 \end{pmatrix}
      = \bm{v}_2 - \bm{v}_1
      .
\]
It holds,
\begin{align*}
   \tilde{G}_2 &= 
   \begin{pmatrix*}[r]
      \tilde{G}_1     & \phantom{-}\bm{0}_1 \\
      -\mathbb{I}_2  &            \bm{1}_2
   \end{pmatrix*}
   =
   \left(\begin{array}{rr|r}
      -1  &  1  & \phantom{-}0 \\
      \hline{}
      -1  &  0  &            1 \\
       0  & -1  &            1 \\
   \end{array}\right),
   \\
   \tilde{G}_3 &=
   \begin{pmatrix*}[r]
      \tilde{G}_2     & \phantom{-}\bm{0}_3 \\
      -\mathbb{I}_3  &            \bm{1}_3
   \end{pmatrix*}
   =
   \left(\begin{array}{rrr|r}
     -1  &  1  &  0  &  \phantom{-}0\\
     -1  &  0  &  1  &             0\\
      0  & -1  &  1  &             0\\
      \hline{}
     -1  &  0  &  0  &  1\\
      0  & -1  &  0  &  1\\
      0  &  0  & -1  &  1\\
   \end{array}\right),
   \\
   \tilde{G}_4 &=
   \begin{pmatrix*}[r]
      \tilde{G}_3     & \phantom{-}\bm{0}_6 \\
      -\mathbb{I}_4  &            \bm{1}_4
   \end{pmatrix*},
   \quad \dotso{}
\end{align*}

$\tilde G_d$ maps vertices to edges.
Correspondingly, $\tilde G_d^\mathrm{T}$ maps edges to vertices and
$\tilde G_d^\mathrm{T} \tilde G_d$ is the \emph{graph-Laplacian}
of the $d$-dimensional simplex.

Denote the length of an edge $\bm{E}_i^1$ by $\lvert
\tilde G_1 \bm{E}{_i^1} \rvert$ and the diagonal matrix
of all edge lengths for a simplex $\bm{E}^d_i$ by $[ \,
\lvert \tilde G_d \bm{E}^d_i \rvert \, ]$.
Then,
\[
   {\left[ \lvert \tilde G_d \bm{E}^d_i \rvert \right]}^{-1}
      \tilde G_d \bm{u}
\]
is the approximation of $\nabla u$ in directions $\partial u/\partial
n_{k}$ orthogonal to the surfaces of the Voronoi cell related to
the simplex $\bm{E}^d_i$, transforming linear functions defined
on each edge into the correct constant.
Further,
\[
   \tilde G_d^T  [ \sigma_{\bm{E}_i^d} ]
      = \tilde G_d^\mathrm{T} {\left[ \lvert \tilde G_d \bm{E}^d_i \rvert \right]}^{-1} 
      {\left\{
         \left[ \lvert \tilde G_d \bm{E}^d_i \rvert \right] 
         \left[ \sigma_{\bm{E}^d_i} \right]
      \right\}}^{\frac{1}{2}}
   \]
is the symmetric form of the edge to vertex map, including a geometric weight 
of the Voronoi facet intersected by the edge times the edge length.
Hence, the metric factor $\{\cdot\}$ has the meaning of the determinant times two 
(the edge length is twice the hight of the $d$-dimensional simplex formed by 
the intersection of the simplex and $\sigma_i$) and not that of the volume
intersected with the corresponding part of each edge of the simplex
$\bm{E}^d_i$.
It has to be handled in all terms of the operator, boundary conditions consistently.
Hence,
\[
   \varepsilon \int_{\bm{E}^d_i} \nabla \cdot \nabla u \, \D{V}
   = \oint_\Sigma \nabla u \, \D{\bm{S}}
   \approx \varepsilon \, \tilde G_d^\mathrm{T} [\sigma_i] \,
   {\left[ \lvert \tilde G_d \bm{E}^d_i \rvert \right]}^{-1} \,
   \tilde G_d \bm{u}  := \varepsilon \,
   G_d^\mathrm{T} G_d \bm{u}
\]
is the approximation of the \emph{volume integrated Laplace operator}
on $\bm{E}^d_i$, where $\varepsilon = const$ on $\bm{E}^d_i$
and $\sigma_j \ge 0$ for each edge $j$.

The boundary condition
\[
   \alpha(u) u +\beta(u) \frac{\partial u}{\partial n} +\gamma(u)=0
\]
is integrated over Voronoi surfaces belonging to surface vertices and the surface
is defining the normal direction and coincides with the edge direction in case of a boundary conforming mesh,
\[
   \varepsilon \frac{\partial u}{\partial n} \partial V_{s,i}
   = - \partial V_{s,i} \frac{\varepsilon}{\beta(u)}
   \left( u_i \alpha(u_i) + \gamma(u_i) \right)
   .
\]
The normal derivative causes a change of the corresponding diagonal
element and the right-hand side. Dirichlet values follow for $\beta
\approx \epsilon_\mathrm{machine}^2$ and are easily homogenized.

A boundary conforming Delaunay mesh allows to split the Voronoi cells along the boundary, which coincides with the simplices, and to introduce the same surface integration on the $d-1$-dimensional Voronoi cells.
In case of boundary layers, the vertices on the boundary can assume the extrema enforced, e.g., by the boundary condition.
This is not the case if the Voronoi cells are used to approximate $\Omega_i^d$.

The essential assumption for the FVM is the validity of the Gauss theorem.
The largest function space where the Gauss theorem hods is the space of functions with bounded variation.
It is a very weak assumption.

Allowing an arbitrary jump in the diffusion coefficient of adjacent regions results in additional angle restrictions.

%--- Compensation effect ------------------------------
\subsection{Compensation and the S-matrix property}

A major advantage of the Voronoi FVM is that the system
matrix is a non-singular \emph{M-matrix}.  An \emph{M-matrix}
(Minkowski matrix) is a Z-matrix (off-diagonal entries are non-positive)
with a positive diagonal whose column sums are all non-negative.
If at least one column sum of an M-matrix is positive, then it is
inverse-positive, i.e.,\ it is invertible and its inverse is a
positive matrix.  A symmetric M-matrix is called an \emph{S-matrix}
(Stieltjes Matrix).

\begin{figure}[t]%
   \centering{}%
   \includegraphics[width=0.3\textwidth,clip]{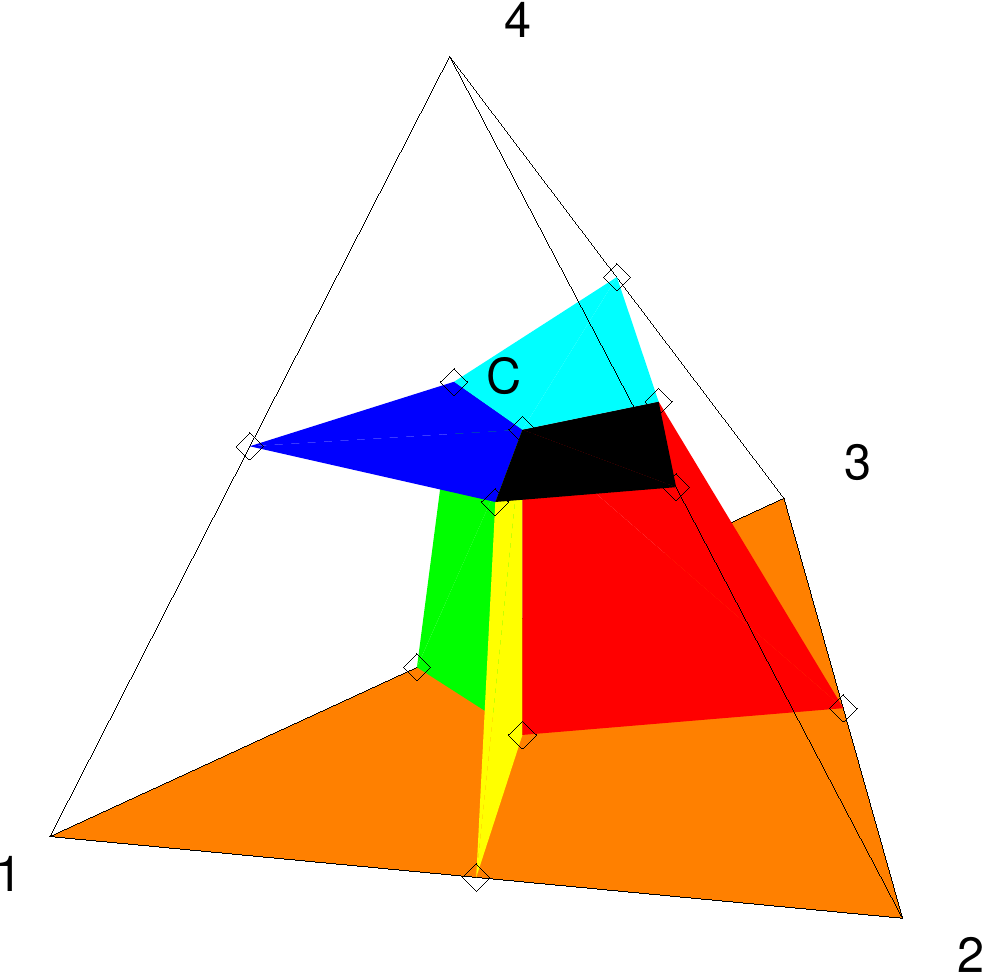}%
   \hspace{0.1\linewidth}%
   \includegraphics[width=0.3\textwidth,clip]{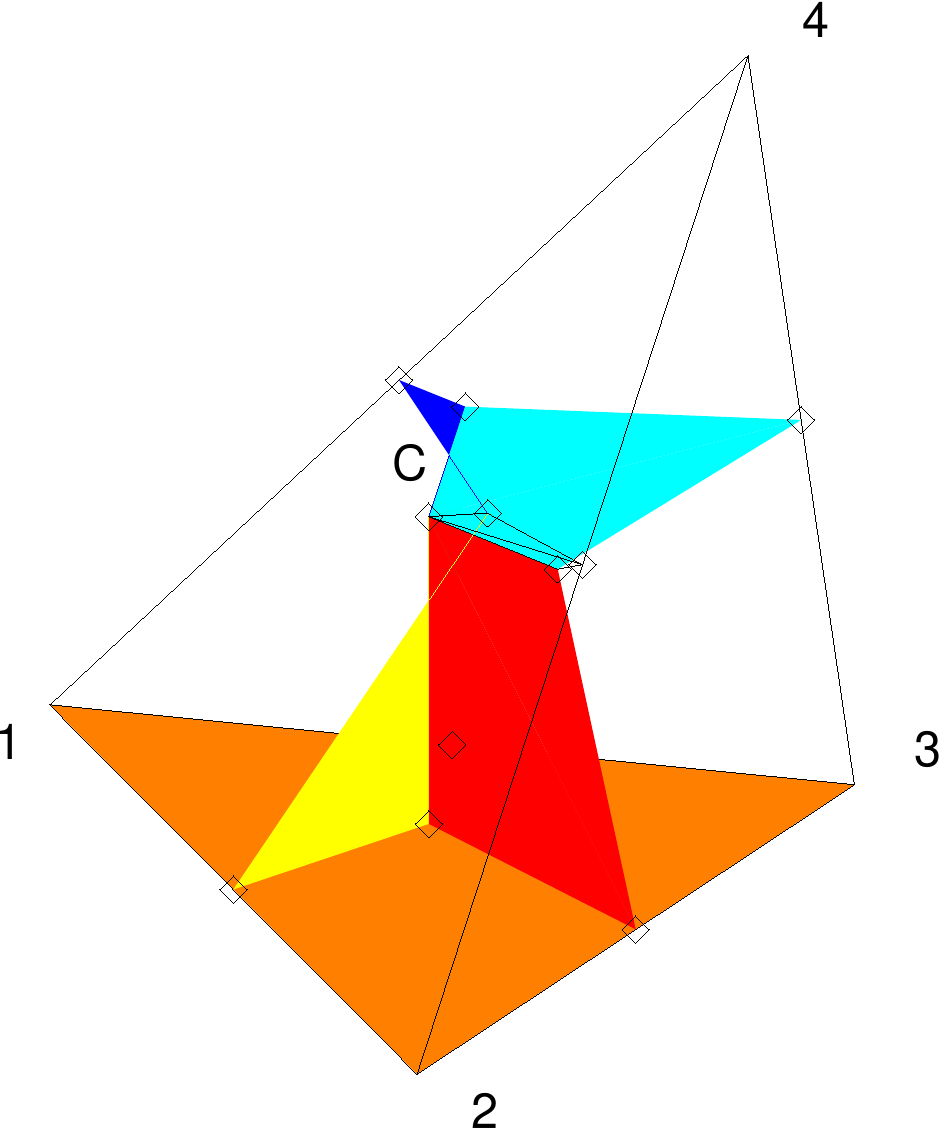}%
   \hfill{}%
   \caption{Left: all circumcenters are inside the corresponding defining
      simplices of dimension $d = 1, 2, 3$.
      Right: the circumcenter $C$ of the tetrahedron and the circumcenter of
      $\triangle 234$ are outside of their defining simplices%
      \label{fig:tet:1}%
      \label{fig:tet:2}%
   }%
\end{figure}

Consider the $\mathbb{R}^3$ case and two possible Delaunay tetrahedra
in \cref{fig:tet:1,fig:tet:2}.  The free fourth vertex opposing
each triangular facet determines the four Delaunay conditions.  In
\cref{fig:tet:1} (left) all circumcenters (hereafter $C$ with an index
describing the corresponding, possibly low-dimensional simplex) are
inside the corresponding defining simplex of dimension one, two,
or three, respectively.

In \cref{fig:tet:2} (right) the tetrahedron's circumcenter $C_{1234}$ and
the circumcenter $C_{234}$ of $\triangle 234$ are outside of their
defining simplex.
Each of the Voronoi faces is defined by two triangles, e.g.,\ the
one defined by the vertices 1 and 2 consists of $\triangle C_{12}
C_{123} C_{1234}$ and $\triangle C_{21} C_{124} C_{1234}$ which are
orthogonal to the edge $12$ and change the sign accordingly to the
positions of the related triangle's and tetrahedron's circumcenters
and have a positive projection on edges 12 and 21 for interior circumcenters.
The neighboring tetrahedron sharing $\triangle 124$ as a facet has
the same circumcenters of the shared edges and the same direction
of the normal to $\triangle 124$.
The Voronoi facet is bounded by a polygon connecting all circumcenters
of tetrahedra sharing one edge and is located in the plane
orthogonal to that edge.
Summing up all signed contributions of tetrahedra sharing an edge
yields the signed Voronoi facet area related to this edge.
The Delaunay condition guaranties that the distance between the
circumcenters of two neighboring tetrahedra is non-negative and,
thus, the signed Voronoi facet area is non-negative.

This effect is called \emph{compensation}: the negative contribution
is compensated by the positive contribution of the neighboring
tetrahedron.  This relation holds in any dimension and guarantees
that the finite volume matrix is an S-matrix or an M-matrix.

If more than $m \ge d + 2$ vertices lie on a common $d$-sphere, $m
- d$ circumcenters will coincide.  A boundary conforming Delaunay
mesh guarantees that all circumcenters are inside the corresponding
defining simplices, the per simplex defined oriented parts of Voronoi
cells and surfaces are non-negative, and the volume of
$\Omega_i$ as well as the surface measures are preserved.
This allows an easy handling of all third kind boundary conditions.

For the P\textsubscript{1} (linear) finite elements method (FEM),
the Laplace matrix in $d$ dimensions is the projection of each of
the $d - 1$-dimensional facets on all others,
\[
   P^\mathrm{-T} P^{-1},
   \quad P = \bm{v}_j - \bm{v}_0,
   \quad j = 1, \dotsc, d
   .
\]
In two dimensions, it is identical to the Voronoi finite volume
matrix $G_d^\mathrm{T} G_d$.
In more than two dimensions, it is not the case anymore (except in
very special cases), the compensation of negative contributions
works differently, and boundary conforming Delaunay meshes are
neither sufficient nor necessary for the S-matrix property for the
FEM (for details see~\cite{BanRos87,Ker96,XuZik99,KoFHPS00}
and the references therein).
The so-called \emph{non-obtuse angle condition} (all angles between
the surface normals are non-obtuse), which guarantees the S-matrix
property for the FE Laplace matrix, is too restrictive and difficult
to fulfill in practice.

%--- Dissipativity ----------------------------------------
\subsection{Maximum principle and dissipativity}

The finite volume solution satisfies the \emph{weak discrete maximum principle}.
Let $\bm{u}$ fulfill $G_d^d G_d \bm{u}= \bm{0}$
with some Dirichlet boundary values $\bm{u}_{D}$, $\hat u=\max u_{D}$,
$\check u=\min u_{D}$, and assume that $\max u = U > \hat u$.
Testing the equation with $(\bm{u}-\hat u)^\mathrm{+T}$ results in
$(\bm{u}-\hat u)^\mathrm{+T} G^\mathrm{T}_d G_d \bm{u} > 0$ 
because $G_d \bm{u}$ and $G_d {(\bm{u}-\hat u)}^{+}$
have either the same signs (and at least one edge value is non zero)
or $G_n {(\bm{u}-\hat u)}^{+} = \bm{0}$, hence a contradiction.

Dissipativity follows from the same argument. In the discrete case both are
equivalent in this sense.

Dissipativity of zero-order terms:
the Voronoi volumes form a diagonal
matrix $[V]$, larger support can not be handled without smallness
assumptions of the variation of the solution on the support.

%--- Non-Delaunay ----------------------------------------
\subsection{The effect of non-Delaunay meshes for the Laplace operator}
\label{sec:nonDelaunayLaplace}

Since Voronoi cells have non-negative surface measures, assembling
\[
   A = G_d^\mathrm{T} [\varepsilon] G_d
\]
results in an S-matrix: $A_{ii} > 0$, $A_{ij} \le 0$, $i\ne j$, and ${\bm 1}^\mathrm{T}_d G_d^\mathrm{T} = \bm{0}_d$ (column and row sum zero at vertices without boundary conditions).
Clearly, $A^{-1}$ is positive for irreducible $A$ or connected meshes.
Bounded $\varepsilon$ and boundary conditions add only positive values to the
diagonal at least at one vertex.
Hence, the solution is bounded by the Dirichlet values and positive for zero
Dirichlet values and positive right-hand sides.

If the underlying mesh is non-Delaunay, summing up the surface contributions per edge
may result in negative surfaces around this edge. It could be an arbitrary edge
between two vertices not corresponding to a Voronoi face (the circumsphere
of one simplex is containing the circumcenter of another).
It follows $A_{ij}>0$ and positivity is not guarantied anymore, at least one may
observe some local extrema at vertices $i$ and $j$ for significant deviations from
the Delaunay condition. Due to the construction,
$G_d \bm{1}_d =\bm{0}_d$. Setting $A_{ij}=0$ consistently with $A_\mathrm{ii}$
corresponds to introducing a local Neumann boundary condition and an inconsistent discretization.

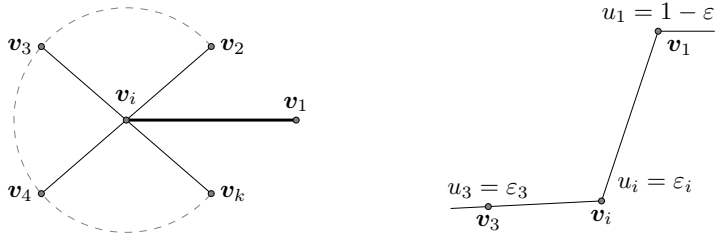
\begin{figure}[t]%
   \centering{}%
   \begin{tikzpicture}[scale = 0.8]%
      % define vertices 1, 2, 3, 4
      \tkzDefPoint( 0.0,  0.0){pi}
      \tkzDefPoint( 3.0,  0.0){p1}
      \tkzDefPoint( 1.5,  1.3){p2}
      \tkzDefPoint(-1.5,  1.3){p3}
      \tkzDefPoint(-1.5, -1.3){p4}
      \tkzDefPoint( 1.5, -1.3){pk}
      % triangles 1-2-4 (dotted)
      \tkzDrawArc[dashed](pi,p2)(pk)
      \tkzDrawSegment[very thick](pi,p1)
      \tkzDrawSegments[thin](pi,p2 pi,p3 pi,p4 pi,pk)
      % draw vertices
      \tkzDrawPoints(pi,p1,p2,p3,p4,pk)
      \tkzLabelPoint[above = 1mm](pi){$\bm{v}_i$}
      \tkzLabelPoint[above](p1){$\bm{v}_1$}
      \tkzLabelPoint[right](p2){$\bm{v}_2$}
      \tkzLabelPoint[left](p3){$\bm{v}_3$}
      \tkzLabelPoint[left](p4){$\bm{v}_4$}
      \tkzLabelPoint[right](pk){$\bm{v}_{k}$}
   \end{tikzpicture}%
   \hspace{0.1\linewidth}%
   \begin{tikzpicture}[scale = 0.8]%
      % define vertices 1, 2, 3, 4
      \tkzDefPoint( 0.0,  0.0){pi}
      \tkzDefPoint( 1.0,  3.0){p1}
      \tkzDefPoint( 2.0,  3.0){p11}
      \tkzDefPoint(-2.0, -0.1){p3}
      % triangles 1-2-4 (dotted)
      \tkzDrawLine[thin, add = 1/3 and 0](p3,pi)
      \tkzDrawSegments[thin](pi,p1 p1,p11)
      % draw vertices
      \tkzDrawPoints(pi,p1,p3)
      \tkzLabelPoint[below](pi){$\bm{v}_i$}
      \tkzLabelPoint[above right = 0 and 1mm](pi){$u_i=\varepsilon_i$}
      \tkzLabelPoint[below right](p1){$\bm{v}_1$}
      \tkzLabelPoint[above](p1){$u_1 = 1 - \varepsilon$}
      \tkzLabelPoint[below](p3){$\bm{v}_3$}
      \tkzLabelPoint[above](p3){$u_3 = \varepsilon_3$}
   \end{tikzpicture}%
   \caption{%
      A jump in the solution along the edge $\bm{v}_i \bm{v}_1$
      from $u(\bm{v}_i) = \varepsilon_i$
      to $u(\bm{v}_1) = 1 - \varepsilon$%
   }\label{fig:jump}
\end{figure}

Assume that a solution $u$ has a moving internal layer with
\[
   u_1 = 1 - \varepsilon,
   \quad u_j = \varepsilon_j,
   \quad \text{with}~ 0< \varepsilon_j \le \varepsilon \ll 1,
   \quad j = i, 2, \dotsc, k,
\]
representing the order of the variation of $u$ outside of the jump (\cref{fig:jump}).
Let $A_\text{non-Del}$ be of type $\tilde G^\mathrm{T} [W] \tilde G$ and the sum over all simplices 
at vertex $i$ (the sum is not indicated, the summation order over edges or
simplices of given type can be chosen for each purpose and the context 
should distinguish between a local and a global meaning).

Let $l_{ij}$ be the edge length for all neighbors
$\bm{v}_j$ of $\bm{v}_i$ and
\[
   \begin{cases}
      W_{ij}   = \frac{\sigma_{ij}}{l_{ij}} = 1,
      & j= 2, \dotsc, k,\\
      W_{i1} = \frac{\sigma_{i1}}{l_{i1}} = -1,
   \end{cases}
\]
i.e.,\ the edge $\bm{v}_i \bm{v}_1$ is non-Delaunay.
Hence, $A_\text{non-Del} {\bm u}|_i$ at vertex $i$ results in
\[
   A_\text{non-Del} {\bm u}|_i
   = \sum^k_{j=2} (u_i - u_j)
   -(u_i - u_1)
   = \sum^k_{j=2} (\varepsilon_i- \varepsilon_j)
   -(\varepsilon_i - (1 - \varepsilon))
  \approx 1 > 0
  ,
\]
whereas in the Delaunay case
\[
   A_\text{Del} {\bm u}|_i
   = \sum^k_{j=2} (u_i - u_j)
      + (u_i - u_1)
   = \sum^k_{j=2} (\varepsilon_i- \varepsilon_j)
      + (\varepsilon_i - (1 - \varepsilon))
   \approx -1 < 0
   .
\]
Hence, the right-hand side at $i$ is very large compared to that
at vertices inside the plateau of level $\varepsilon$ or $1-\varepsilon$.
Because $A \bm{u}$ in Newton's method (or the corresponding
linear residual correction method) is a part of the function $f$
and $A_\text{Del}$ contributes to the Jacobian ($A$) we have
\[
   \delta \tilde u\vert_i
      = (-A_\text{Del}^{-1} f)\vert_i > 0
\]
and therefore $u^+ = u + \delta u$, the next approximation of
$u_i$ in the Newton process, is increased towards the mean
value of its next neighbors by the diffusion part of the Jacobian.

The factor $-1$ and the huge (by modulus) right-hand side is
equivalent to a large source of the wrong sign in the non-Delaunay
case.  Additionally, $A_\text{non-Del}$ can be indefinite.  Clearly,
a constant $A$ cannot be the reason for the jump but other terms may
create jumps in case of a small diffusion.  A small $A_\text{non-Del}$
can be sufficient to violate the bounds for $u$ and, thus, causes
serious problems which can not be cured by simple limiters for $u$
without creating new issues depending on the details of the problem.

In a more general setting, this example corresponds to a locally
negative dissipation rate and therefore to the \emph{local free energy
production} by the discretization scheme.
The concept of a decaying free energy along trajectories will be
violated in general and the possible steady state will deviate in
dependence of the mesh from the minimal energy solution because
the free energy produced by the numerical scheme has to be dissipated
somewhere else.

%--- Mass Matrix -------------------------------------------------
\subsection{Continuous dependence on vertex positions}

Consider the following one-dimensional boundary value problem:
\begin{align*}
   \begin{cases}
      u(y)'' + c u(y) = 0, & y \in (0,1), \quad c = const,\\
      u(0) = 1,\\
      u(1) = 0.
   \end{cases}
\end{align*}
This problem and its discretization on a given one-dimensional mesh
in $y$ can be extended orthogonally in $x$ direction with arbitrary
step sizes and homogeneous Neumann boundary conditions while keeping
the profile of $u$ in $y$ direction (the two-dimensional solution
will be invariant in $x$ direction).

\begin{figure}[t]
   \centering{}
   \begin{tikzpicture}[scale = 1.5]%
      % define vertices 1, 2, 3, 4
      \tkzDefPoint(0, 0){p00}
      \tkzDefPoint(0, 1){p01}
      \tkzDefPoint(0, 2){p02}
      \tkzDefPoint(1, 0){p10}
      \tkzDefPoint(1, 1){p11}
      \tkzDefPoint(1, 2){p12}
      \tkzDefPoint(2, 0){p20}
      \tkzDefPoint(2, 1){p21}
      \tkzDefPoint(2, 2){p22}
      \tkzDefPoint(3, 0){p30}
      \tkzDefPoint(3, 1){p31}
      \tkzDefPoint(3, 2){p32}
      \tkzDefPoint(4, 0){p40}
      \tkzDefPoint(4, 1){p41}
      \tkzDefPoint(4, 2){p42}
      % areas
      \tkzFillPolygon[gray!15](p00,p01,p12,p22,p21,p10)
      \tkzFillPolygon[gray!15](p20,p21,p32,p42,p41,p30)
      % circumball example
      \tkzDefMidPoint(p01,p12)\tkzGetPoint{c}
      \tkzDrawCircle[dotted](c,p01)
      % box + diagonals
      \tkzDrawSegments[thin](p00,p40 p40,p42 p42,p02 p02,p00) % box
      \tkzDrawSegments[thin](p01,p41 p10,p12 p20,p22 p30,p32) % grid
      \tkzDrawSegments[thin](p00,p11 p10,p21 p20,p31 p30,p41) % diagonals
      \tkzDrawSegments[thin](p01,p12 p11,p22 p21,p32 p31,p42) % diagonals
      \tkzDrawSegments[dashed, thin](p02,p11) % diagonals
      %% draw vertices
      \tkzDrawPoints(p11,p31)
      \tkzLabelPoint[below right = 0 and 0.15](p11){$\bm{v}_1$}
      \tkzLabelPoint[below right](p31){$\bm{v}_2$}
   \end{tikzpicture}%
   \hspace{0.1\linewidth}%
   \begin{tikzpicture}[scale = 1.5]%
      % define vertices 1, 2, 3, 4
      \tkzDefPoint(0, 0){p00}
      \tkzDefPoint(0, 1){p01}
      \tkzDefPoint(0, 2){p02}
      \tkzDefPoint(1, 0){p10}
      \tkzDefPoint(1, 1){p11}
      \tkzDefPoint(1, 2){p12}
      \tkzDefPoint(2, 0){p20}
      \tkzDefPoint(2, 1){p21}
      \tkzDefPoint(2, 2){p22}
      \tkzDefPoint(3, 0){p30}
      \tkzDefPoint(3, 1){p31}
      \tkzDefPoint(3, 2){p32}
      \tkzDefPoint(4, 0){p40}
      \tkzDefPoint(4, 1){p41}
      \tkzDefPoint(4, 2){p42}
      % areas
      \tkzFillPolygon[gray!15](p00,p20,p22,p02)
      \tkzFillPolygon[gray!15](p21,p32,p41,p30)
      % box + diagonals
      \tkzDrawSegments[thin](p00,p40 p40,p42 p42,p02 p02,p00) % box
      \tkzDrawSegments[thin](p01,p41 p10,p12 p20,p22 p30,p32) % grid
      \tkzDrawSegments[thin](p00,p11 p11,p20 p21,p30 p30,p41) % diagonals
      \tkzDrawSegments[thin](p02,p11 p11,p22 p21,p32 p32,p41) % diagonals
      %% draw vertices
      \tkzDrawPoints(p11,p31)
      \tkzLabelPoint[below right = 0 and 0.15](p11){$\bm{v}_1$}
      \tkzLabelPoint[below right](p31){$\bm{v}_2$}
   \end{tikzpicture}%
   \caption{%
      The choice of the Delaunay diagonal edges for the regular rectangular
      mesh is not unique since the four vertices lie on one circle.
      `Lower left to upper right' edge assignment rule for a rectangular vertex
      distribution results in a translation-invariant mesh (left).
      Some of the diagonal edges are chosen differently (right),
      resulting in the so-called criss-cross pattern around $\bm{v}_1$%
   }\label{fig:crisscross}%
\end{figure}
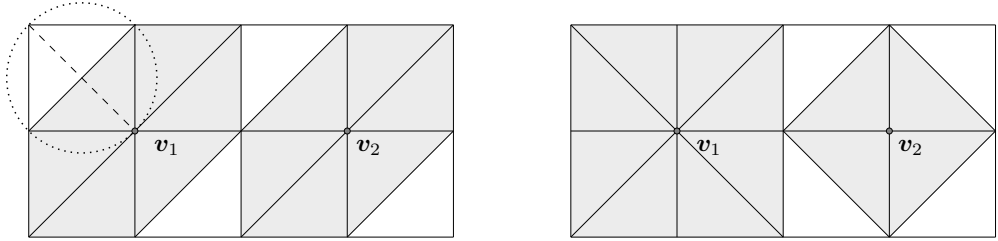

The Delaunay triangulation for the resulting regular rectangular
vertex set is not unique because each four neighboring vertices
sharing one interval in $x$ and $y$ directions lie on one circle
(\cref{fig:crisscross}, left: both diagonal edge directions are
possible).  Therefore, in this case, the unique assignment of the
diagonal edges can be done only by a special rule, e.g., `lower
left to upper right'.  Using exact arithmetics to improve the
Delaunay test will not help, since it will still result in
non-uniqueness in case of exactly represented coordinates or randomly
assign an edge depending on the rounding error in coordinates if
the four vertices are co-circular by construction but have rounded
coordinates in the floating point representation.  For instance,
rotating meshes in \cref{fig:crisscross} by $\pi/4$ will result in
meshes where four co-circular vertices will not be co-circular
anymore in exact arithmetics, because $\sin \pi/4 = \cos \pi/4 =
\sqrt{2}/2$ cannot be represented exactly by floating point numbers.

For the FVM, the choice of the diagonal edges is not
important because the Voronoi facets corresponding to the diagonal
edges have a measure of zero.  The Voronoi surfaces and the Voronoi
volumes are scaled correspondingly to the step size in $x$ direction,
thus, preserving the invariance of the solution in $x$ direction.
Moreover, small perturbations of the vertex coordinates will result
in a comparably small change of facet measures (recall that Voronoi
cell and facet measures are continuous functions of vertex positions).
The same is true for the matrices of type $G^\mathrm{T} G$ in any
dimension $d$ (e.g., the Laplace and the mass parts).  Hence, a
very small change of vertex positions will result in a very small
change of the numerical solution.

For the FEM, the situation is different
%.  The fact, that the element volumes scale correspondingly to the step
%size in $x$ direction is not sufficient
since topology comes
into play as well.  On one hand, for meshes in \cref{fig:crisscross},
the Laplace part of the P\textsubscript{1}
finite element discretization is invariant to the choice of the
diagonal edges, since the finite element Laplace part is equal to
that of the finite volumes for $d=2$ and the latter depends only
on the vertex positions.  On the other hand, the mass matrix $M$
depends on topology, e.g.,
\[
   M_{ii} 
   = \frac{2 \lvert \omega_{i} \rvert}{(d+1)(d+2)}
   ,
\]
where $\omega_{i}$ is the patch associated with the vertex
$\bm{v}_i$ and $\lvert \omega_{i} \rvert$ is its
volume (area).  In \cref{fig:crisscross} (left), the mesh is
translation-invariant (up to the first and the last vertex)
and each patch has the same size.
Hence,
\[
   M_{11} = M_{22} 
   = \frac{2 \times \left(6 \times \frac{h^2}{2}\right)}{12}
   = \frac{1}{2} h^2
   .
\]
In \cref{fig:crisscross}~(right), the mesh is still Delaunay but
some of the diagonal edges are chosen differently.  The patches of
$\bm{v}_1$ and $\bm{v}_2$ are of different sizes and so are the
corresponding entries of the mass matrix,
\[
   M_{11}
   = \frac{2 \times \left(8 \times \frac{h^2}{2}\right)}{12}
   = \frac{2}{3} h^2,
   \quad \text{whereas} \quad
   M_{22}
   = \frac{2 \times \left(4 \times \frac{h^2}{2}\right)}{12}
   = \frac{1}{3} h^2
   ,
\]
so that the mass part is not translation-invariant, while the Laplace
part is.  As a consequence, the numerical solution will loose its
invariance in the $x$ direction.  Moreover, the inconsistency of
the mass matrix may also lead to lower than expected order of
convergence.
Lumped-mass quadrature cannot correct this issue in general because
the mass assignment to vertices depends on the mesh topology and,
thus, is already inconsistent before mass lumping is carried out.
An interested reader can find a discussion on the lost of accuracy
of linear and higher-order finite elements in dependence of
the mesh topology in recent work of Kopteva~\cite{Kop14,Kop20}.

%---------------------------------------------------------------------
\section{Application examples}
%---------------------------------------------------------------------

%--- Image segmentation ----------------------------------------------
\subsection{Phase separation}

Consider diffusing black and white particles which sum
up to 1 everywhere (Fermi statistics), a short range interaction potential
(black attracts black), and the diffusion acting as repulsive force
(that is the qualitative meaning of the well-known Einstein relation).
Gradients can be increased by increasing the attractive interaction potential
relative to the diffusion.
The concentrations $u_{b}$ of
the black and $u_{w}$ of the white particles fulfill
\[
   1 > 1 - \varepsilon > u_{b} > \varepsilon > 0,
   \quad u_{w} = 1 - u_{b}
   .
\]
If the diffusive repulsion is dominating, then the solution is a perfect gray distribution.
Otherwise, the black and the white particles are separated by a minimal surface.

In the following we give a brief mathematical description for this non-local phase segregation problem;
a detailed discussion can be found in~\cite{GajGar05a} and the references therein.
The problem is given by
\begin{align}
   \begin{cases}
      \frac{\partial u}{\partial t}
         - \nabla \cdot \left(
            f (\abs{\nabla v})
            \left( \nabla u + \frac {\nabla w}{\phi''(u)} \right) \right) = 0
      \quad \text{in $\Omega$}
      ,
      \\
      u(0,\cdot) = u_0(\cdot)
      ,
      \\
      v = \phi'(u) + w
      ,
      \\
      w(t,\bm{x}) = \int_\Omega \mathcal{K} \left( \abs{\bm{x}-\bm{y}} \right) 
         \left(1 - 2 u(t,\bm{y}) \right) \D{\bm{y}}
      .
   \end{cases}
   \label{eq:phaseSep}
\end{align}
where $\Omega \subset \R^N$, $1 \le N \le 3$ is a bounded Lipschitz domain, $\phi$ a convex function, the kernel $\mathcal{K}$ represents non-local attracting forces, and $w$ and $v$ are the interaction and chemical potentials.
The initial value $u_0$ satisfies $0< u_0 <1$ and the system is completed with homogeneous Neumann boundary conditions.

The system was rigorously derived in~\cite{GiaLeb97} for the case of a constant $f$
and can be seen as a non-local variant of the Cahn-Hilliard equation
associated with the local Ginzburg-Landau free energy
\[
   F_{gl}(u) = \int_\Omega \left(
      \phi(u) + \kappa  u \left(1 - u\right)
      + \frac{\lambda}{2} \abs{\nabla u}^2
   \right) \D{\bm{x}}
   .
\]
The Lyapunov functional for \cref{eq:phaseSep} turns out to be
\[
   F(u) = \int_\Omega \left(
      \phi(u) + u \int_\Omega \mathcal{K} 
   \left(\abs{x-y}\right) \left(1 - u(y) \right) \D{y} \right) \D{\bm{x}}
   .
\]
We consider the specialized model with a constant $f$,
$\phi(u) =  u \ln u + (1 - u) \ln (1 - u)$,
and $\mathcal{K}$ being the Green's function of the elliptic boundary value problem
\begin{align*}
   \begin{cases}
      - \sigma^2 \nabla \cdot \nabla w + w = m \left(1 - 2u\right)
      & \text{in $\Omega$}
      ,
      \\
      \nu \cdot \nabla w = 0
      & \text{on  $\partial \Omega$}
      .
   \end{cases}
\end{align*}
with constant $\sigma, m > 0$.

In the short time regime the model has been used for image reconstruction and
denoising as well and it is not restricted to just two species~\cite{GajGar05,GajGri06,Gri04}.

An appropriate numerical scheme for this problem is
\begin{align*}
   \left(G^T \sigma^2 G +[V]\right) \bm{w}^+ 
      &= [V] m \left(\bm{1} -2 \bm{u}^+ \right)
      ,
      \\
   [V] \frac{\bm{u} - \bm{u}^+} \tau
      & = G^T f \left(
         G \bm{u}^+ + \left[ \frac{1}{ \Phi''_a} \right]
         \frac{1}{2} G \left( \bm{w}^+ + \bm{w} \right) \right)
      ,
\end{align*}
with the notation $z := z(t)$, $z^+ := z(t+\tau)$, $[V]$ being the diagonal matrix of Voronoi volumes, and $\left[\frac{1}{ \Phi''_a}\right]$ depending on ${\bm{u}^+}$~\cite{GajGar05a}.
This scheme is shown to be dissipative~\cite[Theorem~2.7]{GajGar05a}.
Moreover, if $\bm{v}^+ \ne c$, $0 < u < 1$, then the implicit Euler scheme
preserves the property $0 < u^+ < 1$~\cite[Proposition~2.9]{GajGar05a}

The time step size control in the computation is based on the
predictor corrector differences of the free energy and discrete
dissipation rate.  Jacobian matrices are derived as analytic
derivatives on simplices, and Newton's method is applied without
any damping.
The initial value is
\[
   u_0(\bm{x}) = \varepsilon + \frac{x_1}{128 \left( 1 + 2\varepsilon \right)}
  ,
  \quad \Omega = (0,128) \times (0,128),
  \quad m = 8,
  \quad \varepsilon=10^{-7},
  \quad \sigma = \sqrt{f} = 2
  .
\]

\begin{figure}[t]%
   \includegraphics[width=0.3\textwidth]{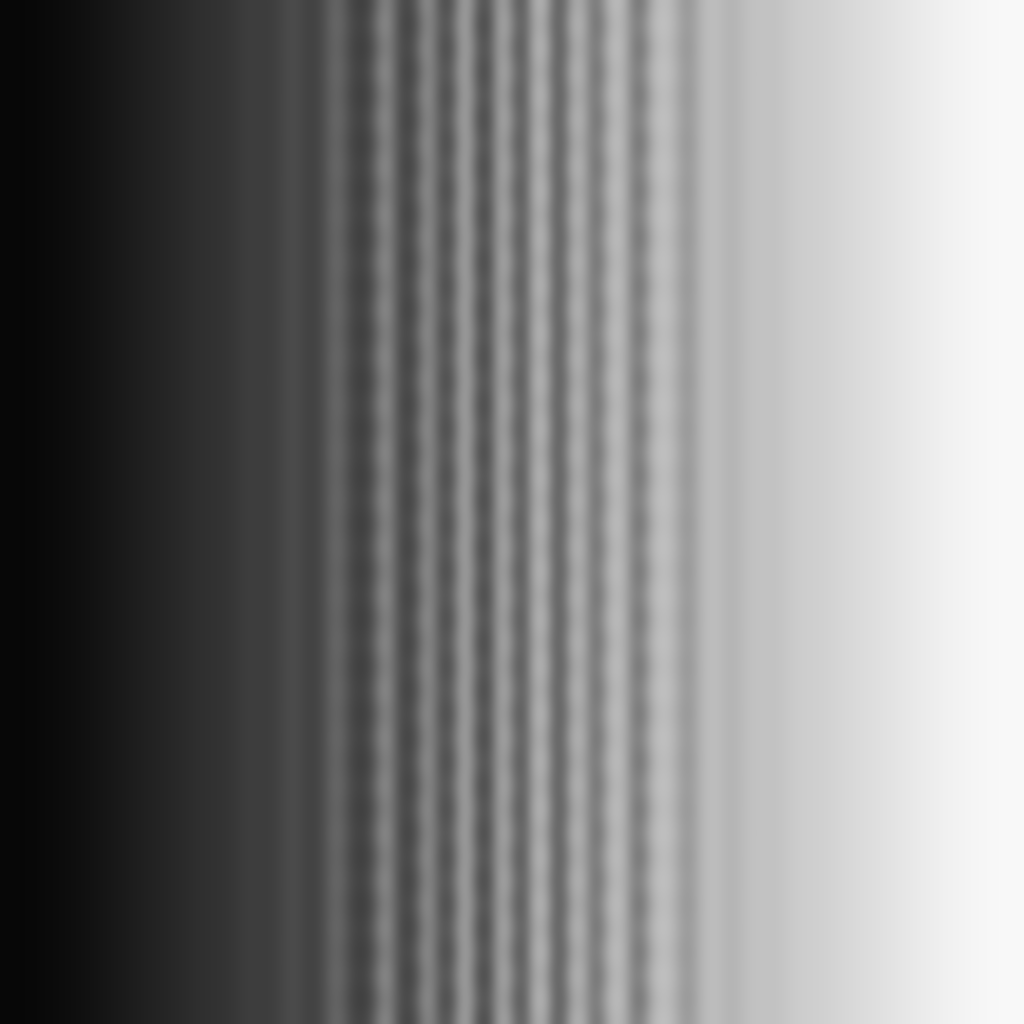}%
   \hfill{}%
   \includegraphics[width=0.3\textwidth]{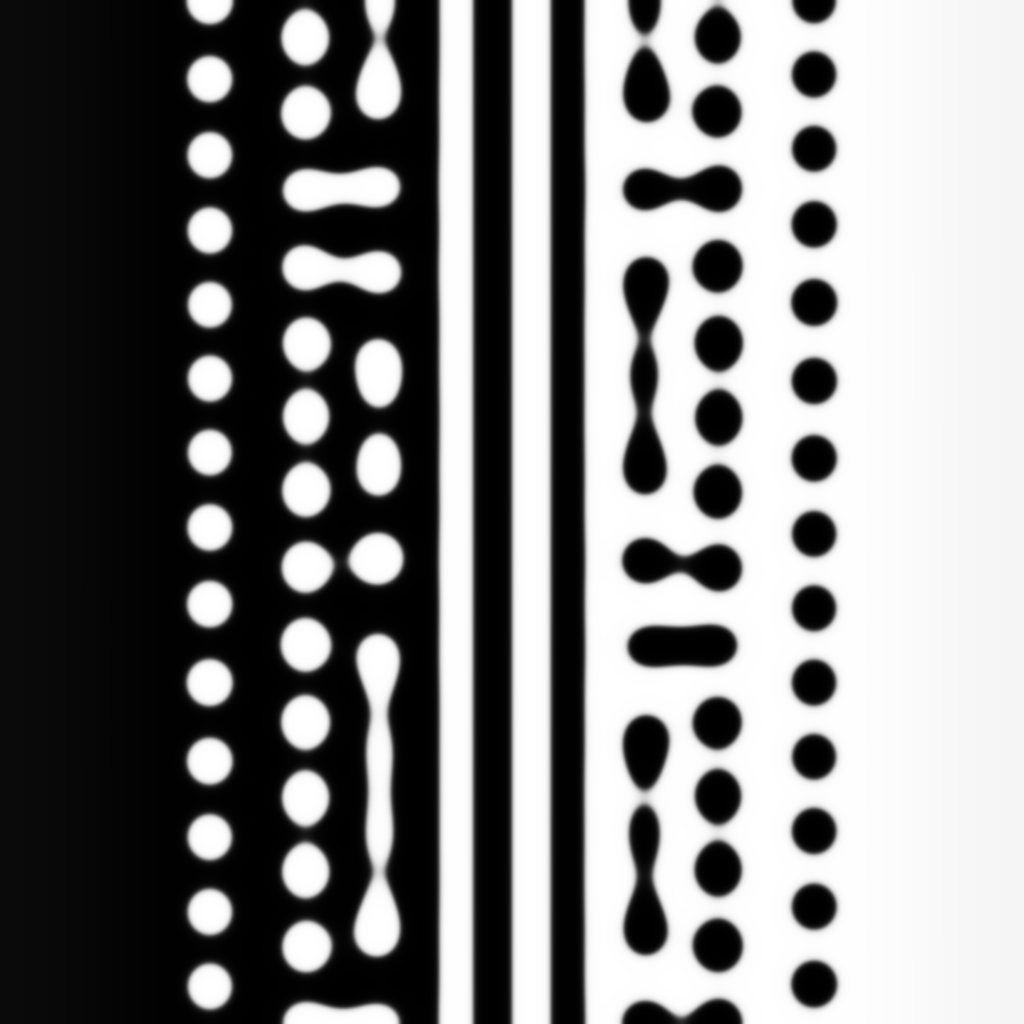}%
   \hfill{}%
   \includegraphics[width=0.3\textwidth]{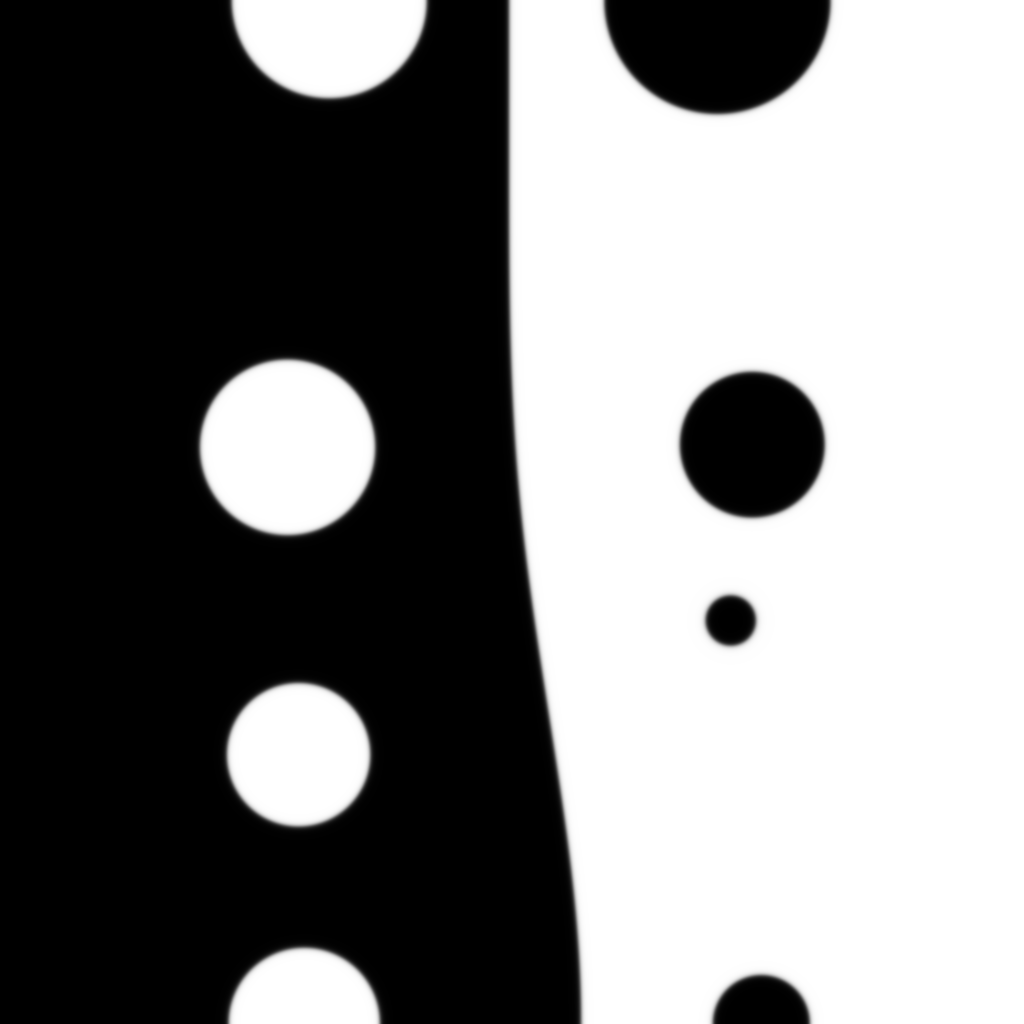}%
   \caption{%
    The black-to-white initial gray distribution forms very fast
    small structures (left).  Rounding breaks the symmetry but
    objects' curvature and distance are defining the evolution
    (center).  The object with the smallest diameter is the next
    to vanish (right);  the steady state straight borderline is still far away
    on the linear time scale%
    }\label{fig:bubbles}%
\end{figure}

\begin{figure}[t]%
   \centering{}%
   \includegraphics[height=0.26\textheight, clip]{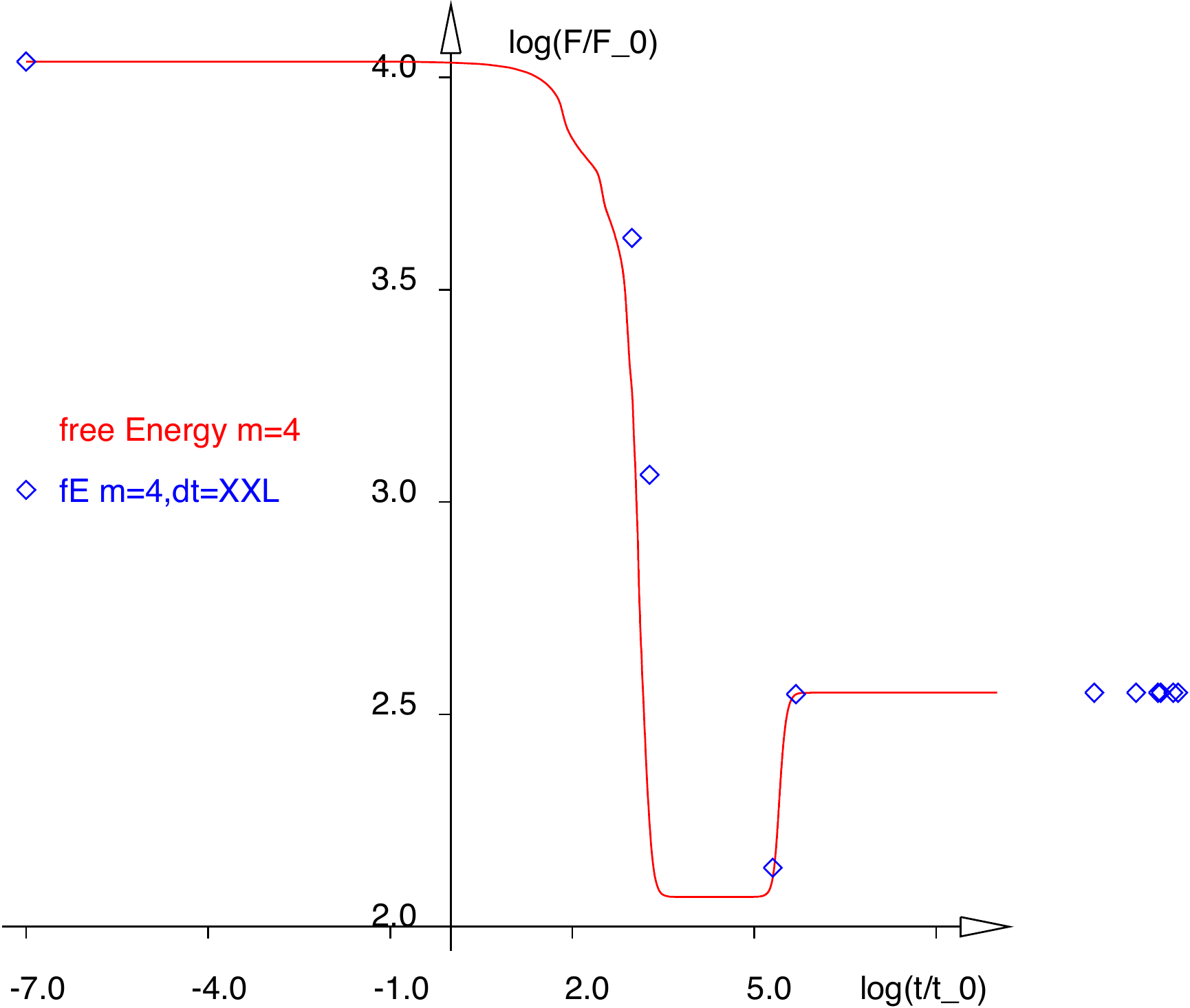}%
   \hspace{0.1\linewidth}%
   \includegraphics[width=0.26\textheight, clip]{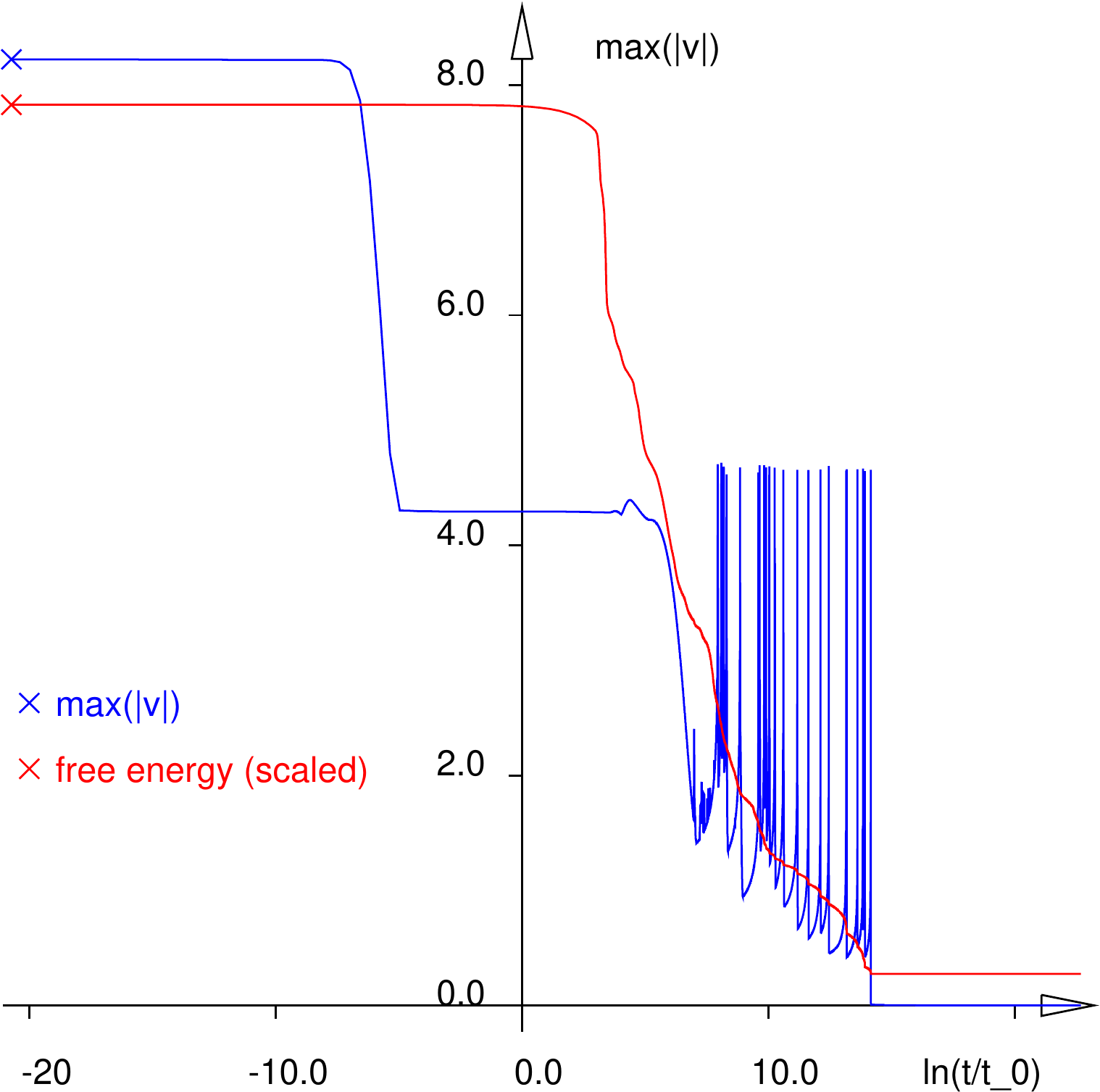}%
   \caption{%
      The free energy is increasing for the non-dissipative
      discretization (for small and large time steps, left) and
      decaying finally in staircase form for the dissipative and
      bounds preserving discretization (right).  Each peek in the
      time derivative (blue, equal hight means equal maximal
      vanishing speed) is indicating the finally very rapidly vanishing,
      compared to all other time scales at this time,
      of a white or black object in the mainly anti-particle domain.
      The stretching of the initially curved final border line
      defines the longest time scale%
   }\label{fig:bubbles:plot}%
\end{figure}

The dissipative scheme approximates the expected minimal surface.
After $10^{20}$ in time, the initial gray distribution is
transformed through a striped pattern with bubbles into a black and
white square halves divided by a straight vertical line in the middle of the square.
\Cref{fig:bubbles,fig:bubbles:plot} show a typical evolution process:
thin stripes (\cref{fig:bubbles}, left) decay into circles with a typical diameter related to the
stripe width (\cref{fig:bubbles}, center). The symmetry is broken because of the rounding, small
perturbations in the initial data, and because $\Omega$ is large compared to $\sigma$ and not compatible with arbitrary stripe width.
The  longest time scale is defined by the
time which is required to straighten the borderline between these
halves to a straight line.

If the discretization is non-dissipative, then
the asymptotic minimal surface will be replaced by a much larger surface.
The free energy is maximal for the initial value.  Close to the
end it decays and,  after some orders of magnitude in time,
looks like a constant.  The remaining gradient in the solution,
very small almost everywhere, is sufficient to produce some free
energy at some position in the mesh in the non-dissipative
discretization. This causes larger gradients and larger free energy
production. In the mostly dissipative volume it is resulting in a
much larger separating surface because the strong dissipation along
the surface is needed to compensate the production of free energy.
At the end, the computed steady state can be far from that of the true minimal energy.

%

%--- Athena-WFI ------------------------------------------------------
\subsection{Active pixel DePFET detectors}

Active pixel DePFET detectors are silicon sensors based on depleted $p$-channel field-effect transistors (DePFETs)~\cite{KemLut87,LutAndEck07}, which are p-channel field-effect transistors (FETs) integrated onto a fully depleted silicon bulk.
Entering particles generate electrons in the bulk, which are then collected at the internal gate below the transistor channel.
This increases the conductivity of the FET proportionally to the amount of the stored charge, providing a measure for the photon energy.
The unique advantages of detectors based on this technology are the amplification and storage of the signal charge directly inside the pixels, non-destructive readouts, very low noise, and high power-efficiency, which makes them highly suitable for X-ray astronomy and particle physics~\cite{LutAndEck07,TreAndHau18}.

An example of such a detector is the Wide Field Imager (WFI) for the Advanced Telescope for High-Energy Astrophysics (Athena)~\cite{athenaUrl,athena} of the European Space Agency.
The WFI is a solid-state imaging spectrometer for X-ray images in the 0.1--\SI{15}{\kilo\electronvolt} energy band with the simultaneous spectrally and time-resolved photon detection.
It is developed by a collaboration of the Halbleiterlabor of the Max-Planck-Society (MPG HLL), Erlangen Centre for Astroparticle Physics, the Institute for Astronomy and Astrophysics Tübingen, the Dresden University of Technology, PNSensor GmbH, the University of Leicester, and the \foreignlanguage{french}{L'Institut de Recherche en Astrophysique et Planétologie} Toulouse~\cite{athena}.
The numerical simulation was carried out in collaboration of the MPG HLL with the Weierstrass Institute for Applied Analysis and Stochastics.

In semiconductors, in comparison to the previous example, the particle interaction is reversed: particles of different types (holes and electrons) attract each other whereas diffusion is the repelling force for the particles of the same type.
The mathematical model is given by the van~Roosbroeck equations
\begin{align*}
   -\nabla \cdot \varepsilon \nabla w 
   &= C - n - p
   ,\\
   \frac{\partial n}{\partial t} - \nabla \cdot n_i \mu_n e^w \nabla e^{-\phi_n} 
   &= R
   ,\\
   \frac{\partial p}{\partial t} - \nabla \cdot n_i \mu_p e^{-w} \nabla e^{\phi_p}
   &= R,
\end{align*}
where $w$ is the electrostatic potential, $n = n_i e^{w - \phi_n}$ and  $p = n_i e^{\phi_p - w}$
are the electron and hole densities, $C$ is doping, $\phi_n$ and $\phi_p$ are quasi-Fermi potentials
and $R = r(\bm{x},n,p) \left( n_i^2 - np \right)$ is the recombination/generation rate with 
$r(\bm{x},n,p)>0$.

\begin{figure}[p]%
   \includegraphics[width=0.47\textwidth]{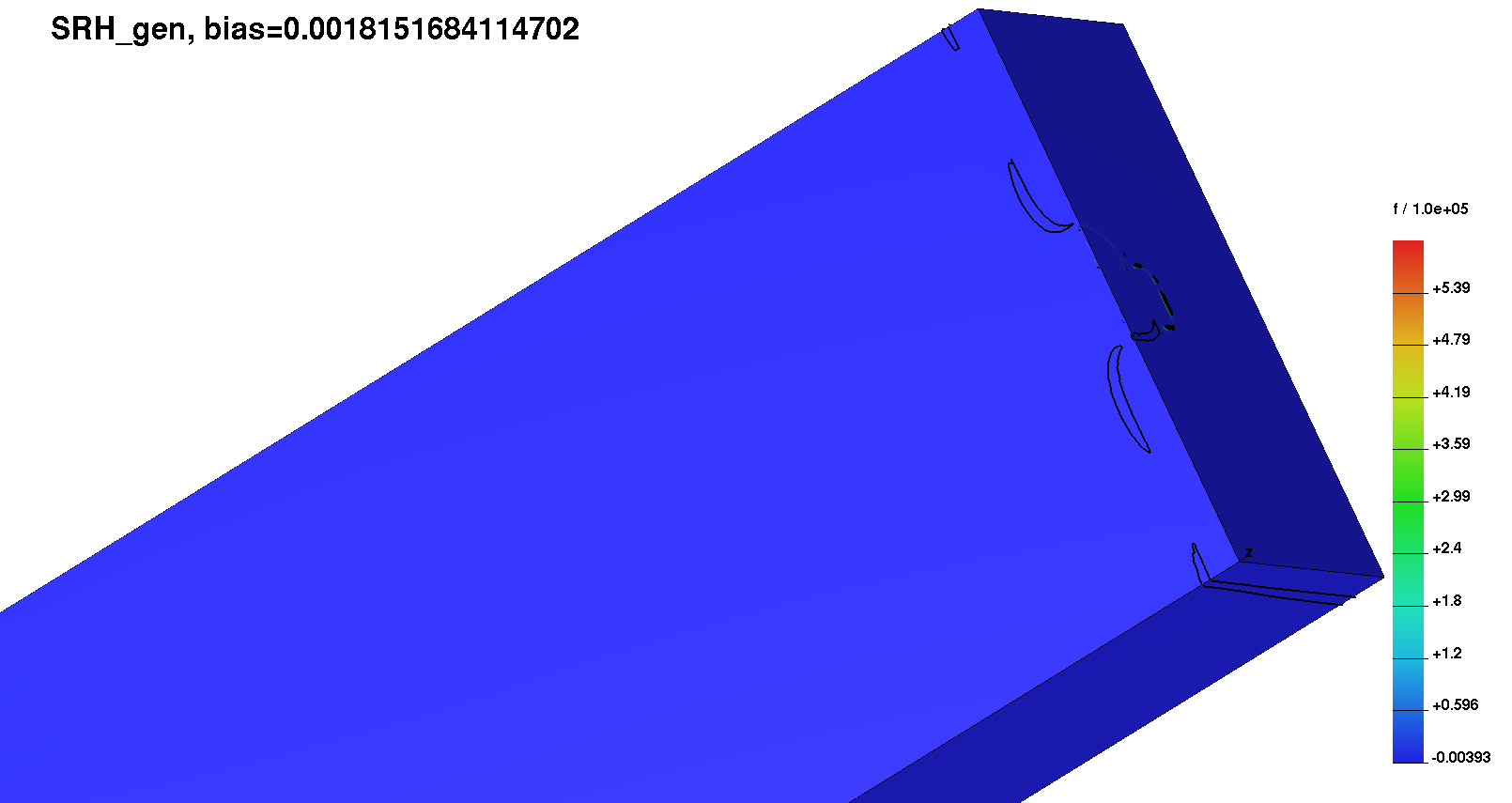}%
   \hfill{}%
   \includegraphics[width=0.47\textwidth]{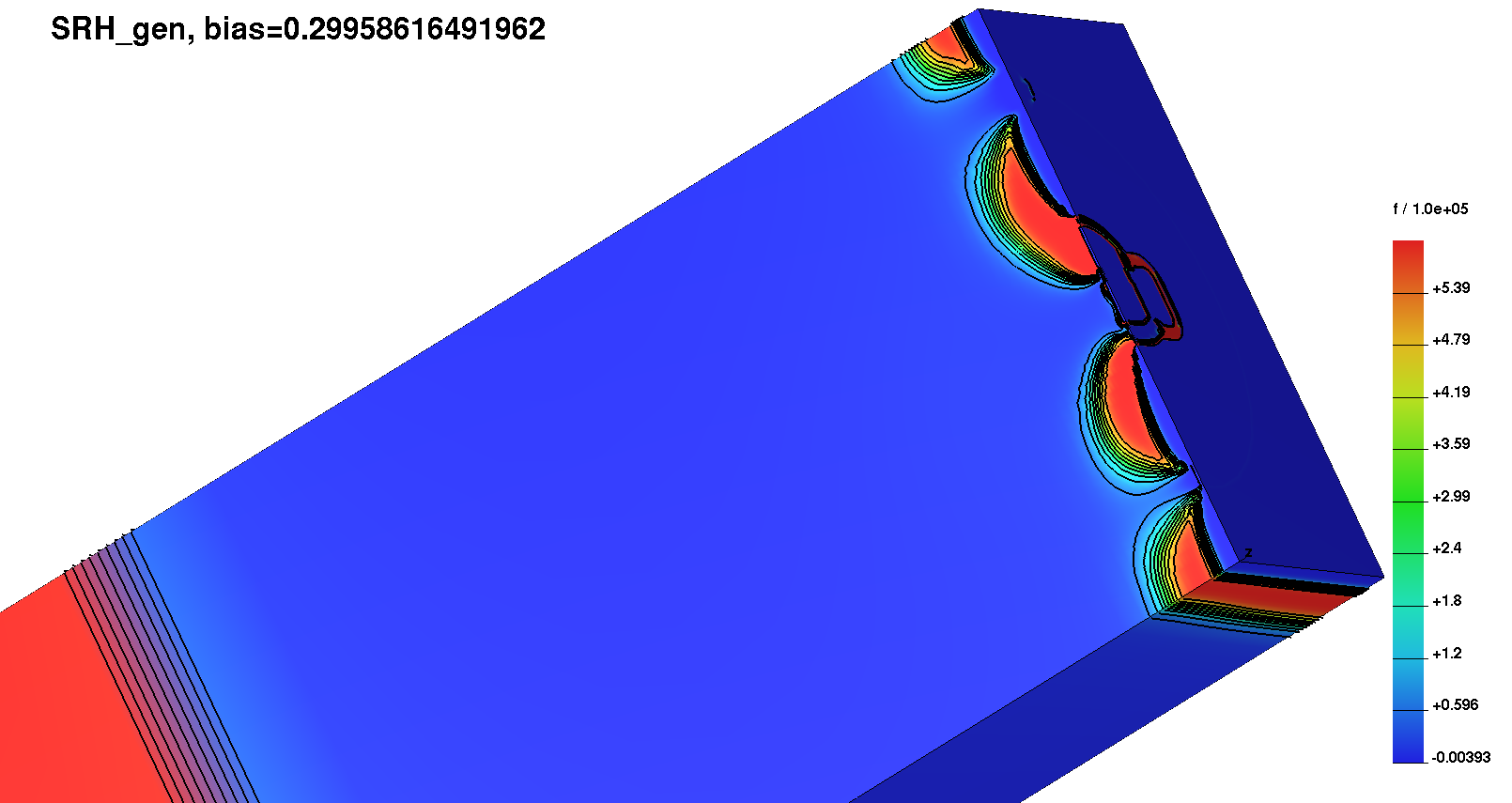}%
   \\[0.25\baselineskip]%
   \includegraphics[width=0.47\textwidth]{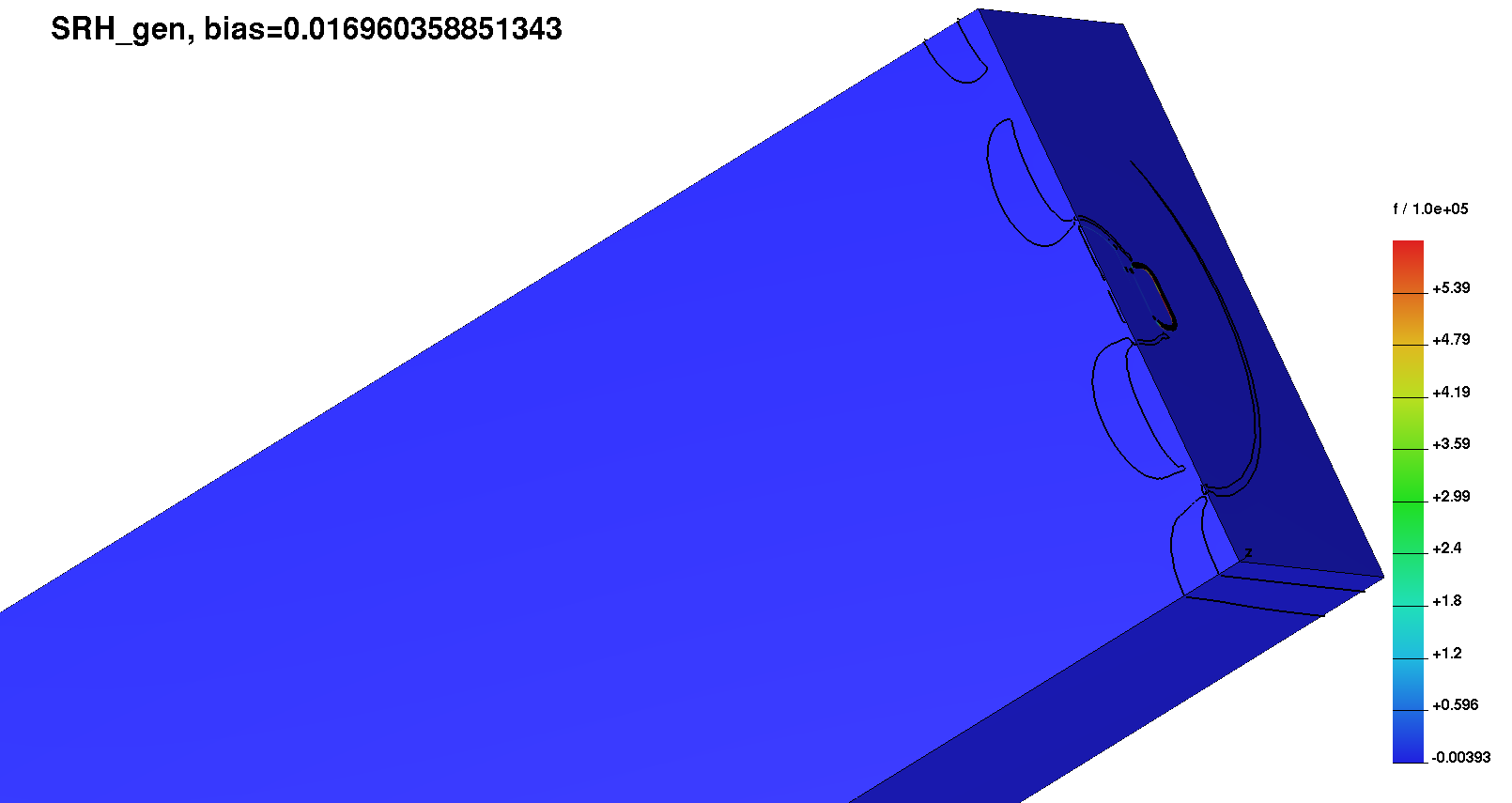}%
   \hfill{}%
   \includegraphics[width=0.47\textwidth]{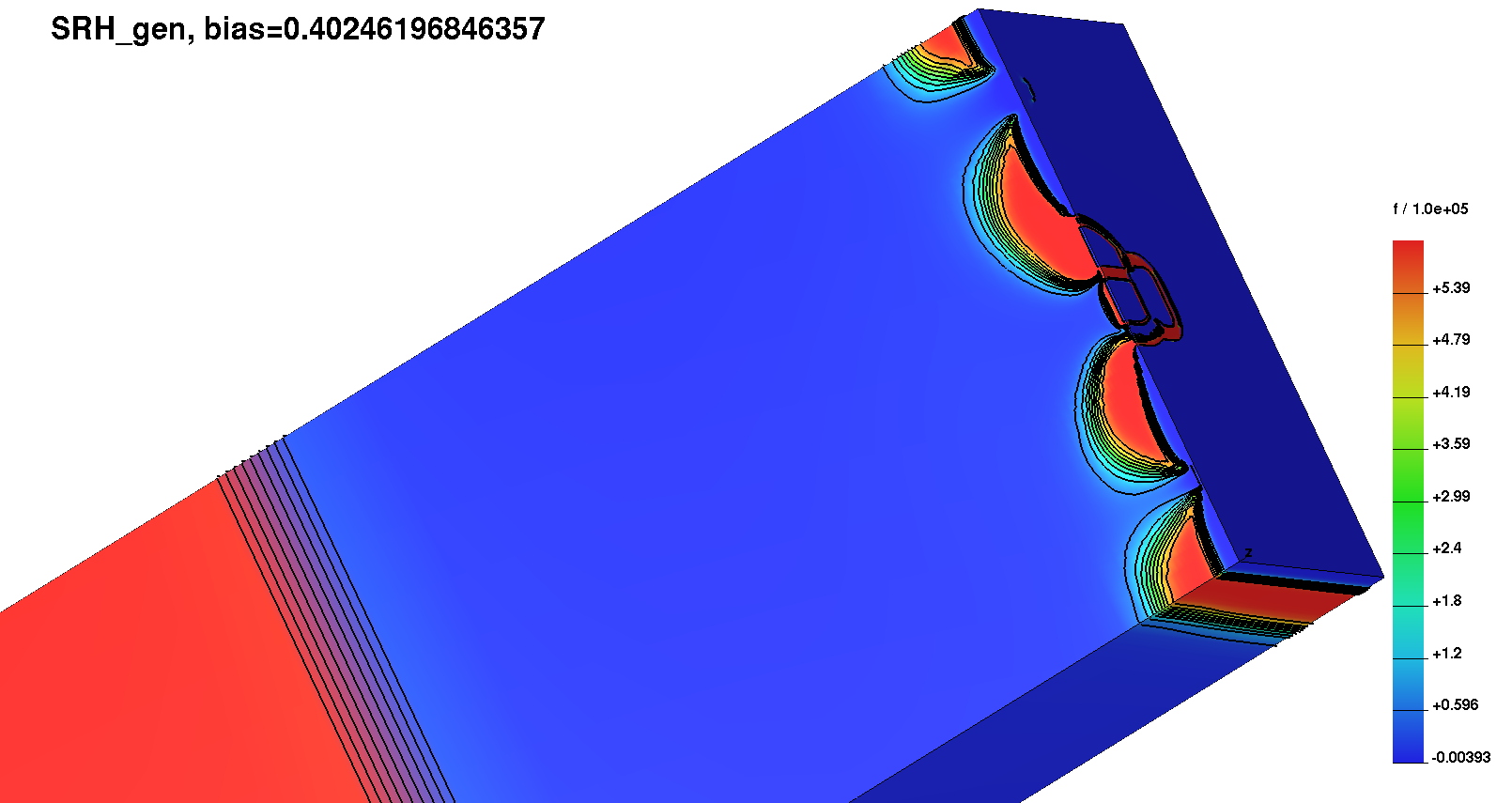}%
   \\[0.25\baselineskip]%
   \includegraphics[width=0.47\textwidth]{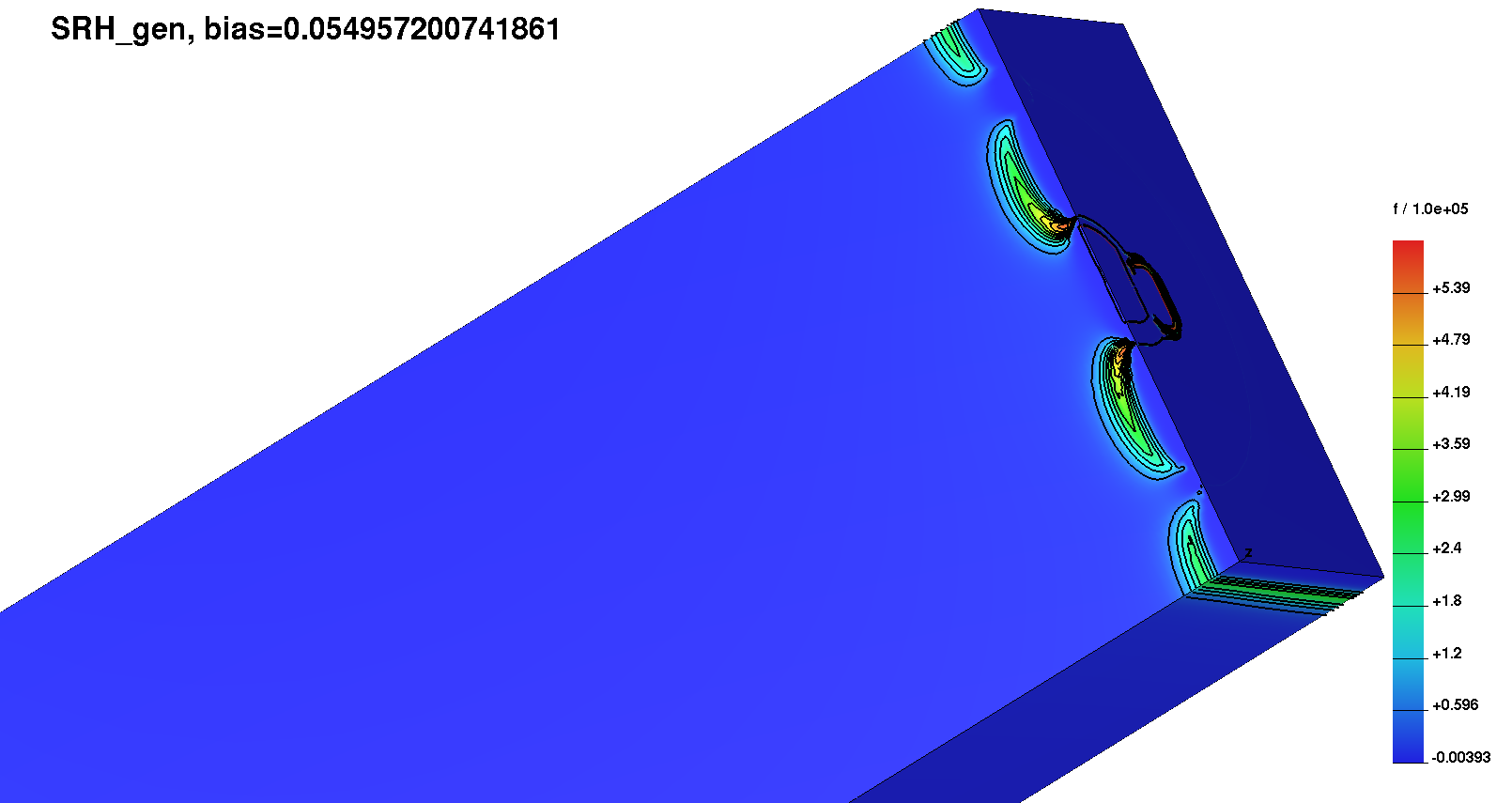}%
   \hfill{}%
   \includegraphics[width=0.47\textwidth]{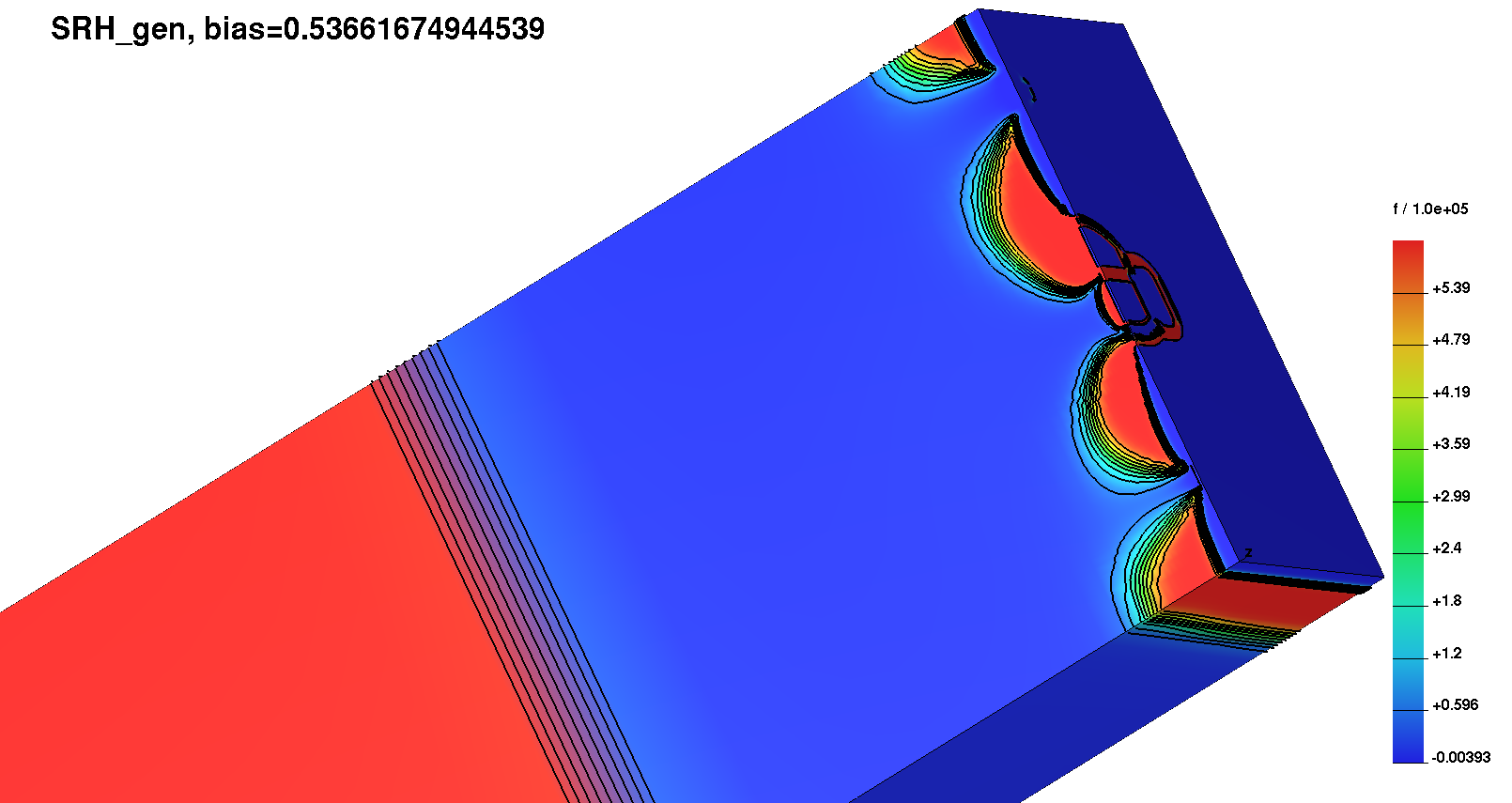}%
   \\[0.25\baselineskip]%
   \includegraphics[width=0.47\textwidth]{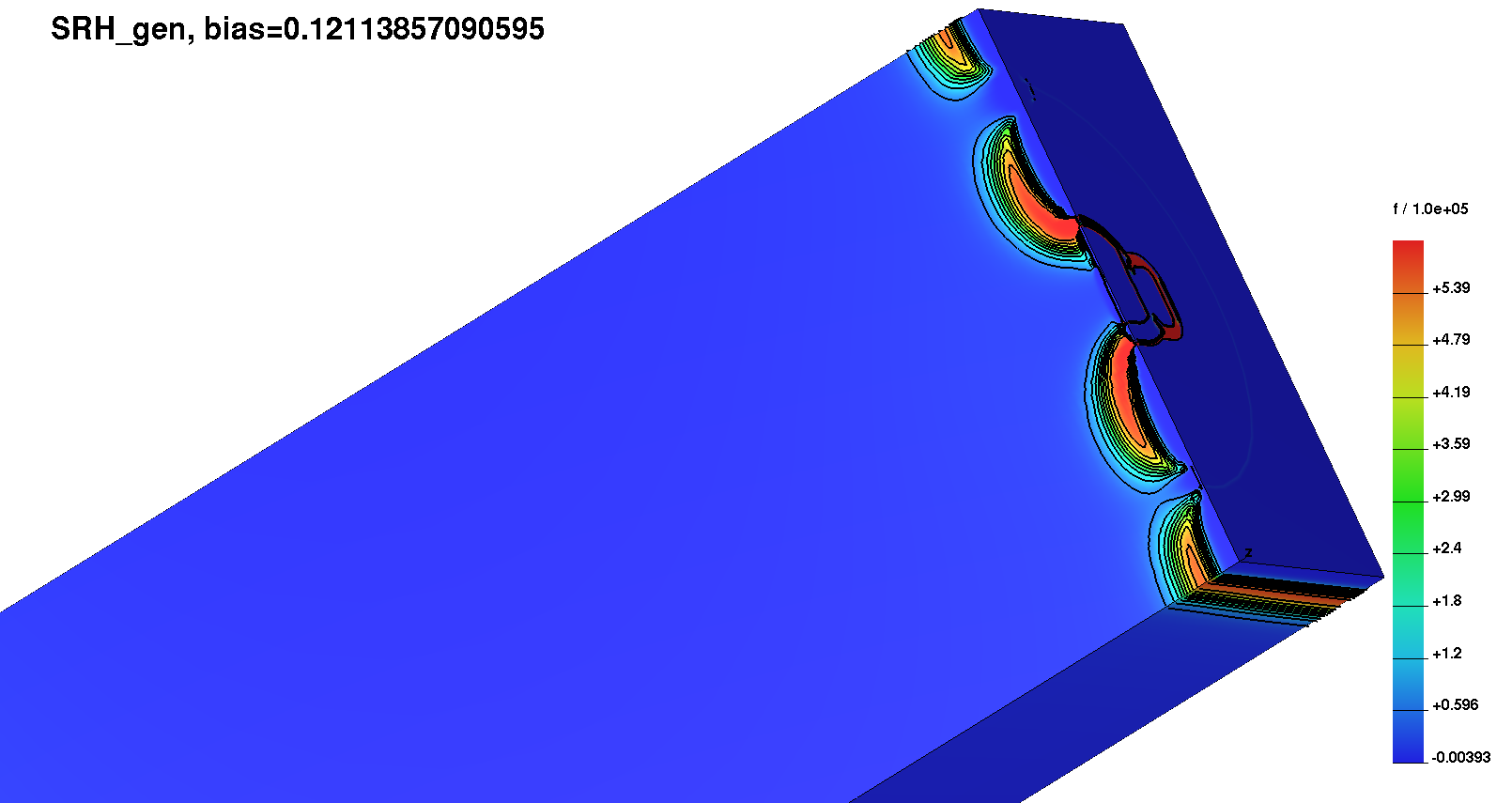}%
   \hfill{}%
   \includegraphics[width=0.47\textwidth]{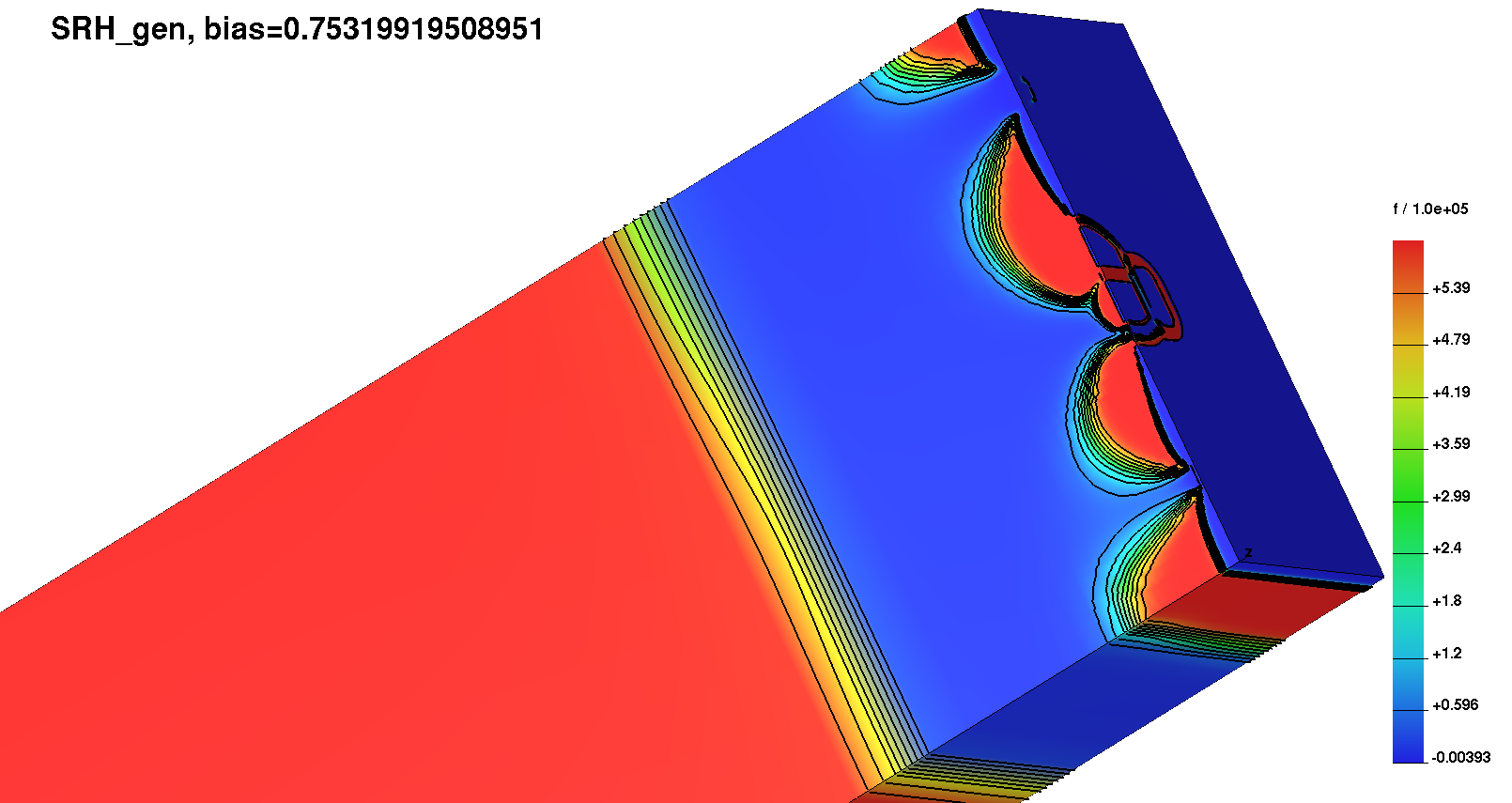}%
   \\[0.25\baselineskip]%
   \includegraphics[width=0.47\textwidth]{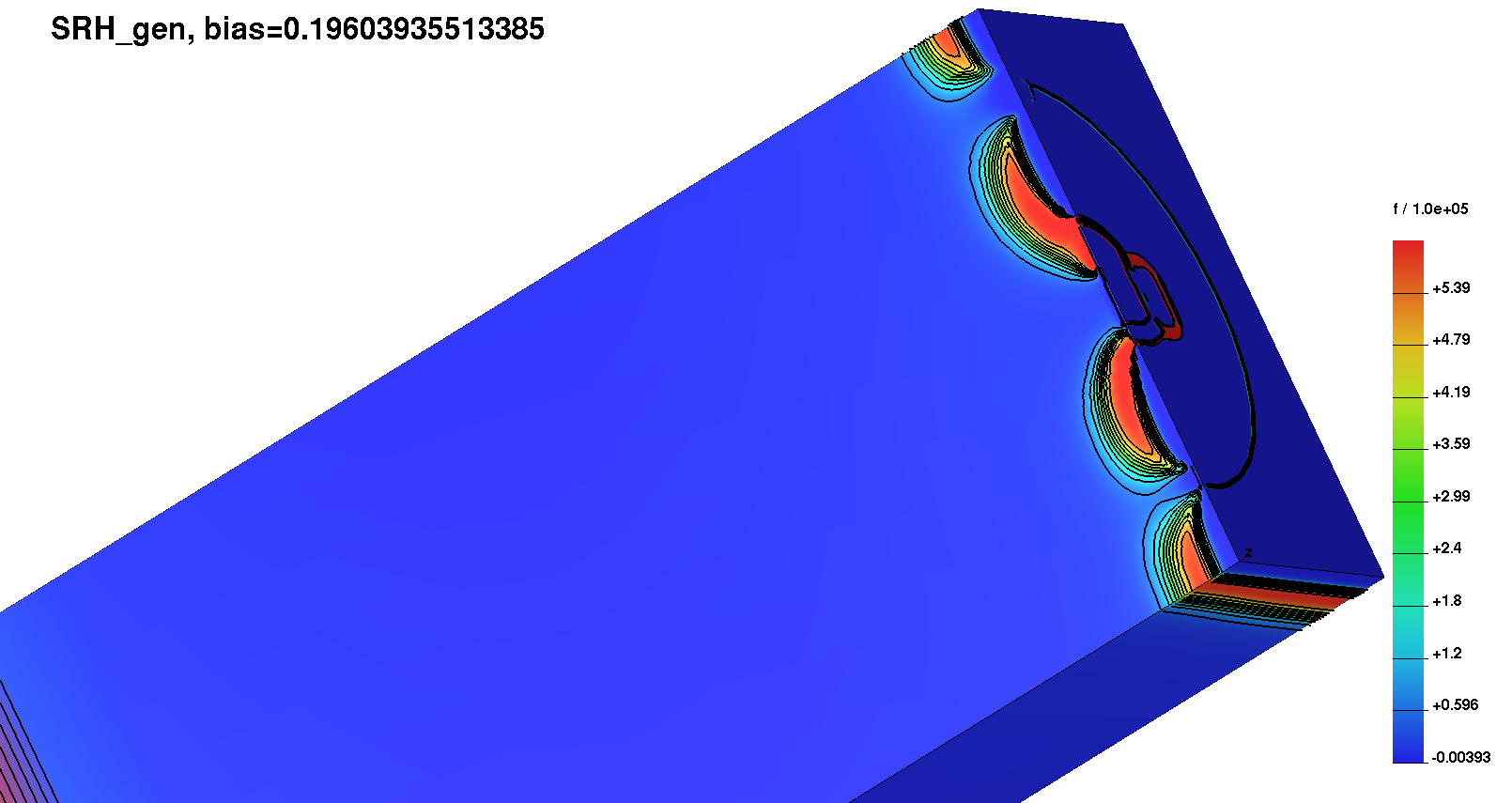}%
   \hfill{}%
   \includegraphics[width=0.47\textwidth]{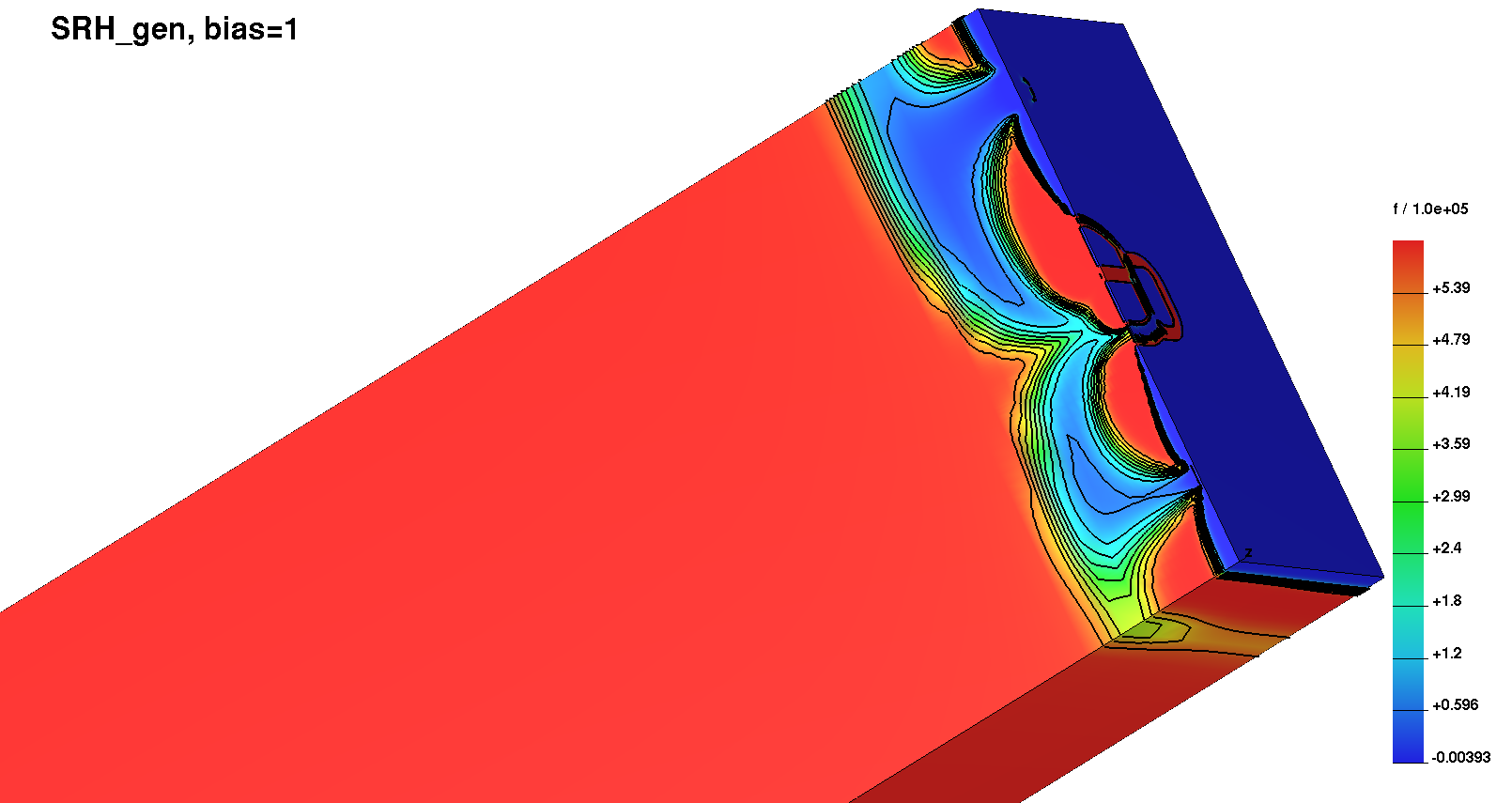}%
   \caption{Athena WFI detector: simulation of the depletion front propagation}\label{fig:athena}
\end{figure}

For the numerical simulation, it is essential for the discrete system to preserve essential properties of the analytical system, such as the maximum principle, positivity of the carrier concentrations, and current conservation.
The Voronoi FVM allows to carry over this properties to the discretized equations.
Numerical computations are done with the well-known Scharfetter–Gummel scheme (for a rigorous treatment see, e.g.,~\cite{Gar09}).
\cref{fig:athena} shows the simulation of the depletion front propagation in a pixel of the sensor.
Because of the symmetry of the domain the computation is carried out on one half a pixel.

%--- Summary ----------------------------------------
\section{Concluding remarks}

Delaunay-Voronoi meshing offers a chance to approximate the geometry
of a problem while preserving the essential qualitative properties
of the analytic solution.  It is difficult to imagine another way
to reach more than just a short time approximations with discrete
schemes that violate the essential stability properties of the
analytic problem such as solution bounds, dissipativity for all
concentrations, and any space and time step sizes.

The knowledge is not always a strictly increasing function in time
and the relations between analytic properties of (nonlinear)
PDE systems, finite volume methods, Delaunay meshes, and mesh
generation may have been stronger in the past.  During the last
decades, boundary conforming Delaunay meshes seem to be an uncommon
species but there is good hope that they will experience a new
revival~\cite{GarKudTsv19}.

%\newpage{}


\begin{thebibliography}{10}

\bibitem{athenaUrl}%
\newblock{}%
Athena --- Advanced Telescope for High-Energy Astrophysics.
\url{https://sci.esa.int/athena}

\bibitem{athena}%
\newblock{}%
Athena. The extremes of the universe: from black holes to large-scale structure.
\newblock{}%
\emph{Assessment study report ESA/SRE(2011)17}, European Space Agency, 2012.
\url{http://sci.esa.int/jump.cfm?oid=49835}

\bibitem{AllenSouthwell55}%
D.~Allen and R.~Southwell.
\newblock{}%
Relaxation methods applied to determine the motion, in two
  dimensions, of a viscous fluid past a fixed cylinder.
\newblock{}%
\emph{Quart. J. Mech.\ and Appl. Math.}, 8:129–145, 1955.

\bibitem{BanRos87}%
R.~E. Bank and D.~J. Rose.
\newblock{}%
Some error estimates for the box method.
\newblock{}%
\emph{SIAM J. Numer. Anal.}, 24(4):777--787, 1987.

\bibitem{Chew89}%
P.~L. Chew.
\newblock{}%
Constrained {D}elaunay triangulations.
\newblock{}%
\emph{Algorithmica}, 4(1):97--108, 1989.

\bibitem{Delaunay}%
B.~Delaunay.
\newblock{}%
\foreignlanguage{french}{Sur la sph\'ere vide. A la mémoire de Georges Voronoï.}
\newblock{}%
\emph{Izvestia Akademii Nauk SSSR.\@ Otd. Matem.\ i Estestv. Nauk}, 6:793--800, 1934.

\bibitem{EymardGallouetHerbin97}%
R.~Eymard, T.~Gallouet, and R.~Herbin.
\newblock{}%
\emph{Finite volume methods}.
\newblock{}%
P. G. Ciarlet and J. L. Lions, Ed., Handbook of Numerical Analysis.
  Elsevier Science B.V., Amsterdam, 1997.

\bibitem{PFleischmann-99}%
P.~Fleischmann.
\newblock{}%
{M}esh generation for {T}echnology {CAD} in three dimensions.
\newblock{}%
\emph{Dissertation}, Technische Universität Wien, 1999.

\bibitem{GabSok69}%
K.~R. Gabriel and R.~R. Sokal.
\newblock{}%
A new statistical approach to geographic variation analysis.
\newblock{}%
\emph{Syst. Biol.}, 18(3):259--278, 1969.

\bibitem{GajGar05a}%
H.~Gajewski and K.~Gärtner.
\newblock{}%
A dissipative discretization scheme for a nonlocal phase segregation model.
\newblock{}%
\emph{Z. Angew. Math. Mech.}, 85:815--822, 2005.

\bibitem{GajGar05}%
H.~Gajewski and K.~Gärtner.
\newblock{} On a nonlocal model of image segmentation.
\newblock{} \emph{Z. Angew. Math. Phys.}, 56:572--591, 2005.

\bibitem{GajGri06}%
H.~Gajewski and J.~A. Griepentrog.
\newblock{}%
A descent method for the free energy of multicomponent systems.
\newblock{}%
\emph{Discrete Contin. Dyn. Syst.}, 15(2):505--528, 2006.

\bibitem{GarKudTsv19}%
V.~Garanzha, L.~Kudryavtseva, and V.~Tsvetkova.
\newblock{}%
Structured orthogonal near-boundary {V}oronoi mesh layers for planar domains.
\newblock{}%
In \emph{Numerical Geometry, Grid Generation and Scientific Computing},
Lect. Notes Comput. Sci. Eng. 131. Springer International Publishing, 2019.
\newblock{}%
Proceedings of the 9th international conference NUMGRID-2018 / Voronoi-150,
celebrating the 150th anniversary of G.F. Voronoi, Moscow, Russia, December 2018.

\bibitem{Gar09}%
K.~Gärtner.
\newblock{}%
Existence of bounded discrete steady-state solutions of the van
{R}oosbroeck system on boundary conforming {D}elaunay grids.
\newblock{}%
\emph{SIAM J. Sci. Comput.}, 31(2):1347--1362, 2009.

\bibitem{ga15}%
K.~Gärtner.
\newblock{}%
Existence of bounded discrete steady state solutions of the van
  {R}oosbroeck system with monotone {F}ermi–{D}irac statistic functions.
\newblock{}%
\emph{J. Comput. Electron.}, 14(3), 2015.

\bibitem{GiaLeb97}%
G.~Giacomin and J.~L. Lebowitz.
\newblock{}%
Phase segregation dynamics in particle systems with long range
interactions. {I}. {M}acroscopic limits.
\newblock{}%
\emph{J. Statist. Phys.}, 87(1-2):37--61, 1997.

\bibitem{glgae-09}%
A.~Glitzky and K.~Gärtner.
\newblock{}%
Energy estimates for continuous and discretized electro-reaction-diffusion systems.
\newblock{}%
\emph{Nonlinear Anal.}, 70:788--805, 2009.

\bibitem{Gri04}%
J.~A. Griepentrog.
\newblock{}%
On the unique solvability of a non-local phase separation problem for multicomponent systems.
\newblock{}%
In \emph{Nonlocal elliptic and parabolic problems}, \emph{Banach Center Publ.}, 66:153--164,
Polish Acad. Sci. Inst. Math., Warsaw, 2004.

\bibitem{Ilin69}
A.~M. Il’in.
\newblock{} A difference scheme for a differential equation with a small
  parameter multiplying the second derivative.
\newblock{} \emph{Mat. Zametki}, 6:237–248, 1969.

\bibitem{KemLut87}%
J.~Kemmer and G.~Lutz.
\newblock{}%
New detector concepts.
\newblock{}%
\emph{Nucl. Instrum. Methods Phys. Res. A}, 253(3):365--377, 1987.

\bibitem{Ker96}%
T.~Kerkhoven.
\newblock{}%
Piecewise linear {P}etrov-{G}alerkin error estimates for the box method.
\newblock{}%
\emph{SIAM J. Numer. Anal.}, 33(5):1864--1884, 1996.

\bibitem{Kop14}%
N.~Kopteva.
\newblock{}%
Linear finite elements may be only first-order pointwise accurate on anisotropic triangulations.
\newblock{}%
\emph{Math. Comp.}, 83(289):2061--2070, 2014.

\bibitem{Kop20}%
N.~Kopteva.
\newblock{}%
How accurate are finite elements on anisotropic triangulations in the
  maximum norm?
\newblock{}%
\emph{J. Comput. Appl. Math.}, 364:112316, 2020.

\bibitem{KoFHPS00}%
R.~Kosik, P.~Fleischmann, B.~Haindl, P.~Pietra, and S.~Selberherr.
\newblock{}%
On the interplay between meshing and discretization in three-dimensional diffusion simulation.
\newblock{}%
\emph{IEEE Trans. Comput.-Aided Des. Integr. Circuits Syst.}, 19(11):1233--1240, 2000.

\bibitem{Lawson77}%
C.~L. Lawson.
\newblock{}%
Software for ${C}^1$ surface interpolation.
\newblock{}%
In \emph{Mathematical Software III}, pages 161--194, 1977.

\bibitem{LutAndEck07}%
G.~Lutz, L.~Andricek, R.~Eckardt, O.~H{\"a}lker, S.~Hermann, P.~Lechner,
  R.~Richter, G.~Schaller, F.~Schopper, H.~Soltau, L.~Str{\"u}der, J.~Treis,
  S.~W{\"o}lfl, and C.~Zhang.
\newblock{}%
{DePFET}--detectors: {N}ew developments.
\newblock{}%
\emph{Nucl. Instrum. Methods Phys. Res. A}, 572(1):311--315, 2007.

\bibitem{Macneal53}%
R.~H. MacNeal.
\newblock{}%
An asymmetrical finite difference network.
\newblock{}%
\emph{Quart. Math. Appl.}, 11:295--310, 1953.

\bibitem{Rip90}%
S.~Rippa.
\newblock{}%
Minimal roughness property of the {D}elaunay triangulation.
\newblock{}%
\emph{Comput. Aided Geom. Design}, 7(6):489--497, 1990.

\bibitem{Sch_Gu}%
D.~L. Scharfetter and H.~K. Gummel.
\newblock{}%
Large-signal analysis of a silicon read diode oscillator.
\newblock{}%
\emph{IEEE Trans. Electr. Dev.}, 16:64--77, 1969.

\bibitem{She02b}%
J.~R. Shewchuk.
\newblock{}%
Delaunay refinement algorithms for triangular mesh generation.
\newblock{}%
\emph{Comput. Geom.}, 22(1-3):21--74, 2002.

\bibitem{She08}%
J.~R. Shewchuk.
\newblock{}%
General-dimensional constrained {D}elaunay and constrained regular
  triangulations. {I}. {C}ombinatorial properties.
\newblock{}%
\emph{Discrete Comput. Geom.}, 39(1-3):580--637, 2008.

\bibitem{Si15}%
H.~Si.
\newblock{}%
{TetGen}, a {D}elaunay-based quality tetrahedral mesh generator.
\newblock{}%
\emph{ACM Trans. Math. Softw.}, 41(2):11:1--11:36, 2015.

\bibitem{SiGarFuh10}%
H.~Si, K.~Gärtner, and J.~Fuhrmann.
\newblock{}%
Boundary conforming {D}elaunay mesh generation.
\newblock{}%
\emph{Comput. Math. Math. Phys.}, 50(1):38--53, 2010.

\bibitem{TreAndHau18}%
W.~Treberspurg, R.~Andritschke, G.~Hauser, P.~Lechner, N.~Meidinger,
J.~Ninkovic, J.~M{\"u}ller-Seidlitz, and F.~Schopper.
\newblock{}%
Measurement results of different options for spectroscopic {X}-ray {DePFET} sensors.
\newblock{}%
\emph{J. Instrum.}, 13(9):P09014--P09014, 2018.

\bibitem{Varga62}%
R.~S. Varga.
\newblock{}%
Matrix Iterative Analysis.
\newblock{}%
Prentice-Hall, 1962.
\newblock{}%
Englewood Cliffs, N. J.

\bibitem{Voronoi07}%
G.~Voronoi.
\newblock{}%
\foreignlanguage{french}{Nouvelles applications des param\`etres continus \`a la th\'eorie des
   formes quadratiques.}%
\newblock{}%
\emph{Reine Angew. Math.}, 133:97–178, 1907.

\bibitem{XuZik99}%
J.~Xu and L.~Zikatanov.
\newblock{}%
A monotone finite element scheme for convection-diffusion equations.
\newblock{}%
\emph{Math. Comp.}, 68(228):1429--1446, 1999.

\end{thebibliography}
\end{document}

%--- Title frame -----------------------------------------------------
\renewcommand*{\today}{\emph{\small{}31 июля 2019г.}}

\begin{frontmatter}
\title{%
   Зачем нужны сетки Вороного-Делоне?
   Основные свойства метода конечных объемов с использованием ячеек Вороного.%
}

\author{Клаус Гертнер}
\author{Леннард Каменски}
\address{m4sim GmbH, Берлин, Германия (\url{https://m4sim.de})}

\begin{abstract}%
В отличие от схем, которые локально нарушают существенные свойства устойчивости параболических и эллиптических задач, методы конечного объема с мультиматериальных областях  с использованием ячеек Вороного и сеток Делоне, согласованных с границей, в которых центры приграничных шаров Делоне принадлежат расчетной области, обеспечивают хорошее приближение геометрии задачи и способны при этом сохранять существенные качественные свойства решения при любом разрешении в пространстве и времени, а также при изменениях временных масштабов на несколько порядков величин.
К сожалению, в последнее время эти методы используются все реже и реже, поскольку отсутствуют автоматические генераторы сеток с нужными свойствами.  Данная работа дает краткое описание ключевых свойств методов конечного объема с использованием сеток Вороного-Делоне с примерами из практики применения и мотивацию для более широкого использования таких методов и для разработки алгоритмов построения сеток Делоне-Вороного.
\\[0.5\baselineskip]
{\footnotesize{}%
\textcopyright~2019г. Авторы.
Препринт Произведения, принятого для публикации в ЖВММФ,
Pleiades Publishing (\url{http://pleiades.online}).}
\\[0.5\baselineskip]
\emph{Ключевые слова:}
метод конечных объемов; сетки Вороного-Делоне; сетки Делоне, согласованные с границей.
\\
2010 MSC: 65M08, 65M50
\end{abstract}

\end{frontmatter}

\section{Introduction}

Многие прикладные проблемы приводят к сильно нелинейным эллиптическим
или параболическим дифференциальным уравнениям в частных производных.  Многие
из них изучались в течение длительного времени и для них имеются
аналитические результаты, зачастую основанные на вариационной
формулировке.  Использование энтропии или функционалов свободной
энергии позволяет оценить свойства решения и вывести  свойство понижения
свободной энергии без ограничительного предположения о малом размере
свободной энергии или типе граничных условий.  Как правило, должны быть
определены такие априорные оценки для решений, как позитивность и
ограниченность; может быть доказано понижение свободной энергии
вдоль траектории.

В дискретном случае это анализ менее очевиден.  Зачастую параметры
дискретизации во времени и пространстве, $h$ и $\tau$, которые на
самом деле являются пределами, рассматриваются как конечные и
универсально достижимые.  Однако это не всегда бывает в ситуациях
с ограниченной гладкостью, например, в состояниях удаленных от
термодинамического равновесия или при более высокой плотности
свободной энергии, и особенно если требуется приближение негладких
границ.  Решения могут иметь плато и внутренние или пограничные
слои, соединяющие плато и граничные значения.  Плотность энергии
часто очень велика по сравнению с привычной средой обитания человека.
Известным примером является перенос заряда в полупроводниках,
где большая разница в силе является результатом того, что разность
электростатического потенциала всего лишь в 1\,В соответствует
разности высот в \num{4000}\,км в гравитационном поле Земли
(механическое отношение эквивалента тепла на \SI{300}{\kelvin});
соответственно, 100\,В, что не так уж много для таких задач,
соответствует разнице высот от Земли до места где-то за Луной.
Многие задачи из области электрохимии и химических реакций, и из
других областей науки приводят к такому же типу проблем и требуют
сохранения качественных свойств в их дискретной форме.

В таких случаях классическая ошибка дискретизации имеет лишь
ограниченную значимость с точки зрения гарантии устойчивости
дискретного решения.  Часто даже физически не имеет значения,
достигается ли стационарное предельное состояние проблемы при
наименьшем масштабе времени в 1\,с за $10^{19}$ или $10^{20}$\,с,
поскольку время реакции и продолжительность жизни человека составляют
соответственно 1\,с и 100\,лет ($100 \text{\,лет}  \approx 100
\times 400 \times 25 \times \num{3600} \text{\,с} =
\num{3.6e+9}\text{\,с}$).  Напротив, большой интерес представляют
стационарные состояния и качественные характеристики, например,
обеспечивает ли модель А лучшее решение в гигагерцовом диапазоне, чем модель
B.  Поэтому существенным является сохранение структурных свойств
задач в пространстве и времени для конечных и не слишком малых
параметров дискретизации.

Более того, входные данные часто бывают неточны, и внутренние
слои могут разделять подобласти  с различными свойствами (генерация на одной
стороне слоя, рекомбинации с другой).  Оценка расположения этих подобластей с точностью
до процента от общего объема может быть достаточным.  Часто требуется просто найти
положение слоя, не разрешая его структуру.  До тех пор, пока дискретная
схема наследует свойства устойчивости аналитической задачи независимо
от размера пространственного и временного шага, требование разрешать все мелкие структуры
может быть неуместным.

Априори не ясно, какое сочетание дискретизации пространства-времени
обеспечит безусловную стабильность дискретных уравнений, достаточно
хорошо приближенных к исходной задаче, или будет ли такая схема
эффективна.  Особую роль здесь играют сетки Вороного-Делоне и неявный
метод  Эйлера.

Идеи устойчивой дискретизации для упрощенных задач восходят к 1950-м
годам~\cite{AllenSouthwell55,Ilin69,Sch_Gu} и в обобщённом виде
остаются актуальными и по сей день.  \emph{Диаграмма
Вороного}~\cite{Voronoi07} и соответствующая ей
\emph{триангуляция Делоне}~\cite{Delaunay}, использовались уже
в~\cite{Macneal53} (и, например, в классической книге
Варги~\cite{Varga62}).  Двумерные триангуляции Делоне являются оптимальными
для различных критериев, например, для максимизации наименьшего угла
триангуляции~\cite{Lawson77}, для минимизации квадрата $L^2$-нормы
градиента кусочно-линейной интерполяции, и для других критериев
(см., например,~\cite[Sect.~4]{She08} и список литературы к статье).
Трудно проследить историю развития сеток Делоне в деталях, но
основные моменты были связаны с $C^1$ поверхностной интерполяцией
и с такими задачами, как перенос нейтронов и перенос зарядов
в полупроводниках.  \emph{Допустимые сетки}~\cite{EymardGallouetHerbin97}
подводят итог для минимально необходимых свойств сеток, позволяющих доказать сходимость аппроксимаций законов сохранения.

По-видимому, местом рождения термина  \emph{сетки
Делоне, согласованные с границами~--- boundary conforming Delaunay meshes} была Вена (см., например,~\cite{PFleischmann-99}). Такие сетки требуются в
случае межматериальных границ (интерфейсов) со скачками свойств материала (например,
коэффициентов диффузии), нелинейных граничных условий третьего рода,
или в некоторых других особых случаях, чтобы гарантировать, что границы
и внутренние интерфейсы области задачи, являются частью  сетки
Вороного-Делоне, в то время как \emph{сетки Делоне с ограничениями~--- constrained Delaunay meshes}
~\cite{Chew89}, ставшие популярными в последнее
время~\cite{She02b,Si15}, для этой цели не подходят.  Локальные
ортогональные координаты были бы лучшим выбором. Несмотря на
то, что они часто не существуют глобально, т.е. во всей области,  локально они
могут стать ключом к сочетанию устойчивости, хороших аппроксимирующих
свойств и невысокой вычислительной сложности дискретной задачи.

Ниже перечислены некоторые желаемые свойства дискретной задачи по сравнению с дифференциальной постановкой:
\setlist{nolistsep}
\begin{itemize}[noitemsep]
\item те же качественные свойства устойчивости, что и у точной задачи;
\item никаких искусственных предположений о малости или гладкости каких-то величин;
\item использование слабой формы дискретной задачи, а также хорошо известных тестовых функций и алгоритмов;
\item следовательно, слабый дискретный энергетический функционал;
\item слабую сходимость решений когда правая часть задается как  распределение;
\item расположение внутренних слоев как конечных разрывов в решении  (случай неразрешенных слоев).
\end{itemize}
Примерами, где список частично выполнен, являются: перенос электронов и дырок в полупроводниках и разделение фаз (не Cahn-Hilliard)~\cite{GajGar05a,GajGar05,Gar09,ga15,glgae-09,GajGri06}.

%--- Definitions ----------------------------------------
\section{Сетки Делоне, ячейки Вороного, и сетки Делоне, согласованные с границей}

Рассмотрим дискретизацию (\emph{сетку}) $\Th(\Omega)$ области $\Omega \subset \mathbb{R}^d$ в предположении, что $\Omega$~--- это объединение конечных $d$-мерных подобластей, $\Omega = \cup_i \Omega_i$.
Каждая подобласть содержит только один материал и состоит из $d$-мерных симплексов. $\bm{E}^d_j$, $\Omega_i = \cup_j \bm{E}^d_j$.
Вершины симплексов сетки обозначим через $\bm{v}_{k}$.
Следовательно, все интерфейсы между материалами и границы области составлены из симплексов меньшей размерности $\bm{E}^m_{k}$, $1 \le m < d$.
Предполагается, что описание геометрии области содержит все критические линии и точки.

%--- Delaunay ----------------------------------------
\subsection{Триангуляция Делоне}
Дискретизация  $\Th(\Omega)$ заданной $d$-мерной области $\Omega$ симплексами называется \emph{сеткой Делоне} или \emph{триангуляцией Делоне}, если она удовлетворяет свойству пустого шара, то есть если никакая вершина сетки не находится внутри описанного сферы любого симплекса, которая является границей шара, называемого \emph{шаром Делоне}:
\[
   \forall \bm{v}_i, \bm{E}^d_j
   \in \Th(\Omega),
   ~
   \bm{v}_i \not\in \mathrm{B^o}(\bm{E}^d_j)
   .
\]
Для заданной сетки, симплекс сетки называется \emph{симплексом Делоне}, если
его шар Делоне пуст (\cref{fig:delaunay}).  Ребро сетки
называется \emph{ребром Делоне}, если существует пустая сфера, содержащая
это ребро, но не имеющая других вершин сетки внутри себя.  Следовательно,
симплекс является Делоне тогда и только тогда, когда все его ребра
являются ребрами Делоне.  Ребра Делоне не обязательно являются самыми
короткими из возможных ребер (см. \cref{fig:delaunay}, справа),   но двумерная
триангуляция Делоне  максимизирует минимальный угол
триангуляции~\cite{Lawson77}, а также квадрат $L^2$-нормы градиента
кусочно-линейной интерполяции~\cite{Rip90}.

Триангуляция Делоне  однозначна, если вершины сетки находятся в общем положении, т.е. никакие $d+2$ вершины не
лежат на общей пустой $d$-сфере.

\begin{figure}[t]%
   \centering{}%
   \begin{tikzpicture}[scale = 0.75]%
      \useasboundingbox (-0.69, -1.67) rectangle (5.15,3.86);
      % define vertices 1, 2, 3, 4
      \tkzDefPoint(0.0,  0.0){p1}
      \tkzDefPoint(4.5, -0.3){p2}
      \tkzDefPoint(2.8,  3.2){p3}
      \tkzDefPoint(0.1,  2.5){p4}
      % triangles 1-2-3 (gray) & 2-3-4
      \tkzFillPolygon[fill = gray!15](p1,p2,p3)
      \tkzDrawPolygon[thin](p1,p2,p3,p4)
      \tkzDrawSegment[thin](p1,p3)
      % circumcircle through 1-2-3
      \begin{pgfinterruptboundingbox}
      \tkzCircumCenter(p1,p2,p3)\tkzGetPoint{c}
      \end{pgfinterruptboundingbox}
      \tkzDrawCircle[dotted](c,p1)
      % draw vertices
      \tkzDrawPoints(p1,p2,p3,p4)
      \tkzLabelPoint[left](p1){$\bm{v}_1$}
      \tkzLabelPoint[below right](p2){$\bm{v}_2$}
      \tkzLabelPoint[above](p3){$\bm{v}_3$}
      \tkzLabelPoint[above left](p4){$\bm{v}_4$}
   \end{tikzpicture}%
   \hfill{}%
   \begin{tikzpicture}[scale = 0.75]%
      \useasboundingbox (-0.69, -1.67) rectangle (5.15,3.86);
      % define vertices 1, 2, 3, 4
      \tkzDefPoint(0.0,  0.0){p1}
      \tkzDefPoint(4.5, -0.3){p2}
      \tkzDefPoint(2.8,  3.2){p3}
      \tkzDefPoint(0.1,  2.5){p4}
      % triangles 1-2-4 (gray) & 2-3-4
      \tkzFillPolygon[fill = gray!15](p1,p2,p4)
      \tkzDrawPolygon[thin](p1,p2,p3,p4)
      \tkzDrawSegment[thin](p2,p4)
      % circumcircle through 1-2-4
      \begin{pgfinterruptboundingbox}
      \tkzCircumCenter(p1,p2,p4)\tkzGetPoint{c}
      \end{pgfinterruptboundingbox}
      \tkzDrawCircle[dotted](c,p1)%
      % draw vertices
      \tkzDrawPoints(p1,p2,p3,p4)
      \tkzLabelPoint[left](p1){$\bm{v}_1$}
      \tkzLabelPoint[below right](p2){$\bm{v}_2$}
      \tkzLabelPoint[above](p3){$\bm{v}_3$}
      \tkzLabelPoint[above left](p4){$\bm{v}_4$}
   \end{tikzpicture}%
   \hfill{}%
   \begin{tikzpicture}[scale = 0.75]%
      %\useasboundingbox (-0.69, -1.67) rectangle (5.15,3.86);
      \useasboundingbox (-0.69, -1.67) rectangle (5.15,4.46);
      % define vertices 1, 2, 3, 4
      \tkzDefPoint( 0.0,  1.6){p1}
      \tkzDefPoint( 4.5, -0.3){p2}
      \tkzDefPoint( 2.3,  4.0){p3}
      \tkzDefPoint( 1.0,  3.3){p4}
      % triangles 1-2-4 (gray) & 2-3-4
      \tkzFillPolygon[fill = gray!15](p1,p2,p4)
      \tkzDrawPolygon[thin](p1,p2,p3,p4)
      \tkzDrawSegment[thin](p2,p4)
      \tkzDrawSegment[dashed](p1,p3)
      % circumcircle through 1-2-4
      \begin{pgfinterruptboundingbox}
      \tkzCircumCenter(p1,p2,p4)\tkzGetPoint{c}
      \end{pgfinterruptboundingbox}
      \tkzDrawCircle[dotted](c,p1)%
      % draw vertices
      \tkzDrawPoints(p1,p2,p3,p4)
      \tkzLabelPoint[left](p1){$\bm{v}_1$}
      \tkzLabelPoint[below right](p2){$\bm{v}_2$}
      \tkzLabelPoint[above](p3){$\bm{v}_3$}
      \tkzLabelPoint[above left](p4){$\bm{v}_4$}
   \end{tikzpicture}%
   \caption{%
      $\triangle \bm{v}_1 \bm{v}_2 \bm{v}_3$ удовлетворяет условию
      Делоне, потому что его окружность пуста (слева), в то время
      как $\triangle \bm{v}_1 \bm{v}_2 \bm{v}_4$ не удовлетворяет
      условию Делоне, потому что его описанная окружность содержит
      в себе $\bm{v}_3$ (в центре).  Ребро Делоне не обязательно
      является самым коротким из возможных: $\bm{v}_2 \bm{v}_4$
      удовлетворяет условию Делоне, но не является самым коротким
      из возможных (справа)%
   }\label{fig:delaunay}%
\end{figure}

%--- Voronoi volumes etc. ----------------------------------------
\subsection{Объемы, поверхности и ячейки Вороного}

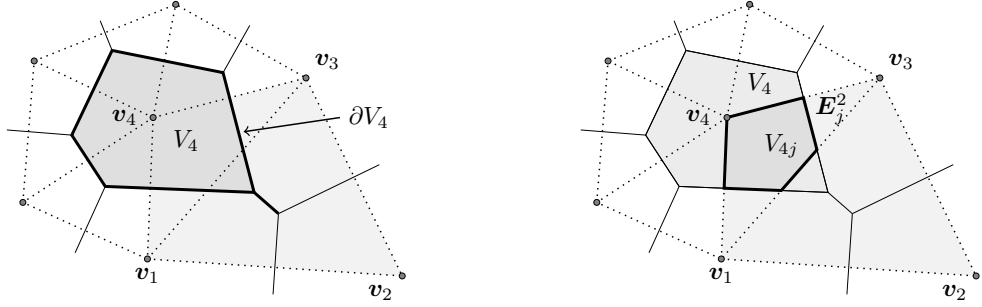
\begin{figure}[t]%
   \centering{}%
   \begin{tikzpicture}[scale = 0.8]%
      % define vertices 1, 2, 3, 4
      \tkzDefPoint(0.0,  0.0){p1}
      \tkzDefPoint(0.1,  2.5){p2}
      \tkzDefPoint(2.8,  3.2){p3}
      \tkzDefPoint(4.5, -0.3){p4}
      \tkzDefPoint( 0.5,  4.5){p5}
      \tkzDefPoint(-2.0,  3.5){p6}
      \tkzDefPoint(-2.2,  1.0){p7}
      \tkzDefPoint( 0.3,  2.1){labelVV}
      \tkzDefPoint( 3.4,  2.5){labelVS}
      % circumcenters
      \begin{pgfinterruptboundingbox}
      \tkzCircumCenter(p1,p2,p3)\tkzGetPoint{c123}
      \tkzCircumCenter(p1,p3,p4)\tkzGetPoint{c134}
      \tkzCircumCenter(p2,p3,p5)\tkzGetPoint{c235}
      \tkzCircumCenter(p2,p5,p6)\tkzGetPoint{c256}
      \tkzCircumCenter(p2,p6,p7)\tkzGetPoint{c267}
      \tkzCircumCenter(p2,p7,p1)\tkzGetPoint{c271}
      \end{pgfinterruptboundingbox}
      % midpoints
      \tkzDefMidPoint(p1,p4)\tkzGetPoint{mid14}
      \tkzDefMidPoint(p3,p4)\tkzGetPoint{mid34}
      \tkzDefMidPoint(p3,p5)\tkzGetPoint{mid35}
      \tkzDefMidPoint(p5,p6)\tkzGetPoint{mid56}
      \tkzDefMidPoint(p6,p7)\tkzGetPoint{mid67}
      \tkzDefMidPoint(p7,p1)\tkzGetPoint{mid71}
      \tkzDefMidPoint(c235,c123)\tkzGetPoint{vs}
      % voronoi volumes
      \tkzFillPolygon[fill = gray!10](p1,p2,p3,p4)
      \tkzFillPolygon[fill = gray!25](c123,c235,c256,c267,c271)
      \tkzDrawPolygon[very thick](c123,c235,c256,c267,c271)
      \tkzDrawLines[very thick, add = 0 and 0](c123,c134)
      \tkzDrawLines[thin, add = 0 and 1/2](c134,mid14)
      \tkzDrawLines[thin, add = 0 and 1/3](c134,mid34)
      \tkzDrawLines[thin, add = 0 and 1/2](c235,mid35)
      \tkzDrawLines[thin, add = 0 and 2/3](c256,mid56)
      \tkzDrawLines[thin, add = 0 and 1/2](c267,mid67)
      \tkzDrawLines[thin, add = 0 and 1/2](c271,mid71)
      % triangles 1-2-4 (dotted)
      \tkzDrawPolygon[dotted](p1,p4,p3,p5,p6,p7)
      \tkzDrawLines[  dotted, add = 0 and 0](p2,p1)
      \tkzDrawLines[  dotted, add = 0 and 0](p2,p3)
      \tkzDrawLines[  dotted, add = 0 and 0](p1,p3)
      \tkzDrawLines[  dotted, add = 0 and 0](p2,p5)
      \tkzDrawLines[  dotted, add = 0 and 0](p2,p6)
      \tkzDrawLines[  dotted, add = 0 and 0](p2,p7)
      % draw vertices
      \tkzDrawPoints(p1,p2,p3,p4,p5,p6,p7)
      \tkzLabelPoint[right](labelVV){$V_4$}
      \tkzLabelPoint[right](labelVS){$\partial V_4$}
      \tkzDrawSegment[add = 0 and -1/20, arrows = ->](labelVS, vs)
      \tkzLabelPoint[below](p1){$\bm{v}_1$}
      \tkzLabelPoint[left = 1mm](p2){$\bm{v}_4$}
      \tkzLabelPoint[above right](p3){$\bm{v}_3$}
      \tkzLabelPoint[below left](p4){$\bm{v}_2$}
   \end{tikzpicture}%
   \hspace{0.15\linewidth}%
   \begin{tikzpicture}[scale = 0.8]%
      % define vertices 1, 2, 3, 4
      \tkzDefPoint(0.0,  0.0){p1}
      \tkzDefPoint(0.1,  2.5){p2}
      \tkzDefPoint(2.8,  3.2){p3}
      \tkzDefPoint(4.5, -0.3){p4}
      \tkzDefPoint( 0.5,  4.5){p5}
      \tkzDefPoint(-2.0,  3.5){p6}
      \tkzDefPoint(-2.2,  1.0){p7}
      \tkzDefPoint( 0.3,  3.1){labelVV}
      \tkzDefPoint( 0.6,  2.0){labelVS}
      \tkzDefPoint( 1.5,  2.65){labelEj}
      % circumcenters
      \begin{pgfinterruptboundingbox}
      \tkzCircumCenter(p1,p2,p3)\tkzGetPoint{c123}
      \tkzCircumCenter(p1,p3,p4)\tkzGetPoint{c134}
      \tkzCircumCenter(p2,p3,p5)\tkzGetPoint{c235}
      \tkzCircumCenter(p2,p5,p6)\tkzGetPoint{c256}
      \tkzCircumCenter(p2,p6,p7)\tkzGetPoint{c267}
      \tkzCircumCenter(p2,p7,p1)\tkzGetPoint{c271}
      \end{pgfinterruptboundingbox}
      % midpoints
      \tkzDefMidPoint(p1,p2)\tkzGetPoint{mid12}
      \tkzDefMidPoint(p1,p4)\tkzGetPoint{mid14}
      \tkzDefMidPoint(p3,p4)\tkzGetPoint{mid34}
      \tkzDefMidPoint(p3,p5)\tkzGetPoint{mid35}
      \tkzDefMidPoint(p5,p6)\tkzGetPoint{mid56}
      \tkzDefMidPoint(p6,p7)\tkzGetPoint{mid67}
      \tkzDefMidPoint(p7,p1)\tkzGetPoint{mid71}
      \tkzDefMidPoint(c235,c123)\tkzGetPoint{vs}
      \tkzDefMidPoint(p2,p3)\tkzGetPoint{mid23}
      \tkzDefMidPoint(p2,p4)\tkzGetPoint{mid24}
      % intersections
      \tkzInterLL(p1,p3)(c271,c123)\tkzGetPoint{i1}
      \tkzInterLL(p1,p3)(c123,c235)\tkzGetPoint{i2}
      % voronoi volume
      \tkzFillPolygon[fill = gray!10](p1,p2,p3,p4)
      \tkzFillPolygon[fill = gray!15](c123,c235,c256,c267,c271)
      \tkzDrawPolygon[thin](c123,c235,c256,c267,c271)
      \tkzDrawPolygon[dotted](p1,p2,p3)
      \tkzDrawPolygon[thin, ](c123,c235,c256,c267,c271)
      \tkzFillPolygon[fill = gray!25](p2,mid12,i1,i2,mid23)
      \tkzDrawPolygon[very thick](p2,mid12,i1,i2,mid23)
      \tkzDrawLines[thin, add = 0 and 0](c123,c134)
      \tkzDrawLines[thin, add = 0 and 1/2](c134,mid14)
      \tkzDrawLines[thin, add = 0 and 1/3](c134,mid34)
      \tkzDrawLines[thin, add = 0 and 1/2](c235,mid35)
      \tkzDrawLines[thin, add = 0 and 2/3](c256,mid56)
      \tkzDrawLines[thin, add = 0 and 1/2](c267,mid67)
      \tkzDrawLines[thin, add = 0 and 1/2](c271,mid71)
      % triangles 1-2-4 (dotted)
      \tkzDrawPolygon[dotted](p1,p4,p3,p5,p6,p7)
      \tkzDrawLines[  dotted, add = 0 and 0](p2,p5)
      \tkzDrawLines[  dotted, add = 0 and 0](p2,p6)
      \tkzDrawLines[  dotted, add = 0 and 0](p2,p7)
      % draw vertices
      \tkzDrawPoints(p1,p2,p3,p4,p5,p6,p7)
      \tkzLabelPoint[right](labelVV){$V_4$}
      \tkzLabelPoint[right](labelVS){$V_{4j}$}
      \tkzLabelPoint[right](labelEj){$\bm{E}^2_j$}
      %\tkzDrawSegment[add = 0 and 1/10, arrows = ->](labelVS, vs)
      %\tkzLabelPoint[left = 1mm](p2){$\bm{v}_i$}
      \tkzLabelPoint[below](p1){$\bm{v}_1$}
      \tkzLabelPoint[left = 1mm](p2){$\bm{v}_4$}
      \tkzLabelPoint[above right](p3){$\bm{v}_3$}
      \tkzLabelPoint[below left](p4){$\bm{v}_2$}
      %\tkzLabelPoint[above](p5){$\bm{v}_5$}
      %\tkzLabelPoint[left](p6){$\bm{v}_6$}
      %\tkzLabelPoint[below left](p7){$\bm{v}_7$}
   \end{tikzpicture}%
   \caption{%
      Объем  $V_4$, поверхность $\partial V_4$ (слева), и элемент
      объема (ячейка) Вороного $V_{4j} = V_4 \cap \bm{E}^2_j$ (справа)%
   }\label{fig:voronoi}%
\end{figure}

\emph{Объем Вороного} $V_i$ вершины $\bm{v}_i$~--- это часть $\mathbb{R}^d$, состоящая из всех точек, для которых вершина $\bm{v}_i$ ближайшая из всех вершин сетки:
\[
   V_i := \{ \bm{x} \in \mathbb{R}^d :
      \|\bm{x} - \bm{v}_i\| < \|\bm{x} - \bm{v}_j \|
      ~
      \forall \bm{v}_j \in \Th(\Omega) \}
      .
 \]
Граница $\partial V_i$ объема Вороного называется соответствующей
\emph{поверхностью Вороного} (\cref{fig:voronoi}, слева).

Разбиение Вороного для заданного набора вершин является двойственной
структурой по отношению к триангуляции Делоне.  Каждая грань Вороного
соответствует ребру триангуляции Делоне.

Пересечение объёма Вороного $V_i$ с симплексом $\bm{E}^d_j$, $V_{ij}
= V_i \cap \bm{E}^d_j$, является \emph{ячейкой или элементом объема
Вороного} вершины $\bm{v}_i$ по отношению к $\bm{E}^d_j$
(\cref{fig:voronoi}, справа).  Будучи пересечением полупространств,
ячейки и грани Вороного всегда выпуклые, и их меры являются
непрерывными функциями координат вершин.

%--- BCD ----------------------------------------
\subsection{Триангуляции (сетки) Делоне, согласованные с границей}

Введение граничных или межфазных условий в сетку может нарушить свойство Делоне.
Например, рассмотрим однородные граничные условия Неймана для границы области (внутренней или внешней) \cref{fig:voronoi} вдоль ребер $\bm{v}_1 \bm{v}_3$ и $\bm{v}_3 \bm{v}_2$ (\cref{fig:bcd} (слева)).
В схемах, в которых величины хранятся в ячейках, распространенной техникой для граничных условий Неймана является введение дополнительной виртуальной вершины $\bm{v}_4^*$ как зеркального отображения $\bm{v}_4$ через $\bm{v}_1 \bm{v}_3$.
Если виртуальная вершина $\bm{v}_4^*$ находится внутри описанного круга для треугольника $\triangle \bm{v}_1 \bm{v}_3 \bm{v}_4$ (что произойдет тогда и только тогда, если $\bm{v}_4$ находится внутри описанного круга для $\bm{v}_1 \bm{v}_3$), то треугольник $\triangle \bm{v}_1 \bm{v}_4 \bm{v}_3$ нарушит свойство Делоне.
Эта проблема присутствует и в более высоких измерениях и мотивирует следующее определение.

%--- Problem with Delaunay and interfaces ------------------------------
\begin{figure}[t]%
   \begin{tikzpicture}[scale = 0.85]%
      \useasboundingbox (-1.23, -1.07) rectangle (4.78, 4.27);
      % define vertices
      \tkzDefPoint(0.0,  0.0){p1}
      \tkzDefPoint(0.1,  2.5){p2}
      \tkzDefPoint(2.8,  3.2){p3}
      \tkzDefPoint(4.5, -0.3){p4}
      \begin{pgfinterruptboundingbox}
      \tkzDefPointBy[reflection=over p1--p3](p2)\tkzGetPoint{p22}
      \tkzCircumCenter(p1,p2,p3) \tkzGetPoint{c}
      \end{pgfinterruptboundingbox}
      \tkzInterLL(p1,p3)(p2,p22) \tkzGetPoint{i13}
      % triangles, boundary segments
      \tkzDrawPolygon[thin, fill = gray!15](p1,p2,p3)
      \tkzDrawLines[very thick, add = 0 and 1/7](p3,p1)
      \tkzDrawLines[very thick, add = 0 and 1/7](p3,p4)
      % reflected point
      \tkzDrawLines[thin, dashed,  add = 0 and 0](p2,p22)
      \tkzDrawSegments[thin, dashed](p1,p22 p22,p3)
      \tkzDrawSegments[thin, dashed](p2,i13 i13,p22)
      \tkzMarkSegments[mark = ||](p2,i13 p22,i13)
      \tkzMarkRightAngle(p2,i13,p3)
      % circumcircle
      \tkzDrawCircle[dotted](c,p1)
      % draw vertices
      \tkzDrawPoints(p1,p2,p3,p4,p22)
      \tkzLabelPoint[left = 0.07](p1){$\bm{v}_1$}
      \tkzLabelPoint[above left](p2){$\bm{v}_4$}
      \tkzLabelPoint[above](p3){$\bm{v}_3$}
      \tkzLabelPoint[below left](p4){$\bm{v}_2$}
      \tkzLabelPoint[right](p22){$\bm{v}_4^*$}
      \end{tikzpicture}%
   \hfill{}%
   \begin{tikzpicture}[scale = 0.85]%
      \useasboundingbox (-1.23, -1.07) rectangle (4.78, 4.27);
      % vertices
      \tkzDefPoint(0.0,  0.0){p1}
      \tkzDefPoint(0.1,  2.5){p2}
      \tkzDefPoint(2.8,  3.2){p3}
      \tkzDefPoint(4.5, -0.3){p4}
      % circumceters of 1-2-3 & 1-3-4
      \begin{pgfinterruptboundingbox}
      \tkzCircumCenter(p1,p2,p3)\tkzGetPoint{c123}
      \tkzCircumCenter(p1,p3,p4)\tkzGetPoint{c134}
      \tkzDefPointBy[projection=onto p1--p3](p2)\tkzGetPoint{p22}
      \end{pgfinterruptboundingbox}
      % edge midpoints
      \tkzDefMidPoint(p1,p2)\tkzGetPoint{mid12}
      \tkzDefMidPoint(p1,p3)\tkzGetPoint{mid13}
      \tkzDefMidPoint(p2,p3)\tkzGetPoint{mid23}
      \tkzDefMidPoint(p3,p4)\tkzGetPoint{mid34}
      \tkzDefMidPoint(p1,p4)\tkzGetPoint{mid14}
      % tdraw riangles and boundary segments
      \tkzDrawPolygon[thin, fill = gray!15](p1,p2,p3)
      \tkzDrawLines[very thick, add = 0 and 1/7](p3,p1)
      \tkzDrawLines[very thick, add = 0 and 1/7](p3,p4)
      % circumball of 1-3
      \tkzDrawCircle[dotted](mid12,p2)
      \tkzDrawCircle[dotted](mid23,p2)
      % admissible vertex
      \tkzDrawLines[dashed,  add = 0 and 0](p2,p22)
      \tkzMarkRightAngle(p2,p22,p3)
      % draw Voronoi volumes
      \tkzDrawLines[thin, add = 0 and 0](c123,c134)
      \tkzDrawLines[thin, add = 0 and 1/2](c134,mid14)
      \tkzDrawLines[thin, add = 0 and 1/2](c134,mid34)
      \tkzDrawLines[thin, add = 0 and 1/3](c123,mid12)
      \tkzDrawLines[thin, add = 0 and 1/3](c123,mid23)
      % draw vertices
      \tkzDrawPoints(p1,p2,p3,p4)
      \tkzLabelPoint[left = 0.07](p1){$\bm{v}_1$}
      \tkzLabelPoint[above left](p2){$\bm{v}_4$}
      \tkzLabelPoint[above](p3){$\bm{v}_3$}
      \tkzLabelPoint[below left](p4){$\bm{v}_2$}
   \end{tikzpicture}%
   \hfill{}%
   \begin{tikzpicture}[scale = 0.85]%
      \useasboundingbox (-1.23, -1.07) rectangle (4.78, 4.27);
      % vertices
      \tkzDefPoint(0.0,  0.0){p1}
      \tkzDefPoint(0.1,  2.5){p2}
      \tkzDefPoint(2.8,  3.2){p3}
      \tkzDefPoint(4.5, -0.3){p4}
      % edge midpoints and projection
      \tkzDefMidPoint(p1,p2)\tkzGetPoint{mid12}
      \tkzDefMidPoint(p2,p3)\tkzGetPoint{mid23}
      \tkzDefMidPoint(p3,p4)\tkzGetPoint{mid34}
      % circumceters of 1-2-3 & 1-3-4
      \begin{pgfinterruptboundingbox}
      \tkzDefPointWith[orthogonal](mid12,p2)\tkzGetPoint{h12}
      \tkzDefPointWith[orthogonal](mid23,p3)\tkzGetPoint{h23}
      \tkzDefPointWith[orthogonal](mid34,p4)\tkzGetPoint{h34}
      \tkzDefPointBy[projection=onto p1--p3](p2)\tkzGetPoint{p22}
      \tkzDefPointBy[projection=onto p1--p3](mid12)\tkzGetPoint{h122}
      \tkzDefPointBy[projection=onto p1--p3](mid23)\tkzGetPoint{h232}
      \end{pgfinterruptboundingbox}
      % tdraw riangles and boundary segments
      \tkzDrawPolygon[thin, fill = gray!15](p1,p2,p3)
      \tkzDrawLines[         add = 0 and 0](p2,p22)
      \tkzDrawLines[very thick, add = 0 and 1/2](h122,p1)
      \tkzDrawLines[very thick, add = 0 and 0](h232,p3)
      \tkzDrawLines[very thick, dashed, add = 0 and 0](h122,h232)
      \tkzDrawLines[very thick, add = 0 and 1/7](p3,p4)
      % circumballs
      \tkzDrawCircle[dotted](h122,p1)
      \tkzDrawCircle[dotted](h232,p3)
      % draw Voronoi volumes
      \tkzDrawLines[thin, add = 0 and -0.5](mid12,h12)
      \tkzDrawLines[thin, add = 0 and -0.6](mid23,h23)
      \tkzDrawLines[thin, add = 0 and -0.5](mid34,h34)
      \tkzDrawLines[thin, add = 0 and  0  ](mid12,mid23)
      \tkzDrawLines[thin, add = 0 and  0  ](mid12,h122)
      \tkzDrawLines[thin, add = 0 and  0  ](mid23,h232)
      \tkzMarkRightAngle(p2,p22,p3)
      % draw vertices
      \tkzDrawPoints(p1,p2,p3,p4,p22)
      \tkzLabelPoint[left = 0.07](p1){$\bm{v}_1$}
      \tkzLabelPoint[above left](p2){$\bm{v}_4$}
      \tkzLabelPoint[above](p3){$\bm{v}_3$}
      \tkzLabelPoint[below left](p4){$\bm{v}_2$}
      \tkzLabelPoint[above right = -2mm and 1mm](p22){$\tilde{\bm{v}}_4$}
   \end{tikzpicture}%
   \caption{%
      Если угол $\angle \bm{v}_1 \bm{v}_4 \bm{v}_3$ тупой, то
      виртуальная вершина $\bm{v}_4^*$ для граничного условия Неймана
      нарушит условие Делоне (слева).
      Вставка новой вершины как ортогональной проекции $\bm{v}_4$
      на $\bm{v}_1 \bm{v}_3$  восстанавливает граничное соответствие сетки
      и не нарушает свойства Делоне ребер $\bm{v}_1 \bm{v}_4$ и $\bm{v}_4
      \bm{v}_3$  (в центре);  жирная  пунктирная линия соответствует
      граничному условию, имеющему отношение к $\tilde{\bm{v}}_4$%
   }\label{fig:bcd}%
\end{figure}

Сетка Делоне называется \emph{сеткой Делоне, согласованной с границей}, если ни одна вершина сетки не находится внутри самой малой (диаметральной, экваториальной) сферы какого угодно маломерного симплекса границы или внутреннего интерфейса:
\[
   \forall \bm{v}_i \in \Th(\Omega),
   \forall \bm{E}^m_j \in \partial \Omega_i,
   ~
   \bm{v}_i \not\in \mathrm{B^o}(\bm{E}^m_j)
   \quad(1 \le m < d)
   .
\]
Свойство пустого диаметрального шара иногда называют свойством \emph{Габриэля}~\cite{GabSok69}.
Иными словами, сетка Делоне, согласованная с границей, представляет собой сетку, где каждый симплекс удовлетворяет критерию Делоне, а каждый граничный симплекс~--- критерию Габриэля.

В ситуации \cref{fig:bcd} (слева), согласование с границей е можно достичь,  разделив $\triangle \bm{v}_1 \bm{v}_3 \bm{v}_4$ новой
вершиной на граничном ребре $\bm{v}_1 \bm{v}_3$.
Пересечение диаметральных дисков $\bm{v}_1 \bm{v}_4$
и $\bm{v}_4 \bm{v}_3$ с ребром границы $\bm{v}_1 \bm{v}_3$ соответствует проекции $\bm{v}_4$ на $\bm{v}_1 \bm{v}_3$
(\cref{fig:bcd}, в центре) и гарантирует восстановлении согласования сетки  с границей (\cref{fig:bcd}, справа) так, что при этом что все остальные симплексы сетки сохранят свойство Делоне. 

Подробности о сетках Делоне, согласованных с границей, можно найти в  работе~\cite{SiGarFuh10} и в ссылках в ней.

%--- Finite Volumes ----------------------------------------
\section{Основные свойства метода конечных объемов с использованием ячеек Вороного (VFVM)}

Если применить теорему Гаусса к ячейке Вороного или ее пересечению
с $\partial \Omega$, то каждая грань двойственна  ребру и, следовательно, ее можно представить как
набор симплексов.
Более удобно ассемблировать уравнения конечных объемов с помощью симплексов двойственной сетки Делоне, в частности,
гораздо удобнее определять свойства материалов на симплексах, чем на объемах Вороного.

\subsection{Уравнения конечных объемов с использованием ячеек Вороного}

%--- Differences along edges ----------------------------------------
\emph{Матрицы смежности} $\tilde G_i$ для $i$-мерного симплекса,
отображающая вершины на ребра может быть определена рекурсивно через
\[
   \tilde{G}_{i+1} =
   \begin{pmatrix*}[l]
     \phantom{-}\tilde{G}_i & \bm{0}_{\frac{i(i-1)}{2}} \\
      -\mathbb{I}_{i+1}       & \bm{1}_{i+1}
   \end{pmatrix*},
   \quad i = 1, 2, \dotsc,
\]
где $\mathbb{I}_j$ обозначает  единичную матрицу $j \times j$,
$\bm{0}_j = ( \overbrace{ 0,\dotsc, 0}^{j}{)}^\mathrm{T}$,
$\bm{1}_j = ( \overbrace{1, \dotsc, 1}^{j}{)}^\mathrm{T}$,
а $\tilde{G}_1 = (-1, 1)$ определяет \emph{разность на интервале},
\[
   \tilde{G}_1 \begin{pmatrix} \bm{x}_1\\\bm{x}_2 \end{pmatrix}
   =
      \begin{pmatrix} -1 & 1 \end{pmatrix}
      \begin{pmatrix} \bm{v}_1\\\bm{v}_2 \end{pmatrix}
      = \bm{v}_2 - \bm{v}_1
      .
\]
Получаем
\begin{align*}
   \tilde{G}_2 &= 
   \begin{pmatrix*}[r]
      \tilde{G}_1     & \phantom{-}\bm{0}_1 \\
      -\mathbb{I}_2  &            \bm{1}_2
   \end{pmatrix*}
   =
   \left(\begin{array}{rr|r}
      -1  &  1  & \phantom{-}0 \\
      \hline{}
      -1  &  0  &            1 \\
       0  & -1  &            1 \\
   \end{array}\right),
   \\
   \tilde{G}_3 &=
   \begin{pmatrix*}[r]
      \tilde{G}_2     & \phantom{-}\bm{0}_3 \\
      -\mathbb{I}_3  &            \bm{1}_3
   \end{pmatrix*}
   =
   \left(\begin{array}{rrr|r}
     -1  &  1  &  0  &  \phantom{-}0\\
     -1  &  0  &  1  &             0\\
      0  & -1  &  1  &             0\\
      \hline{}
     -1  &  0  &  0  &  1\\
      0  & -1  &  0  &  1\\
      0  &  0  & -1  &  1\\
   \end{array}\right),
   \\
   \tilde{G}_4 &=
   \begin{pmatrix*}[r]
      \tilde{G}_3     & \phantom{-}\bm{0}_6 \\
      -\mathbb{I}_4  &            \bm{1}_4
   \end{pmatrix*},
   \quad \dotso{}
\end{align*}

Матрица $\tilde G_d$ отображает вершины на ребра.
Соответственно, $\tilde G_d^\mathrm{T}$ отображает ребра на вершины, и 
$\tilde G_d^\mathrm{T} \tilde G_d$ является \emph{матрицей Кирхгофа (Лапласианом графа)}
$d$-мерного симплекса.

Обозначим длину ребра $\bm{E}_i^1$ как $\lvert
\tilde G_1 \bm{E}{_i^1} \rvert$ и диагональную матрицу длин всех ребер симплекса  $\bm{E}^d_i$ как $[ \,
\lvert \tilde G_d \bm{E}^d_i \rvert \, ]$.
Тогда
\[
   {\left[ \lvert \tilde G_d \bm{E}^d_i \rvert \right]}^{-1}
      \tilde G_d \bm{u}
\]
является приближением значения $\nabla u$ в направлении $\partial u/\partial
n_{k}$, перпендикулярном грани ячейки Вороного, относительно
симплекса $\bm{E}^d_i$, которое преобразовывает заданные линейные функции
на каждом ребре в корректную постоянную.
Далее,
\[
   \tilde G_d^T  [ \sigma_{\bm{E}_i^d} ]
      = \tilde G_d^\mathrm{T} {\left[ \lvert \tilde G_d \bm{E}^d_i \rvert \right]}^{-1} 
      {\left\{
         \left[ \lvert \tilde G_d \bm{E}^d_i \rvert \right] 
         \left[ \sigma_{\bm{E}^d_i} \right]
      \right\}}^{\frac{1}{2}}
   \]
является симметричной формой отображения с ребер на вершины, включающей геометрический вес грани Вороного, пересекаемой ребром, умноженным на длину ребра.
Следовательно, метрический фактор $\{\cdot\}$ имеет смысл удвоенного детерминанта 
(длина ребра в два раза больше высоты $d$-мерного симплекса, образованного 
пересечением симплекса и $\sigma_i$), а не объема,
пересеченного с соответствующей частью каждого ребра симплекса $\bm{E}^d_i$.
Это правило надо последовательно соблюдать во всех частях оператора и граничных условий.
Следовательно,
\[
   \varepsilon \int_{\bm{E}^d_i} \nabla \cdot \nabla u \, \D{V}
   = \oint_\Sigma \nabla u \, \D{\bm{S}}
   \approx \varepsilon \, \tilde G_d^\mathrm{T} [\sigma_i] \,
   {\left[ \lvert \tilde G_d \bm{E}^d_i \rvert \right]}^{-1} \,
   \tilde G_d \bm{u}  := \varepsilon \,
   G_d^\mathrm{T} G_d \bm{u}
\]
является приближением \emph{оператора Лапласа, проинтегрированного по объему}
на $\bm{E}^d_i$, где $\varepsilon = const$ on $\bm{E}^d_i$
и $\sigma_j \ge 0$ для каждого ребра $j$.

Граничное условие
\[
   \alpha(u) u +\beta(u) \frac{\partial u}{\partial n} +\gamma(u)=0
\]
интегрируется на гранях Вороного, принадлежащих к поверхностным вершинам; грань определяет нормальное направление, которое совпадает с направлением ребра, трансверсального к границе, в случае сетки Делоне согласованной с границей:
\[
   \varepsilon \frac{\partial u}{\partial n} \partial V_{s,i}
   = - \partial V_{s,i} \frac{\varepsilon}{\beta(u)}
   \left( u_i \alpha(u_i) + \gamma(u_i) \right)
   .
\]
Нормальная производная вызывает изменение соответствующего диагонального элемента
и правой части дискретных уравнений. Условие Дирихле получаются для $\beta
\approx \epsilon_\mathrm{machine}^2$ и легко гомогенизируются.

Сетка Делоне, согласованная с границей, позволяет корректно разрезать ячейки Вороного вдоль границы, задаваемой симплексами, и ввести то же поверхностное интегрирование на $d-1$-мерных ячейках Вороного.
В случае пограничных слоев, вершины на границе могут быть носителями экстремумов, введенных граничным условием.
Это будет не так, если для приближения подобласти $\Omega_i^d$ использовать ячейки Вороного.

Основным условием для метода конечных объемов является справедливость теоремы Гаусса.
Наибольшее функциональное пространство, где действует теорема Гаусса~--- это пространство функций ограниченной вариации.
Это очень слабое условие на класс функций.

Если допускать произвольный скачок коэффициента диффузии между смежными областями, то может возникнуть необходимость ввести дополнительные ограничения на углы сетки около интерфейса.

%--- Compensation effect ------------------------------
\subsection{Механизм компенсации и свойство S-матрицы}

Основное преимущество метода конечных объемов с использованием ячеек
Вороного (VFVM) состоит в том, что матрица системы  является невырожденной
\emph{M-матрицей}.  \emph{M-матрица} (матрица Минковского)~--- это
Z-матрица (матрица с неположительными внедиагональными элементами)
с положительной диагональю, сумма всех столбцов которой неотрицательна.
Если хотя бы одна столбцовая сумма М-матрицы положительна, то эта
матрица обратима и обратная матрица является положительной.
Симметричная М-матрица называется \emph{S-матрицей} (матрица
Стильтеса).

\begin{figure}[t]%
   \centering{}%
   \includegraphics[width=0.3\textwidth,clip]{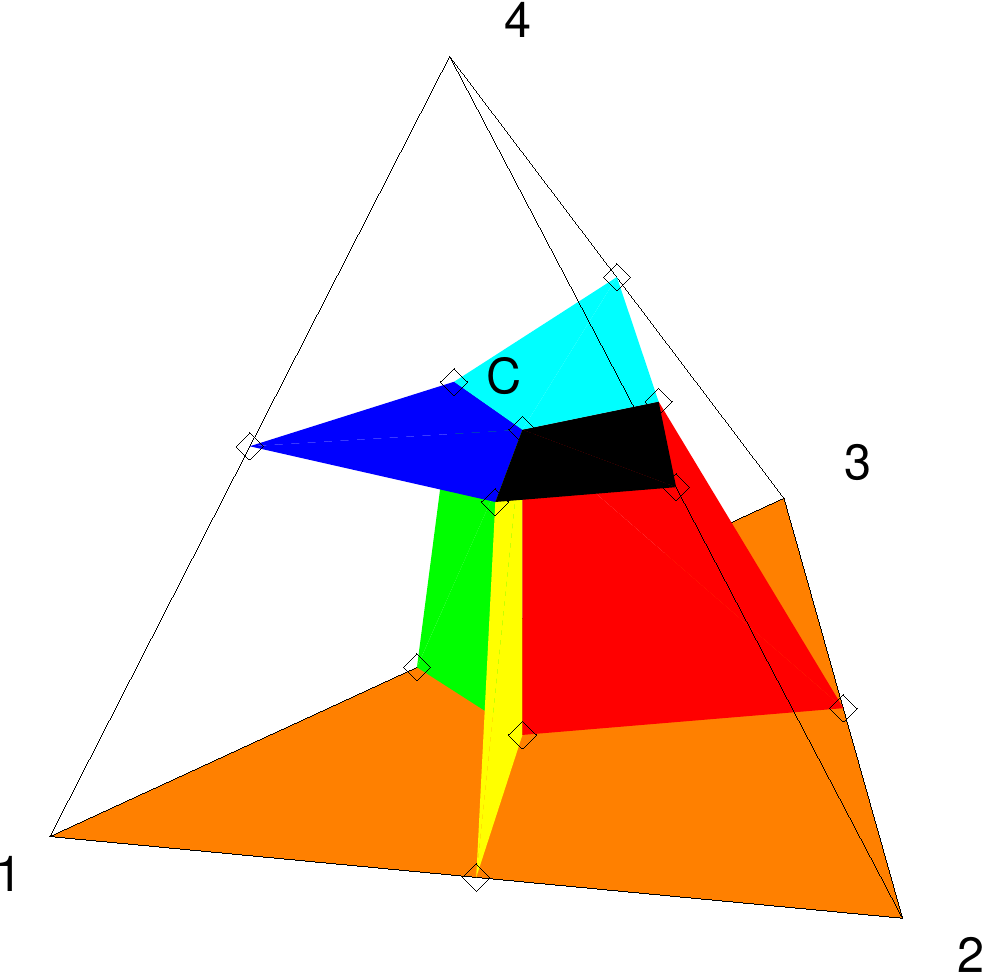}%
   \hspace{0.1\linewidth}%
   \includegraphics[width=0.3\textwidth,clip]{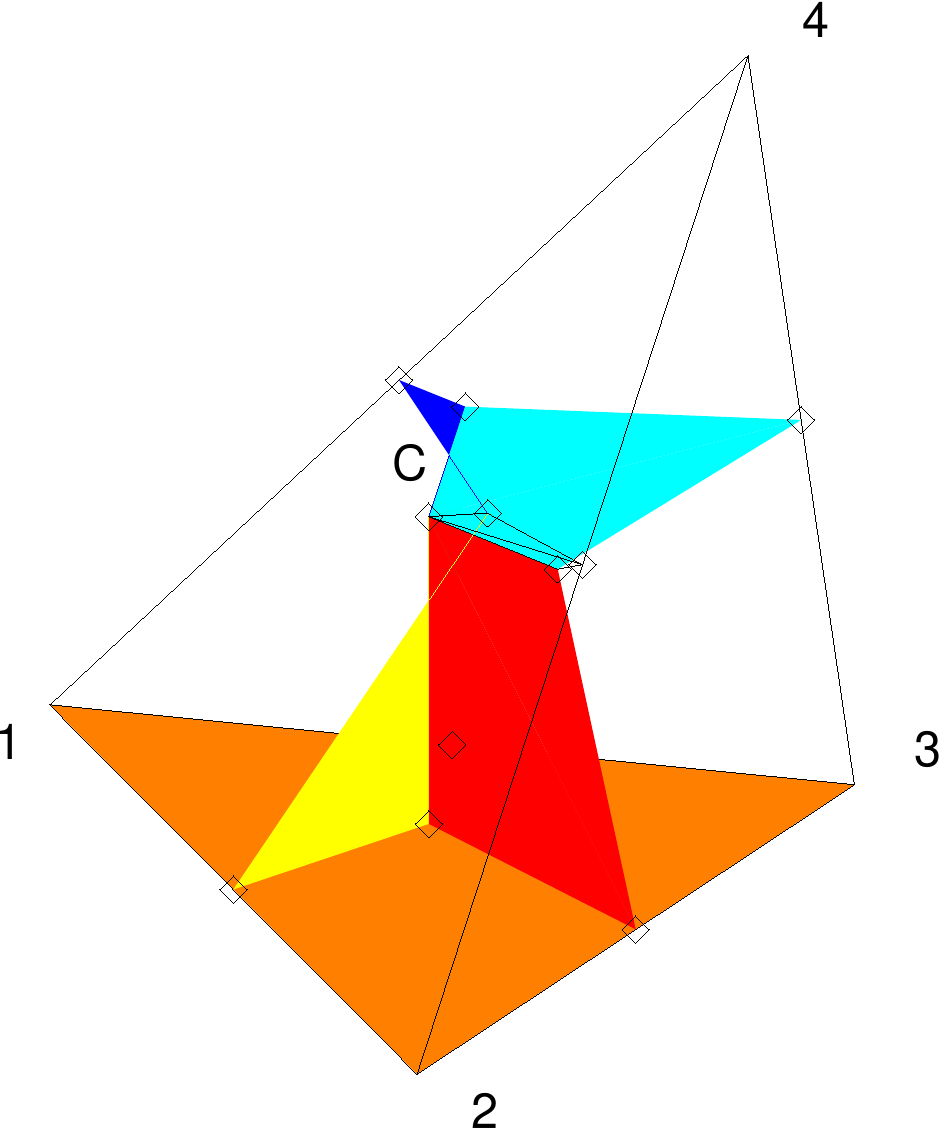}%
   \hfill{}%
   \caption{%
      Слева: все центры окружностей находятся внутри их определяющих
      $d$-симплексов, $d = 1, 2, 3$.
      Справа: центр $C$ окружности тетраэдра  и центр окружности
      треугольника $\triangle 234$ выходят за пределы их определяющих
      симплексов%
      \label{fig:tet:1}%
      \label{fig:tet:2}%
   }%
\end{figure}

Рассмотрим случай $\mathbb{R}^3$ и два допустимых тетраэдра
Делоне (\cref{fig:tet:1,fig:tet:2}).
Для удобства изложения скажем, что свободная четвертая вершина,
противолежащая каждой треугольной грани определяет четыре условия
Делоне.  На~\cref{fig:tet:1} (слева) все центры описанных окружностей
(далее $C$ с индексом соответствующего, возможно низкоразмерного
симплекса) находятся внутри соответствующего определяющего симплекса
(соответственно  одно-, двух- или трехмерного).

На~\cref{fig:tet:2} (справа) центр $C_{1234}$ окружности тетраэдра и
центр $C_{234}$ окружности треугольника $\triangle 234$ находятся вне их
определяющих симплексов.
Каждая из граней Вороного определяется двумя треугольниками, например:
грань, определяемая вершинами 1 и 2, состоит из $\triangle C_{12}
C_{123} C_{1234}$ и $\triangle C_{21} C_{124} C_{1234}$, которые
перпендикулярны ребру $12$ и изменяют знак в соответствии с положением
соответствующих центров окружности треугольника и тетраэдра и имеют
положительную проекцию на ребра 12 и 21 для внутренних центров
окружностей.
Соседний тетраэдр, который тоже имеет $\triangle 124$ в качестве
грани, имеет те же центры описанных окружностей общих ребер и то
же направление нормали к $\triangle 124$.
Грань Вороного ограничена полигоном, соединяющим все центры окружностей
тетраэдров, разделяющих одно ребро и расположенных в плоскости,
перпендикулярной этому ребру.
Суммирование всех вкладов тетраэдров, разделяющих общее ребро, дает
площадь грани Вороного (со знаком), относящуюся к этому ребру.
Условие Делоне гарантирует, что расстояние между центрами окружности
двух соседних тетраэдров является неотрицательным, Таким образом,
значение площади грани Вороного неотрицательно.

Этот эффект называется \emph{компенсацией}: отрицательный вклад
компенсируется положительным вкладом соседнего тетраэдра.  Он
присутствует в любом измерении и гарантирует, что матрица в методе конечного
объема является S-матрицей или M-матрицей.

Если на общей $d$-сфере лежат более $m \ge d + 2$ вершин, то центры
$m- d$ окружностей будут совпадать.  Сетка Делоне с ограничениями,
гарантирует, что все центры окружности находятся внутри соответствующих
определяющих симплексов, части ячеек и поверхностей Делоне являются
неотрицательными, и сохраняются ориентированные объем подобласти
$\Omega_i$ и площади поверхностей.  Это позволяет легко справляться
со всеми граничными условиями третьего рода.

Для P\textsubscript{1} (линейного) метода конечных элементов (МКЭ)
в $d$-мерном пространстве, матрица Дирихле является проекцией каждой
из $d-1$-мерной грани на все остальные,
\[
   P^\mathrm{-T} P^{-1},
   \quad P = \bm{v}_j - \bm{v}_0,
   \quad j = 1, \dotsc, d
   .
\]
В двух измерениях она идентична матрице $G_d^\mathrm{T} G_d$ для VFVM.

В трех- и более измерениях это уже не так (за исключением особых
случаев) и компенсация  негативных частей работает по-другому, так
что граничные сетки Делоне не являются ни достаточными, ни необходимыми
для получения S-матрицы для МКЭ (подробнее
см.~\cite{BanRos87,Ker96,XuZik99,KoFHPS00} и ссылки в них).

Так называемое \emph{условие остроугольности} (все углы между граничными нормалями являются острыми), гарантируещее
свойство S-матрицы для матрицы Дирихле для МКЭ, является слишком
ограничительным и сложным для использования на практике.

%--- Dissipativity ----------------------------------------
\subsection{Принцип максимума и диссипативность}

Решение VFVM удовлетворяет  \emph{слабому дискретному принципу максимума}.
Пусть $\bm{u}$ удовлетворяет $G_d^d G_d \bm{u}= \bm{0}$
с заданным граничным условием Дирихле $\bm{u}_{D}$, $\hat u=\max u_{D}$,
$\check u=\min u_{D}$. Предположим, что $\max u = U > \hat u$.
Проверка уравнения с помощью $(\bm{u}-\hat u)^\mathrm{+T}$ дает
$(\bm{u}-\hat u)^\mathrm{+T} G^\mathrm{T}_d G_d \bm{u} > 0$,
потому что $G_d \bm{u}$ и $G_d {(\bm{u}-\hat u)}^{+}$
имеют либо одинаковый знак (и, как минимум, одно значение отлично от нуля) или
$G_n {(\bm{u}-\hat u)}^{+} = \bm{0}$, что приводит к противоречию.

Диссипативность вытекает из этого же аргумента. В дискретном случае принцип максимума и диссипативность в этом смысле эквивалентны.

Замечание и диссипативности членов нулевого порядка: объёмы Вороного образуют диагональную матрицу масс $[V]$, больший носитель
невозможно получить  без допущения малости вариации решения на носителе.

%--- Non-Delaunay ----------------------------------------
\subsection{Невыполнения условия Делоне для оператора Лапласа}%
\label{sec:nonDelaunayLaplace}

Так как ячейки Вороного имеют неотрицательные площади поверхности,
ассемблирование
\[
   A = G_d^\mathrm{T} [\varepsilon] G_d
\]
приводит к S-матрице: $A_{ii} > 0$, $A_{ij} \le 0$, $i\ne j$, и
${\bm 1}^\mathrm{T}_d G_d^\mathrm{T} = \bm{0}_d$  (сумма в столбцах
и строках равна нулю для вершин без граничных условий).
Очевидно, что матрица $A^{-1}$ поэлементно положительна для неприводимой $A$
или связных сеток.  Конечное значение $\varepsilon$ и граничные условия
добавляют только положительные значения к диагонали по крайней мере
для одной вершины.  Следовательно, решение ограничено значениями
Дирихле и положительно для нулевых значений Дирихле и положительной
правой стороны уравнения.

Если заданная сетка не является сеткой Делоне, суммирование
ориентированной площади поверхностей для каждого ребра может привести
к образованию отрицательных поверхностей вокруг этого ребра. Это
может быть произвольное ребро между двумя вершинами, не соответствующими
грани Вороного (описанная сфера одного симплекса содержит центр
другой).  Далее следует $A_{ij}>0$ и позитивность больше не
гарантируется, по крайней мере, можно случиться, что появятся
локальные экстремумы в вершинах $i$ и $j$ для значительных отклонений
от условия Делоне. По построению, $G_d \bm{1}_d =\bm{0}_d$.  Введение
$A_{ij}=0$ для соответствия с $A_\mathrm{ii}$ эквивалентно введению
локального пограничного условия Неймана и несогласованной
дискретизации.

\begin{figure}[t]%
   \centering{}%
   \begin{tikzpicture}[scale = 0.8]%
      % define vertices 1, 2, 3, 4
      \tkzDefPoint( 0.0,  0.0){pi}
      \tkzDefPoint( 3.0,  0.0){p1}
      \tkzDefPoint( 1.5,  1.3){p2}
      \tkzDefPoint(-1.5,  1.3){p3}
      \tkzDefPoint(-1.5, -1.3){p4}
      \tkzDefPoint( 1.5, -1.3){pk}
      % triangles 1-2-4 (dotted)
      \tkzDrawArc[dashed](pi,p2)(pk)
      \tkzDrawSegment[very thick](pi,p1)
      \tkzDrawSegments[thin](pi,p2 pi,p3 pi,p4 pi,pk)
      % draw vertices
      \tkzDrawPoints(pi,p1,p2,p3,p4,pk)
      \tkzLabelPoint[above = 1mm](pi){$\bm{v}_i$}
      \tkzLabelPoint[above](p1){$\bm{v}_1$}
      \tkzLabelPoint[right](p2){$\bm{v}_2$}
      \tkzLabelPoint[left](p3){$\bm{v}_3$}
      \tkzLabelPoint[left](p4){$\bm{v}_4$}
      \tkzLabelPoint[right](pk){$\bm{v}_{k}$}
   \end{tikzpicture}%
   \hspace{0.1\linewidth}%
   \begin{tikzpicture}[scale = 0.8]%
      % define vertices 1, 2, 3, 4
      \tkzDefPoint( 0.0,  0.0){pi}
      \tkzDefPoint( 1.0,  3.0){p1}
      \tkzDefPoint( 2.0,  3.0){p11}
      \tkzDefPoint(-2.0, -0.1){p3}
      % triangles 1-2-4 (dotted)
      \tkzDrawLine[thin, add = 1/3 and 0](p3,pi)
      \tkzDrawSegments[thin](pi,p1 p1,p11)
      % draw vertices
      \tkzDrawPoints(pi,p1,p3)
      \tkzLabelPoint[below](pi){$\bm{v}_i$}
      \tkzLabelPoint[above right = 0 and 1mm](pi){$u_i=\varepsilon_i$}
      \tkzLabelPoint[below right](p1){$\bm{v}_1$}
      \tkzLabelPoint[above](p1){$u_1 = 1 - \varepsilon$}
      \tkzLabelPoint[below](p3){$\bm{v}_3$}
      \tkzLabelPoint[above](p3){$u_3 = \varepsilon_3$}
   \end{tikzpicture}%
   \caption{%
      Скачок решения вдоль ребра  $\bm{v}_i \bm{v}_1$
      от $u(\bm{v}_i) = \varepsilon_i$
      до $u(\bm{v}_1) = 1 - \varepsilon$%
   }\label{fig:jump}
\end{figure}
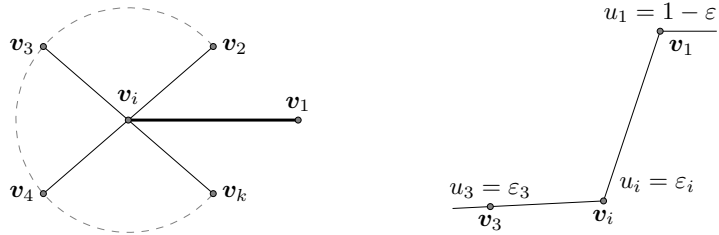

Предположим, что решение $u$ имеет неразрешенный движущийся внутренний слой c параметрами
\[
   u_1 = 1 - \varepsilon,
   \quad u_j = \varepsilon_j,
   \quad \text{с}~ 0< \varepsilon_j \le \varepsilon \ll 1,
   \quad j = i, 2, \dotsc, k,
\]
задающими порядок изменения $u$ за пределами скачка (\cref{fig:jump}). Здесь индексы~--- это просто номера точек сетки, в которых задана скалярная функция $u$.
Пусть $A_\text{non-Del}$ отвечает типу $\tilde G^\mathrm{T} [W]
\tilde G$ a сумма берется по всем симплексам относящимся к вершине $i$
(знак суммирования не указан, порядок сложения по ребрам или симплексам заданного
типа может быть выбран соответственно каждой цели, и контекст
определяет локальное или глобальное значение).

Пусть $l_{ij}$ будет длиной ребра для всех соседей
$\bm{v}_j$ вершины $\bm{v}_i$ и
\[
   \begin{cases}
      W_{ij}   = \frac{\sigma_{ij}}{l_{ij}} = 1,
      & j= 2, \dotsc, k,\\
      W_{i1} = \frac{\sigma_{i1}}{l_{i1}} = -1,
   \end{cases}
\]
т.е.\ ребро $\bm{v}_i \bm{v}_1$ не удовлетворяет условию Делоне.
Таким образом, $A_\text{non-Del} {\bm u}|_i$ при вершине $i$ приводит к тому, что
\[
   A_\text{non-Del} {\bm u}|_i
   = \sum^k_{j=2} (u_i - u_j)
   -(u_i - u_1)
   = \sum^k_{j=2} (\varepsilon_i- \varepsilon_j)
   -(\varepsilon_i - (1 - \varepsilon))
  \approx 1 > 0
  ,
\]
в то время как в случае Делоне
\[
   A_\text{Del} {\bm u}|_i
   = \sum^k_{j=2} (u_i - u_j)
      + (u_i - u_1)
   = \sum^k_{j=2} (\varepsilon_i- \varepsilon_j)
      + (\varepsilon_i - (1 - \varepsilon))
   \approx -1 < 0
   .
\]
Таким образом, правая часть для $i$ велика по сравнению с правой частью
для вершин, находящихся на плато уровня $\varepsilon$ или $1-\varepsilon$.
Поскольку $A \bm{u}$ в методе Ньютона (или другом
методе линейной коррекции невязки) является частью функции $f$,
и $A_\text{Del}$ вносит свой вклад в матрицу Якоби ($A$), то
\[
   \delta \tilde u\vert_i
      = (-A_\text{Del}^{-1} f)\vert_i > 0
\]
и поэтому $u^+ = u + \delta u$ (следующее приближение $u_i$ в
методе Ньютона) увеличивается диффузионной частью матрицы Якоби 
до среднего значения ближайших соседей.

В случае неудовлетворения условия Делоне, коэффициент $-1$ и огромная
(по абсолютному значению) правая сторона эквивалентна сильному
источнику с неправильным знаком.  Кроме того, $A_\text{non-Del}$
может быть неограниченным.  Очевидно, постоянные $A$ не могут быть
причиной скачка, но другие составляющие могут создать скачки
в случае небольшой диффузии.  Небольшая $A_\text{non-Del}$ может
быть уже достаточной для того, чтобы нарушить предельные значения
$u$ и таким образом привести к серьезным проблемам, которые не
могут быть решены простыми ограничителями для $u$, не создавая при
этом новых проблем в зависимости от специфики задачи.

В более общем случае, этот пример соответствует локально отрицательному
коэффициенту диссипации и, соответственно,  \emph{локальному
производству свободной энергии} схемой дискретизации.  Концепция
убывания свободной энергии по траекториям будет нарушена  в целом
и, в зависимости от сетки, возможное устойчивое состояние будет
отличаться от решения с минимальной энергией, так как свободная
энергия, произведенная схемой, должна быть рассеяна в где-то другом
месте.

%--- Mass Matrix -------------------------------------------------
\subsection{Непрерывная зависимость от положения вершин}

Рассмотрим одномерную краевую задачу:
\begin{align*}
   \begin{cases}
      u(y)'' + c u(y) = 0, & y \in (0,1), \quad c = const,\\
      u(0) = 1,\\
      u(1) = 0.
   \end{cases}
\end{align*}
Эта задачу и ее дискретизацию на заданной одномерной сетке по $y$
можно расширить двумерную ортогональную сетку по $x$, используя произвольное значения
шага сетки и однородные граничные условия Неймана.  При этом
сохраняется профиль $u$ по $y$ (двухмерное решение будет инвариантным
в направлении по $x$).

Триангуляция Делоне для получившегося регулярного прямоугольного
распределения узлов не однозначна, потому что каждая из четырех
соседних вершин одной ячейки в направлениях $x$ и $y$ лежат на одной
окружности (\cref{fig:crisscross},слева: оба направления диагональных
ребер удовлетворяют условию Делоне).  В данном случае, однозначное
положение диагональных ребер может получиться только при использовании
специального правила, например, 'снизу слева~-- наверх справа'.
Использование точной арифметики для теста на выполнение условия
Делоне в этом случае не поможет, так как он все равно приведет к
произвольному выбору ребра в случае если четыре узла имеют точные
координаты и лежат на одной окружности или по построению должны
быть на одной окружности, но были сохранены с округлением в
представлении с плавающей точкой.  Например, поворот сетки в
\cref{fig:crisscross} на $\pi/4$ приведет к тому, что четыре узла,
лежавшие вначале на одной окружности, в новой сетке не будут больше
лежать на одной окружности, т.к. $\sin \pi/4 = \cos \pi/4 = \sqrt{2}/2$
не может быть представлен точно числами в арифметике с плавающей
запятой без округления.

Для VFVM выбор диагональных ребер неважен, потому что грань Вороного,
соответствующая диагональному ребру,  имеет площадь, равную нулю.
Поверхности и объемы Вороного масштабируются в соответствии с
размером шага по $x$. Таким образом, сохраняется неизменность решения
в направлении $x$.  Кроме того, небольшие искажения координат узлов
приведут к соответственно небольшому изменению площадей граней
(объемы ячеек и граней Вороного  являются непрерывными функциями
координат узлов).  То же справедливо и для матриц типа $G^\mathrm{T}
G$ в любом $d$-мерном пространстве  (например, матрицы масс и матрицы
Дирихле).  Следовательно, очень небольшое изменение положения узлов
сетки приведет к очень небольшому изменению дискретного решения.

\begin{figure}[t]
   \centering{}
   \begin{tikzpicture}[scale = 1.5]%
      % define vertices 1, 2, 3, 4
      \tkzDefPoint(0, 0){p00}
      \tkzDefPoint(0, 1){p01}
      \tkzDefPoint(0, 2){p02}
      \tkzDefPoint(1, 0){p10}
      \tkzDefPoint(1, 1){p11}
      \tkzDefPoint(1, 2){p12}
      \tkzDefPoint(2, 0){p20}
      \tkzDefPoint(2, 1){p21}
      \tkzDefPoint(2, 2){p22}
      \tkzDefPoint(3, 0){p30}
      \tkzDefPoint(3, 1){p31}
      \tkzDefPoint(3, 2){p32}
      \tkzDefPoint(4, 0){p40}
      \tkzDefPoint(4, 1){p41}
      \tkzDefPoint(4, 2){p42}
      % areas
      \tkzFillPolygon[gray!15](p00,p01,p12,p22,p21,p10)
      \tkzFillPolygon[gray!15](p20,p21,p32,p42,p41,p30)
      % circumball example
      \tkzDefMidPoint(p01,p12)\tkzGetPoint{c}
      \tkzDrawCircle[dotted](c,p01)
      % box + diagonals
      \tkzDrawSegments[thin](p00,p40 p40,p42 p42,p02 p02,p00) % box
      \tkzDrawSegments[thin](p01,p41 p10,p12 p20,p22 p30,p32) % grid
      \tkzDrawSegments[thin](p00,p11 p10,p21 p20,p31 p30,p41) % diagonals
      \tkzDrawSegments[thin](p01,p12 p11,p22 p21,p32 p31,p42) % diagonals
      \tkzDrawSegments[dashed, thin](p02,p11) % diagonals
      %% draw vertices
      \tkzDrawPoints(p11,p31)
      \tkzLabelPoint[below right = 0 and 0.15](p11){$\bm{v}_1$}
      \tkzLabelPoint[below right](p31){$\bm{v}_2$}
   \end{tikzpicture}%
   \hspace{0.1\linewidth}%
   \begin{tikzpicture}[scale = 1.5]%
      % define vertices 1, 2, 3, 4
      \tkzDefPoint(0, 0){p00}
      \tkzDefPoint(0, 1){p01}
      \tkzDefPoint(0, 2){p02}
      \tkzDefPoint(1, 0){p10}
      \tkzDefPoint(1, 1){p11}
      \tkzDefPoint(1, 2){p12}
      \tkzDefPoint(2, 0){p20}
      \tkzDefPoint(2, 1){p21}
      \tkzDefPoint(2, 2){p22}
      \tkzDefPoint(3, 0){p30}
      \tkzDefPoint(3, 1){p31}
      \tkzDefPoint(3, 2){p32}
      \tkzDefPoint(4, 0){p40}
      \tkzDefPoint(4, 1){p41}
      \tkzDefPoint(4, 2){p42}
      % areas
      \tkzFillPolygon[gray!15](p00,p20,p22,p02)
      \tkzFillPolygon[gray!15](p21,p32,p41,p30)
      % box + diagonals
      \tkzDrawSegments[thin](p00,p40 p40,p42 p42,p02 p02,p00) % box
      \tkzDrawSegments[thin](p01,p41 p10,p12 p20,p22 p30,p32) % grid
      \tkzDrawSegments[thin](p00,p11 p11,p20 p21,p30 p30,p41) % diagonals
      \tkzDrawSegments[thin](p02,p11 p11,p22 p21,p32 p32,p41) % diagonals
      %% draw vertices
      \tkzDrawPoints(p11,p31)
      \tkzLabelPoint[below right = 0 and 0.15](p11){$\bm{v}_1$}
      \tkzLabelPoint[below right](p31){$\bm{v}_2$}
   \end{tikzpicture}%
   \caption{%
      Выбор диагональных ребер для сеток Делоне для регулярной
      прямоугольной сетки не однозначен, так как четыре вершины
      прямоугольника принадлежат к одной окружности.  Жесткое правило
      для прямоугольной сетки, например 'снизу слева~-- наверх
      справа', позволяет простроить трансляционно-инвариантную
      сетку  (слева).  Некоторые диагональные ребра выбраны по
      другому (справа), что приводит к так называемой перекрестной
      схеме вокруг $\bm{v}_1$%
   }\label{fig:crisscross}%
\end{figure}

Для метода конечных элементов (МКЭ) ситуация иная, поскольку  здесь
топология сетки имеет значение.  С одной стороны, для сеток на
\cref{fig:crisscross}, матрица жесткости для  дискретизации с помощью
P\textsubscript{1} (линейных) конечных элементов инвариантна по
отношению к выбору диагональных ребер, так как в двухмерном случае,
матрица жесткости МКЭ идентична соответствующей ей матрице метода
конечных объемов, которая зависит только от координат узлов.  С
другой стороны, матрица масс МКЭ зависит от топологии сетки, например,
\[
   M_{ii} = \frac{2 \lvert \omega_{i} \rvert}{(d+1)(d+2)}
   ,
\]
где $\omega_{i}$~--- подобласть, относящаяся к узлу $\bm{v}_i$, а
$\abs{\omega_{i}}$ ее объем (площадь).  Сетка на
\cref{fig:crisscross}~(слева)~--- трансляционно-инвариантна (исключая
первый и последний узел), и каждая подобласть $\omega_{i}$ имеет
одинаковый размер.
Следовательно, справедливо
\[
   M_{11} = M_{22} 
   = \frac{2 \times \left(6 \times \frac{h^2}{2}\right)}{12}
   = \frac{1}{2} h^2
   .
\]
Сетка на \cref{fig:crisscross}~(справа), хотя и удовлетворяет условию
Делоне, но некоторые диагональные ребра были выбраны по другому.
Подобласти, принадлежащие $\bm{v}_1$ и $\bm{v}_2$, имеют разные
размеры, так же, как и соответствующие им элементы в матрице масс.
\[
   M_{11}
   = \frac{2 \times \left(8 \times \frac{h^2}{2}\right)}{12}
   = \frac{2}{3} h^2,
   \quad \text{тогда как} \quad
   M_{22}
   = \frac{2 \times \left(4 \times \frac{h^2}{2}\right)}{12}
   = \frac{1}{3} h^2
   ,
\]
так что матрица масс не будет трансляционно-инвариантной, в отличие
от матрицы жесткости.  Как следствие, численное решение потеряет
инвариантность в направлении $x$.  Более того, несогласованность
матрицы масс также может привести и к более низкому, чем ожидаемому,
порядку сходимости.  Одноточечная квадратура для построения матрицы массы не исправит
этот недостаток в общем, потому что распределение массы по узлам
зависит от топологии сетки и, таким образом, уже является несогласованной
еще до того, как будет проведено вычисление матрицы масс.  Подробное
обсуждение потери точности линейных МКЭ, а также конечных элементов
более высокого порядка в зависимости от топологии сеток содержится
в недавних работах Коптевой~\cite{Kop14,Kop20}.

%---------------------------------------------------------------------
\section{Примеры применения}
%---------------------------------------------------------------------

%--- Image segmentation ----------------------------------------------
\subsection{Разделение фаз}

Рассмотрим диффундирующие черные и белые частицы, сумма которых
везде равна 1 (статистика Ферми), короткодействующий потенциал взаимодействия 
(черное притягивает черное) и диффузию, действующую как
расталкивающую силу (такова качественная интерпретация
известного соотношения Эйнштейна).  Градиенты можно увеличить 
за счет увеличения потенциала притягивания по отношению к диффузии.

Концентрации $u_{b}$ черных и $u_{w}$ белых частиц удовлетворяют соотношению
\[
   1 > 1 - \varepsilon > u_{b} > \varepsilon > 0,
   \quad u_{w} = 1 - u_{b}
   .
\]
Если доминирует диффузионное расталкивание, то решением является
идеальное распределение серого.  В противном случае, черные и белые
частицы будут разделены минимальной поверхностью.

Ниже приводится краткое математическое описание этой нелокальной
задачи фазовой сегрегации; подробное обсуждение дано в
статье~\cite{GajGar05a} и в литературе к ней.  Задача описывается следующей системой уравнений
\begin{align}
   \begin{cases}
      \frac{\partial u}{\partial t}
         - \nabla \cdot \left(
            f (\abs{\nabla v})
            \left( \nabla u + \frac {\nabla w}{\phi''(u)} \right) \right) = 0
      \quad \text{в $\Omega$}
      ,
      \\
      u(0,\cdot) = u_0(\cdot)
      ,
      \\
      v = \phi'(u) + w
      ,
      \\
      w(t,\bm{x}) = \int_\Omega \mathcal{K} \left( \abs{\bm{x}-\bm{y}} \right) 
         \left(1 - 2 u(t,\bm{y}) \right) \D{\bm{y}}
      .
   \end{cases}
   \label{eq:phaseSep}
\end{align}
где $\Omega \subset \R^N$, $1 \le N \le 3$~--- ограниченная Липшицева
область, $\phi$~--- выпуклая функция, ядро $\mathcal{K}$ представляет
нелокальные силы притяжения, а $w$ и $v$~--- потенциал взаимодействия
и химический потенциал, соответственно.  Первоначальное значение
$u_0$ удовлетворяет неравенству $0< u_0 <1$. Система дополняется
однородными граничными условиями Неймана.

Эта система была выведена в~\cite{GiaLeb97} для случая постоянной
функции $f$ и может быть рассмотрена как нелокальный вариант уравнения
Кан-Хиллиарда, связанного с локальной свободной энергией Гинзбурга-Ландау
\[
   F_{gl}(u) = \int_\Omega \left(
      \phi(u) + \kappa  u \left(1 - u\right)
      + \frac{\lambda}{2} \abs{\nabla u}^2
   \right) \D{\bm{x}}
   .
\]
Функционалом Ляпунова для \cref{eq:phaseSep} является
\[
   F(u) = \int_\Omega \left(
      \phi(u) + u \int_\Omega \mathcal{K} 
   \left(\abs{x-y}\right) \left(1 - u(y) \right) \D{y} \right) \D{\bm{x}}
   .
\]
Рассмотрим специальную модель,  в которой функция $f$ постоянна,
$\phi(u) =  u \ln u + (1 - u) \ln (1 - u)$,
и $\mathcal{K}$~--- это функция Грина для эллиптической краевой задачи
\begin{align*}
   \begin{cases}
      - \sigma^2 \nabla \cdot \nabla w + w = m \left(1 - 2u\right)
      & \text{в $\Omega$}
      ,
      \\
      \nu \cdot \nabla w = 0
      & \text{на  $\partial \Omega$}
      .
   \end{cases}
\end{align*}
с константами $\sigma, m > 0$.
При расчете на малые времена,  эта модель также используется для реконструкции
изображений и подавления шума, и не ограничивается только двумя видами
частиц.\cite{GajGar05,GajGri06,Gri04}.

Подходящим численным методом решения этой проблемы является схема
\begin{align*}
   \left(G^T \sigma^2 G +[V]\right) \bm{w}^+ 
      &= [V] m \left(\bm{1} -2 \bm{u}^+ \right)
      ,
      \\
   [V] \frac{\bm{u} - \bm{u}^+} \tau
      & = G^T f \left(
         G \bm{u}^+ + \left[ \frac{1}{ \Phi''_a} \right]
         \frac{1}{2} G \left( \bm{w}^+ + \bm{w} \right) \right)
      ,
\end{align*}
Доказано, что данная схема является диссипативной~\cite[Теорема~2.7]{GajGar05a}.
Более того, если $\bm{v}^+ \ne c$, $0 < u < 1$, то неявный метод Эйлера
сохраняет свойство $0 < u^+ < 1$~\cite[Утверждение~2.9]{GajGar05a}

Контроль размера шага по времени в вычислениях основан на разностях типа предиктор-корректор для 
свободной энергии и дискретной скорости
диссипации.  Матрицы Якоби выводятся при помощи аналитического дифференцирования
на каждом симплексе, а метод Ньютона применяется без демпфирования.
В качестве начального решения берется
\[
   u_0(\bm{x}) = \varepsilon + \frac{x_1}{128 \left( 1 + 2\varepsilon \right)}
  ,
  \quad \Omega = (0,128) \times (0,128),
  \quad m = 8,
  \quad \varepsilon=10^{-7},
  \quad \sigma = \sqrt{f} = 2
  .
\]

\begin{figure}[t]%
   \includegraphics[width=0.3\textwidth]{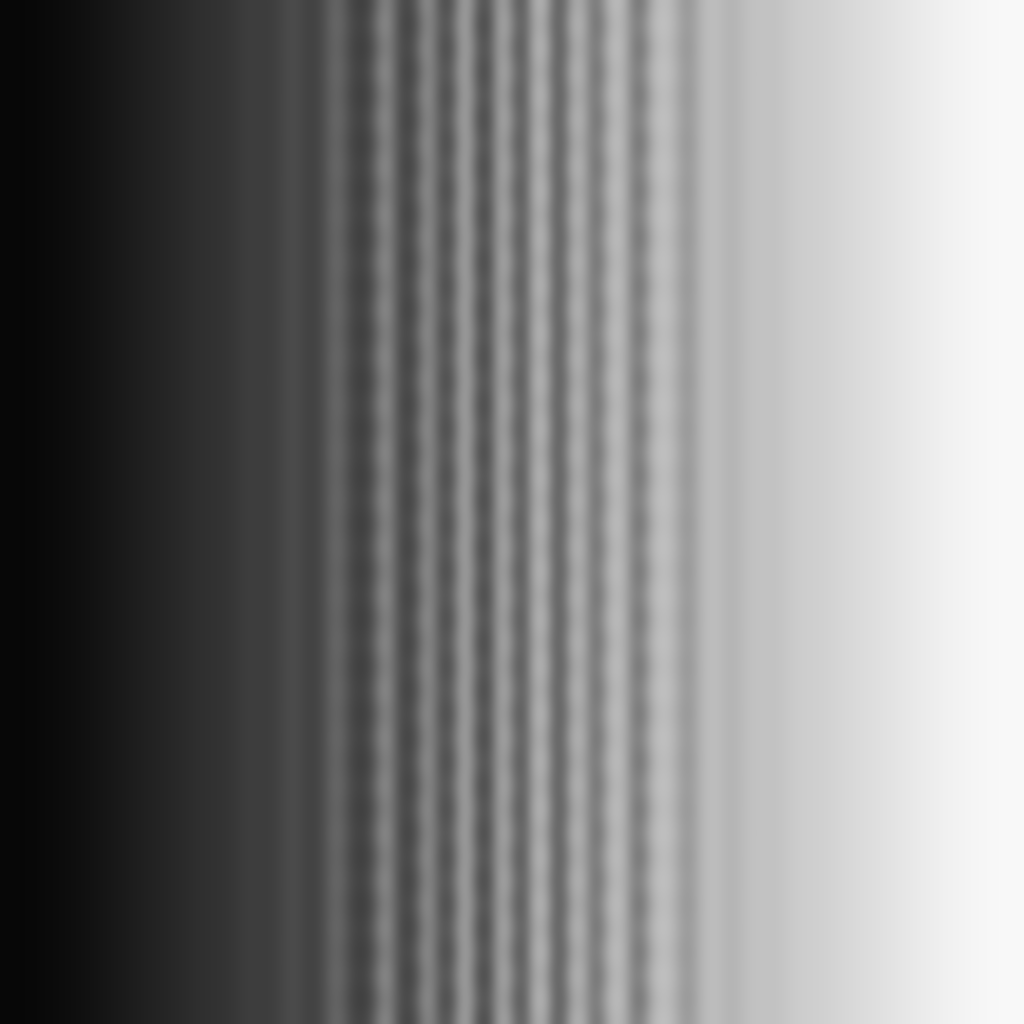}%
   \hfill{}%
   \includegraphics[width=0.3\textwidth]{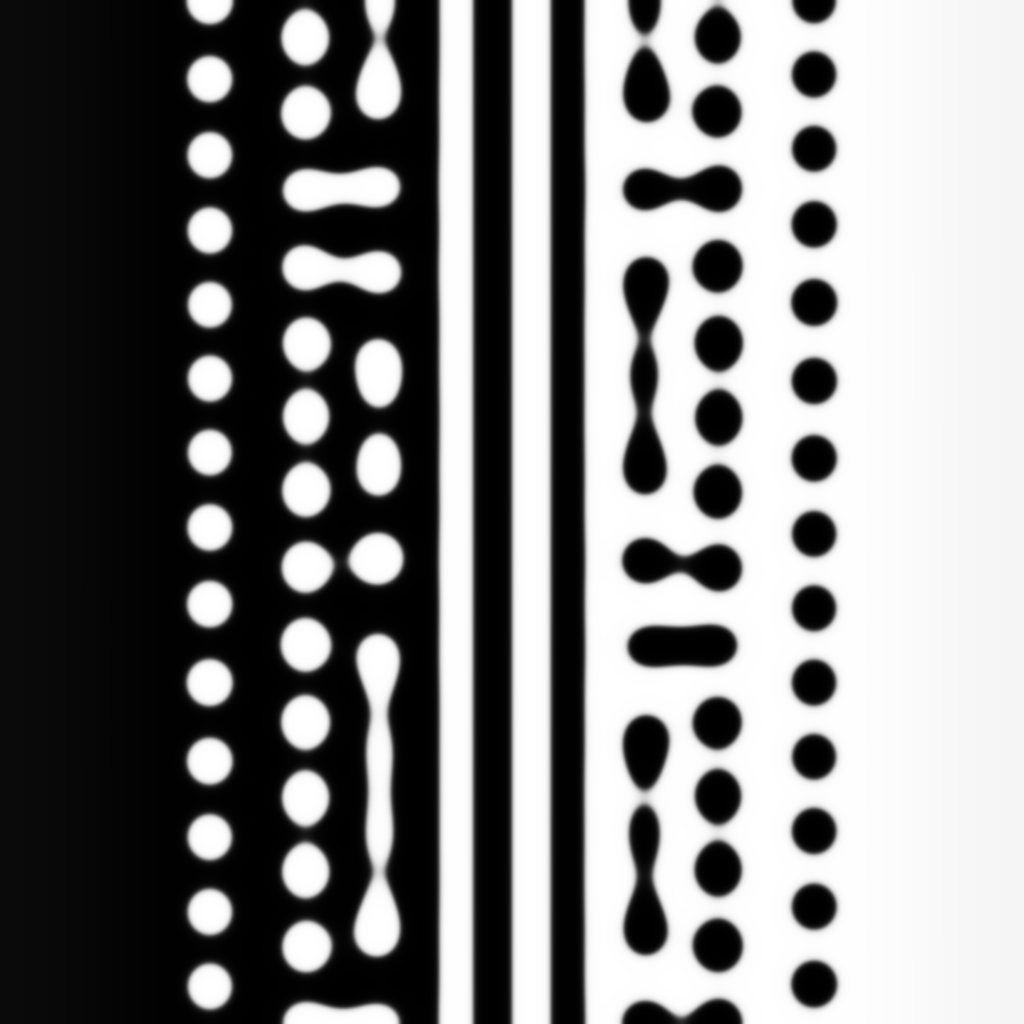}%
   \hfill{}%
   \includegraphics[width=0.3\textwidth]{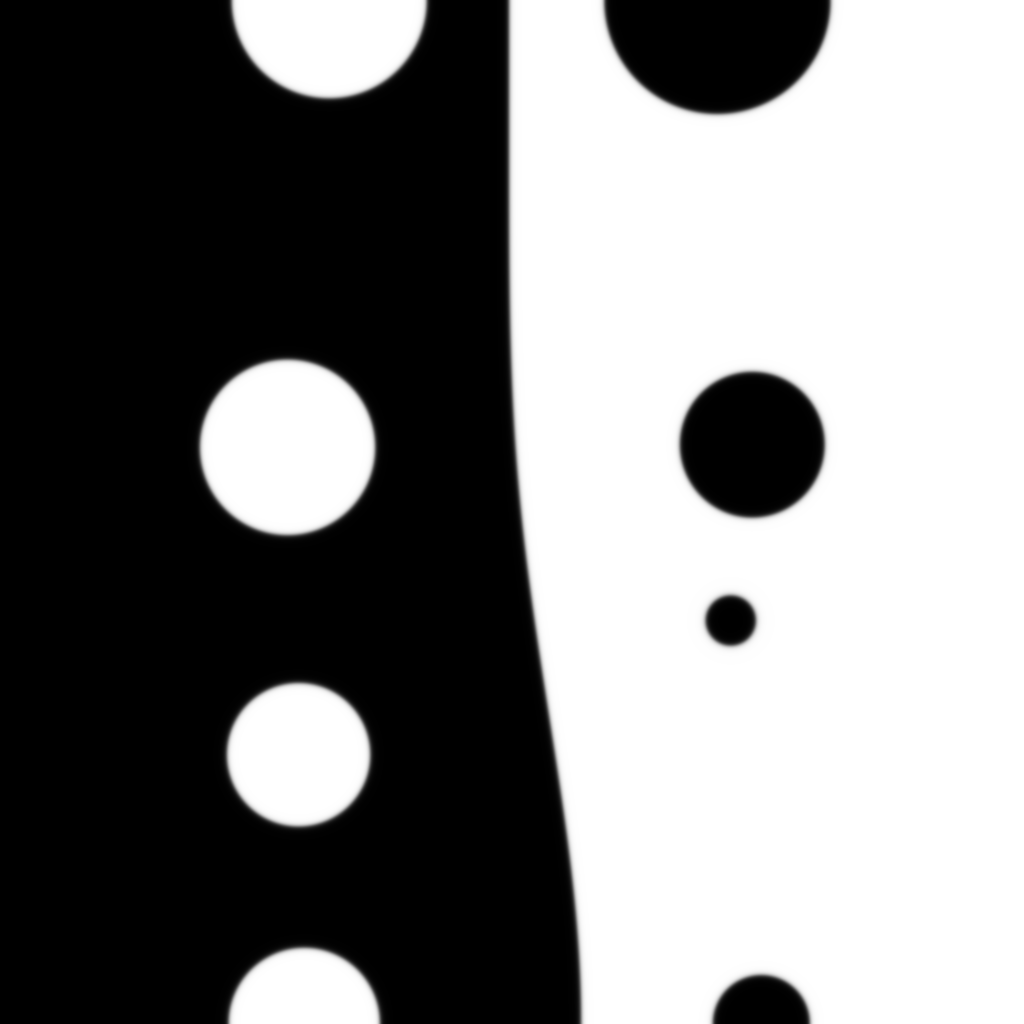}%
   \caption{%
   Изначальное серое распределение черного и белого очень быстро
   переходит в мелкие структуры (слева).  Ошибки округления
   нарушают симметрию,   кривизна объектов и расстояние между ними
   определяют эволюцию (в центре).  Объект с наименьшим диаметром
   исчезнет следующим (справа); устойчивая прямая (конечная) граница
   между черным и белым все еще далека на линейной шкале времени%
   }\label{fig:bubbles}%
\end{figure}

\begin{figure}[t]%
   \centering{}%
   \includegraphics[height=0.26\textheight, clip]{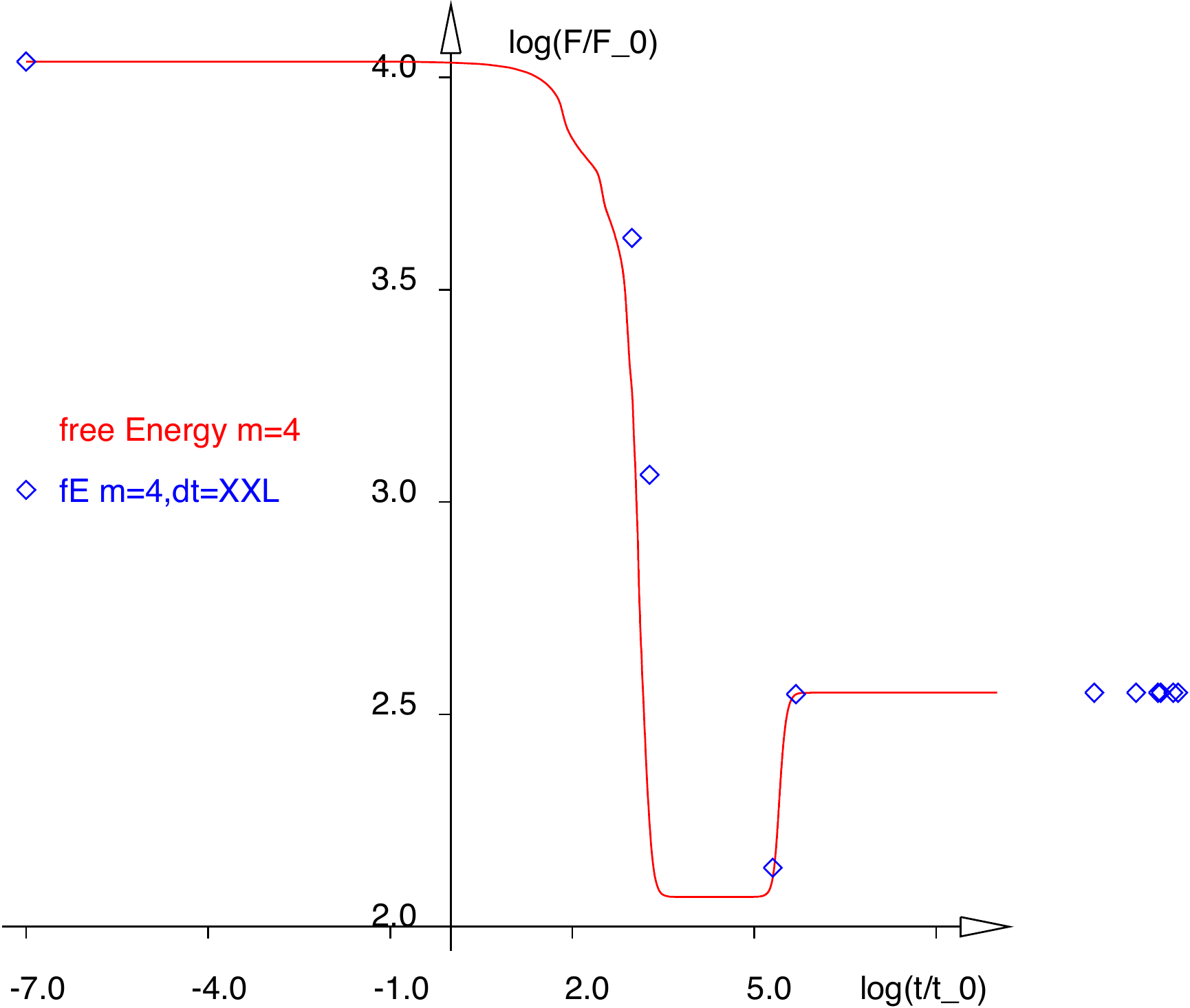}%
   \hspace{0.1\linewidth}%
   \includegraphics[width=0.26\textheight, clip]{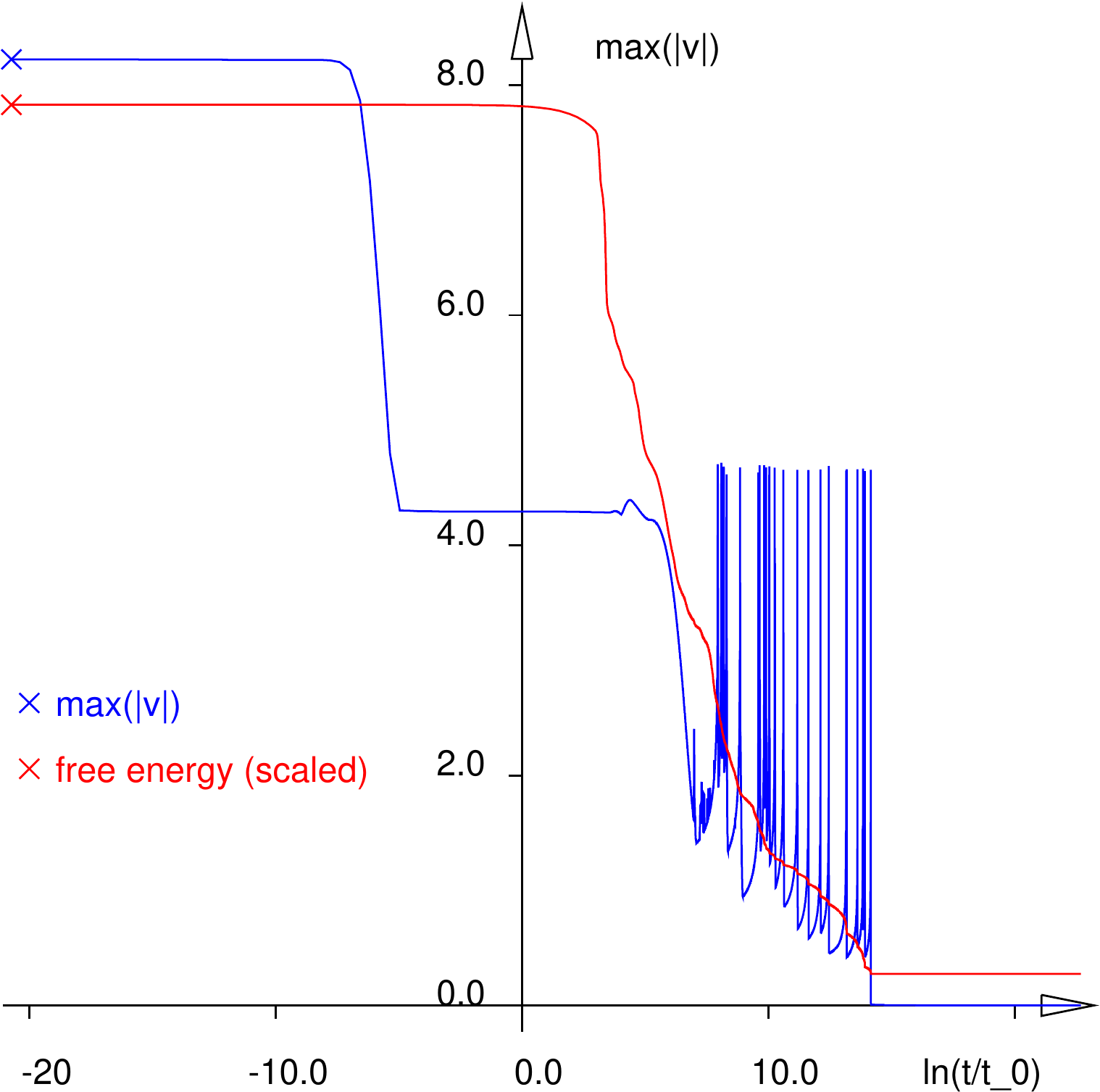}%
   \caption{%
      Свободная энергия увеличивается для недиссипативной дискретизации
      (для малых и больших шагов по времени, слева) и окончательно
      падает ступенчато для дискретизации, являющейся диссипативной
      и сохраняющей ограничения значения решения (справа).  Каждый
      пик производной по времени (синий цвет, одинаковая высота
      означает одинаковую максимальную скорость исчезновения)
      показывает исчезновение белого или черного объекта в области
      другого цвета, к конце очень быстрое по сравнению к другим
      шкалам времени.  Выравнивание первоначально изогнутой конечной
      разделительной линии определяет самую большую шкалу времени%
   }\label{fig:bubbles:plot}%
\end{figure}

Диссипативная схема, как и ожидалось, выдает приближение минимальной
поверхности.  После $10^{20}$ по времени, исходное распределение
серого преображается, проходя через узор из полосок с пузырьками, в
черные и белые половинки квадрата, разделенные прямой вертикальной
линией посередине квадрата.  \Cref{fig:bubbles,fig:bubbles:plot}
показывают типичный процесс эволюции: тонкие полосы (\cref{fig:bubbles},
слева) распадаются на пузырьки с диаметром, связанным с шириной
полосы (\cref{fig:bubbles}, в центре).  Симметрия нарушается из-за
ошибок округления, небольших неточностей в исходных данных, а также
из-за того, что область $\Omega$ имеет большой размер по сравнению
с $\sigma$ и не совместима с произвольной шириной полос.  Самая
длинная шкала времени определяется временем, необходимым для
выравнивания границы между этими половинами до прямой линии.

Если дискретизация не является диссипативной, то асимптотическая
минимальная поверхность будет заменена поверхностью большей площади.
В начальный момент времени свободная энергия максимальна.
Ближе к концу она падает и,
через несколько порядков во времени, выглядит постоянной.  В случае
недиссипативной схемы, оставшийся градиент решения, очень маленький почти
всюду, достаточен для того, чтобы производить немного свободной
энергии в каком-либо месте сетки. Последствием этого являются бóльшие
градиенты и большая генерация свободной энергии. 
Поскольку объем области в основном диссипативен,  это приводит к разделяющей поверхности большей площади,
потому что сильное рассеивание вдоль поверхность необходимо для
компенсации производства свободной энергии.  В конце концов,
полученное устойчивое состояние может быть весьма далеко от истинного
состояния минимальной энергии.


%--- Athena-WFI ------------------------------------------------------
\subsection{DePFET детекторы с активными пикселями}

Детекторы DePFET с активными пикселями~--- это кремниевые сенсоры
на базе обедненных p-канальных полевых транзисторов
(DePFETs)~\cite{KemLut87,LutAndEck07}, которые представляют собой
полевые транзисторы (FETs), интегрированные в обедненный кремниевый
объем.  Входящие частицы (фотоны) генерируют электроны в объеме,
которые затем собираются на внутреннем затворе под транзисторным
каналом.  Это увеличивает проводимость транзистора пропорционально
количеству накопленного заряда, обеспечивая этим измерение энергии
фотонов.  Уникальными преимуществами детекторов, основанных на этой
технологии, являются усиление и хранение заряда сигнала непосредственно
внутри пикселей, неразрушающее считывание, очень низкий уровень
шума и высокая энергоэффективность, что делает их идеальными для
рентгеновской астрономии и физики частиц~\cite{LutAndEck07,TreAndHau18}.

Примером такого детектора является широкоугольная камера (\foreignlanguage{english}{Wide Field Imager, WFI}) для усовершенствованного телескопа для высокоэнергетической астрофизики (\foreignlanguage{english}{Advanced Telescope for High-Energy Astrophysics, Athena})~\cite{athenaUrl,athena} Европейского Космического Агентства.
WFI~--- это твердотельный рентгеновский спектрометр для рентгеновских снимков в полосе 0.1--15\,кВ с одновременным спектральным и временным разрешением фотонов.
Он разработан в совместной работе Лаборатории Полупроводников Общества Макса-Планка (MPG HLL), Эрлангенского Центра Физики Астрочастиц, Института Астрономии и Астрофизики Тюбингена, Дрезденского Технического Университета, PNSensor GmbH, Лейстерского Университета и Научно-исследовательского Института Астрофизики и Планетологии в Тулузе~\cite{athena}.
Численное моделирование проводилось в сотрудничестве MPG HLL с Институтом Вейерштрасса по Прикладному Анализу и Стохастике.

В полупроводниках, по сравнению с предыдущим примером, взаимодействие частиц обратное: частицы разного типа (отверстия и электроны) притягивают друг друга, а диффузия является силой отталкивания для частиц одного и того же типа.
Математическая модель задается уравнениями ван~Росбрека:
\begin{align*}
   -\nabla \cdot \varepsilon \nabla w 
   &= C - n - p
   ,\\
   \frac{\partial n}{\partial t} - \nabla \cdot n_i \mu_n e^w \nabla e^{-\phi_n} 
   &= R
   ,\\
   \frac{\partial p}{\partial t} - \nabla \cdot n_i \mu_p e^{-w} \nabla e^{\phi_p}
   &= R,
\end{align*}
где $w$~--- электростатический потенциал, $n = n_i e^{w - \phi_n}$ и  $p = n_i e^{\phi_p - w}$~--- плотности электронов и дырок, $C$~--- легирование, $\phi_n$ and $\phi_p$~--- квазиуровни Ферми, и $R = r(\bm{x},n,p) \left( n_i^2 - np \right)$ описывает рекомбинацию/генерацию ($r(\bm{x},n,p)>0$).

\begin{figure}[p]%
   \includegraphics[width=0.47\textwidth]{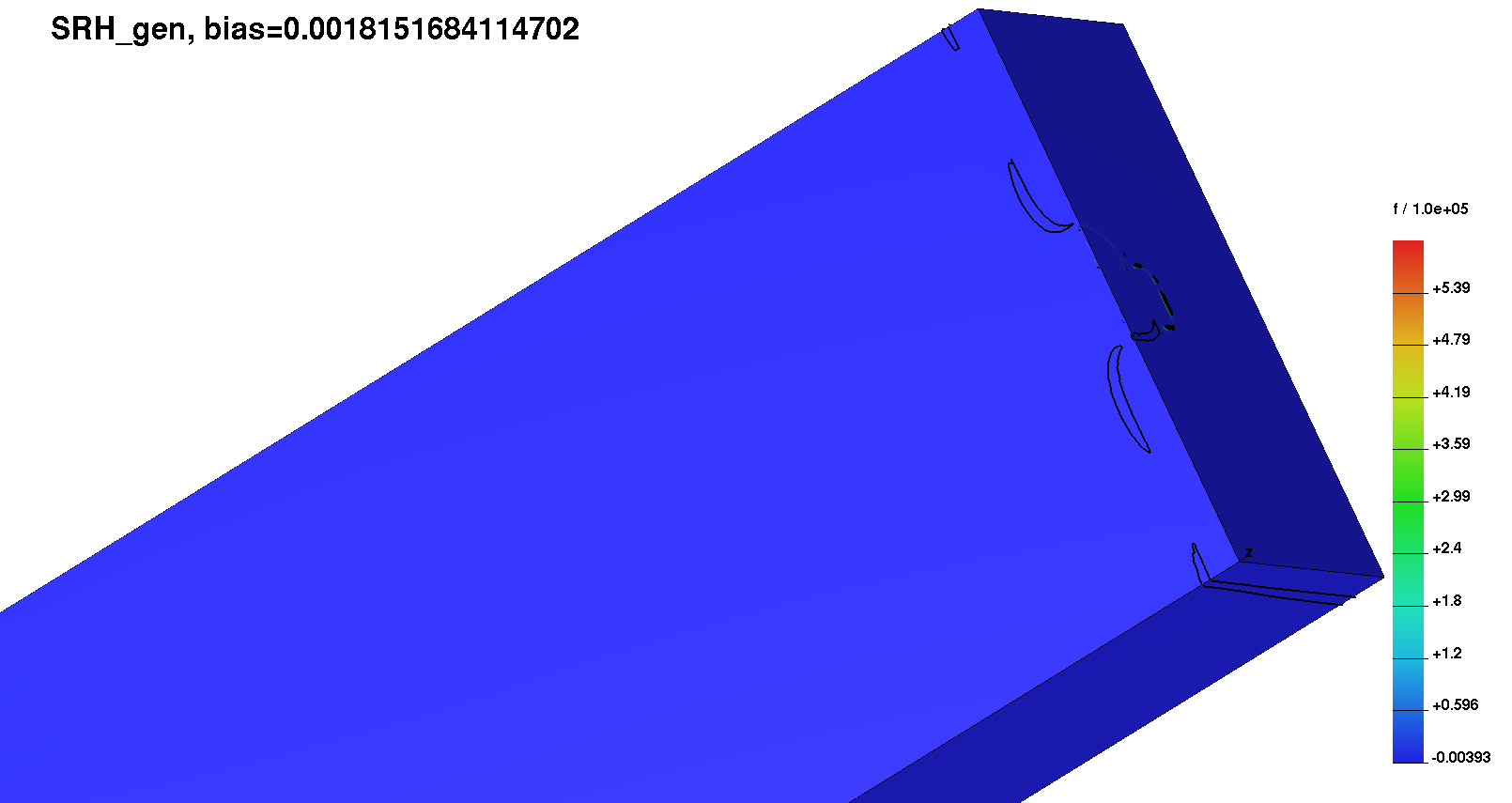}%
   \hfill{}%
   \includegraphics[width=0.47\textwidth]{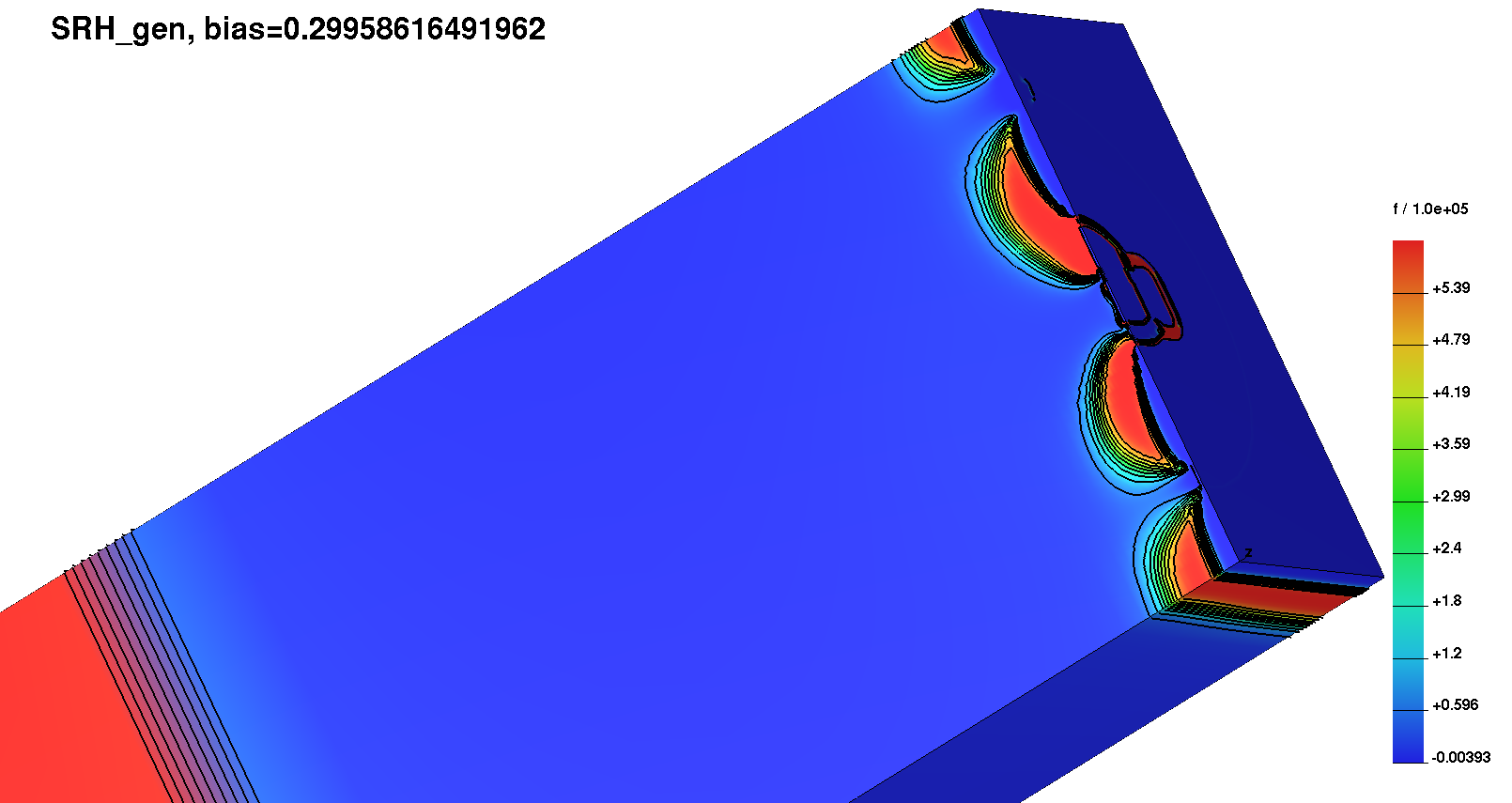}%
   \\[0.25\baselineskip]%
   \includegraphics[width=0.47\textwidth]{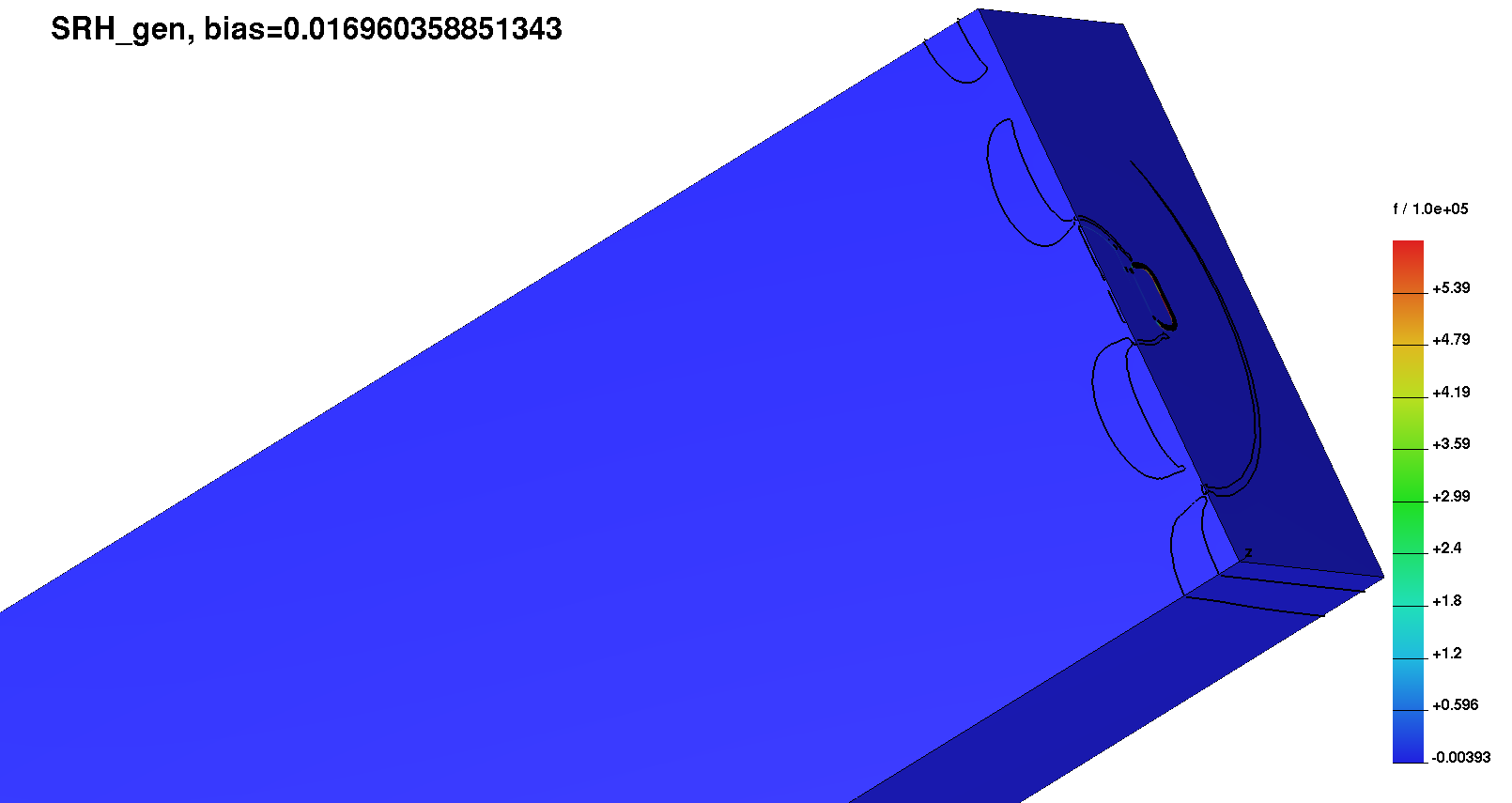}%
   \hfill{}%
   \includegraphics[width=0.47\textwidth]{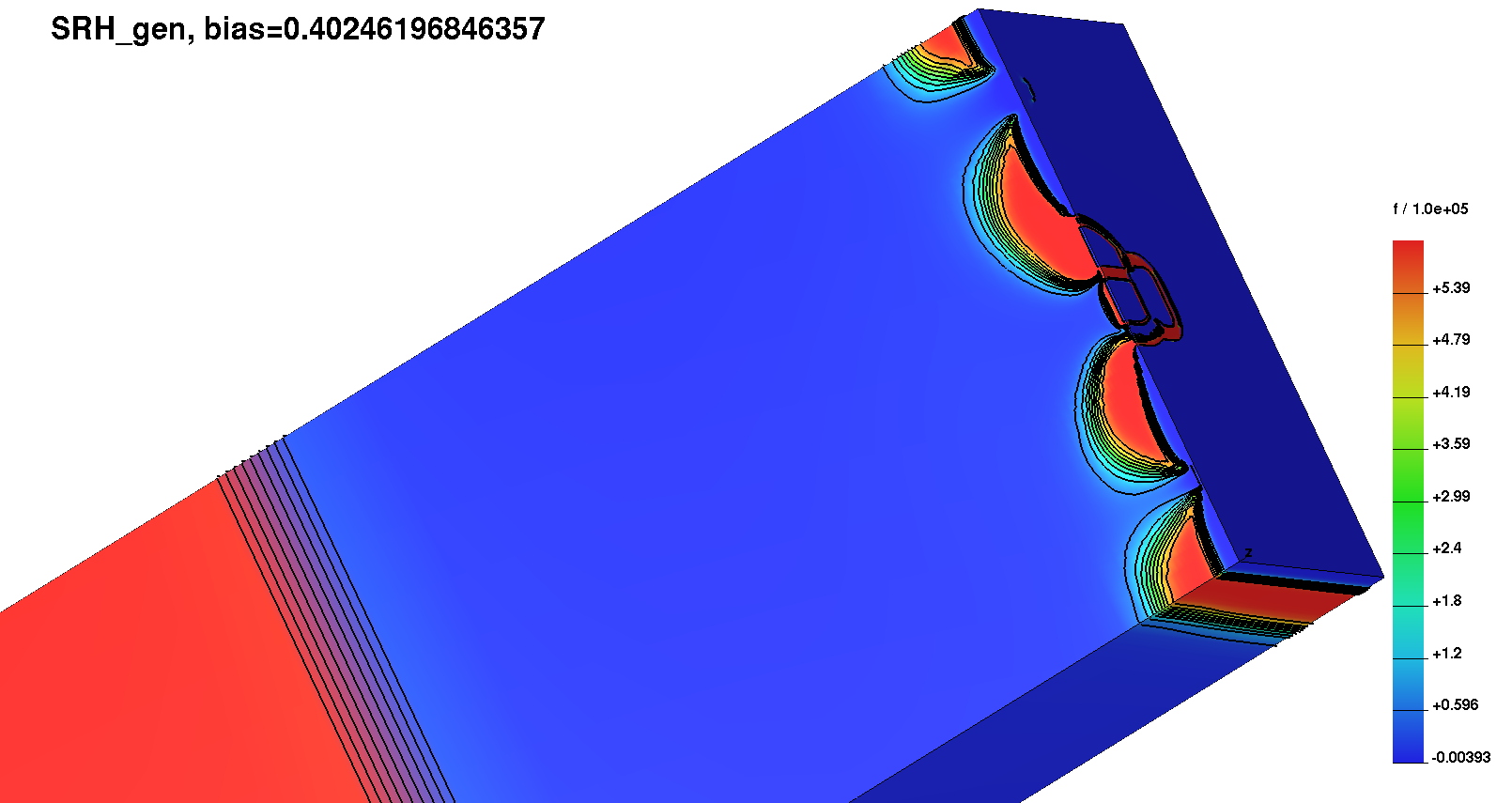}%
   \\[0.25\baselineskip]%
   \includegraphics[width=0.47\textwidth]{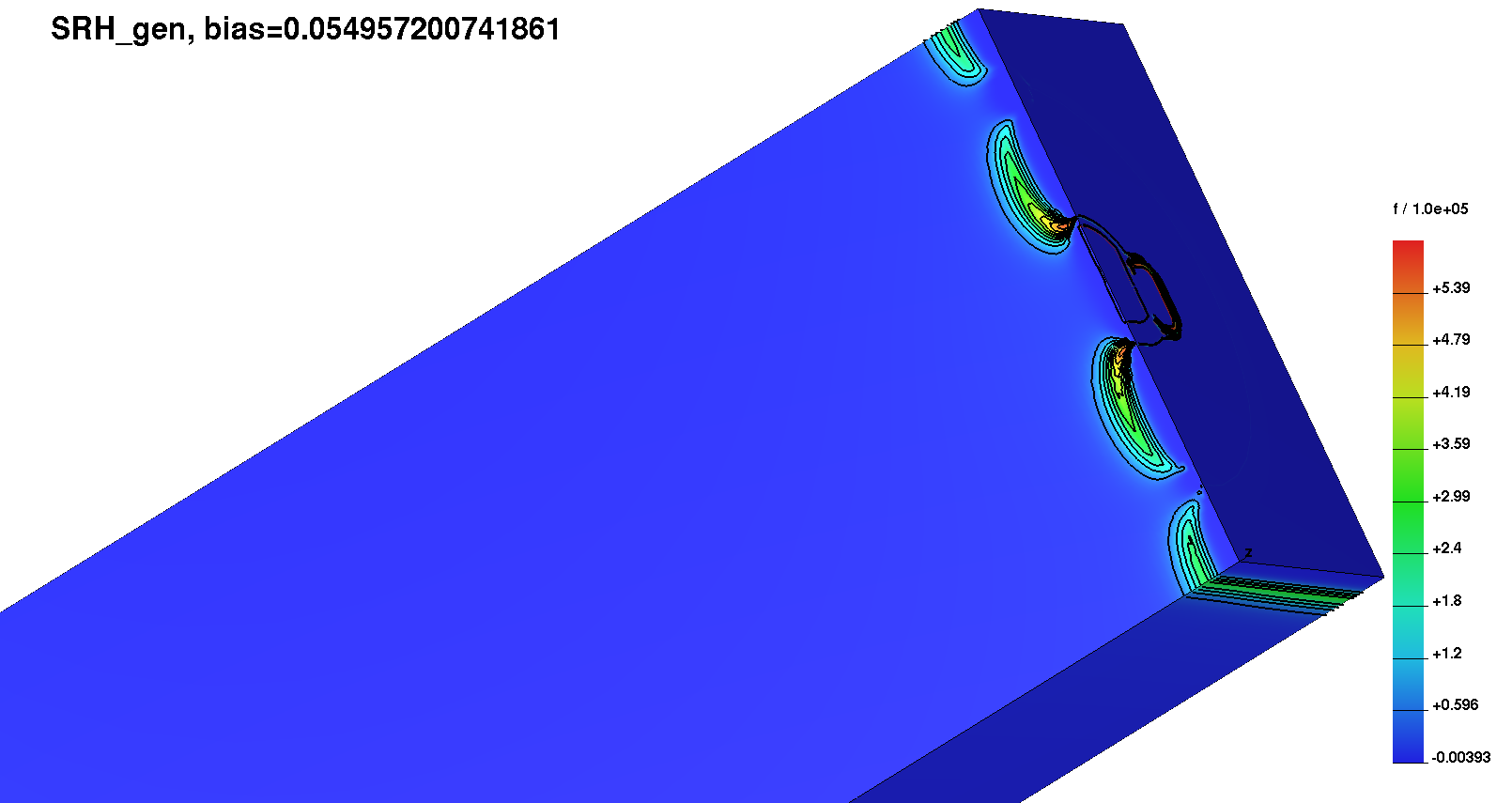}%
   \hfill{}%
   \includegraphics[width=0.47\textwidth]{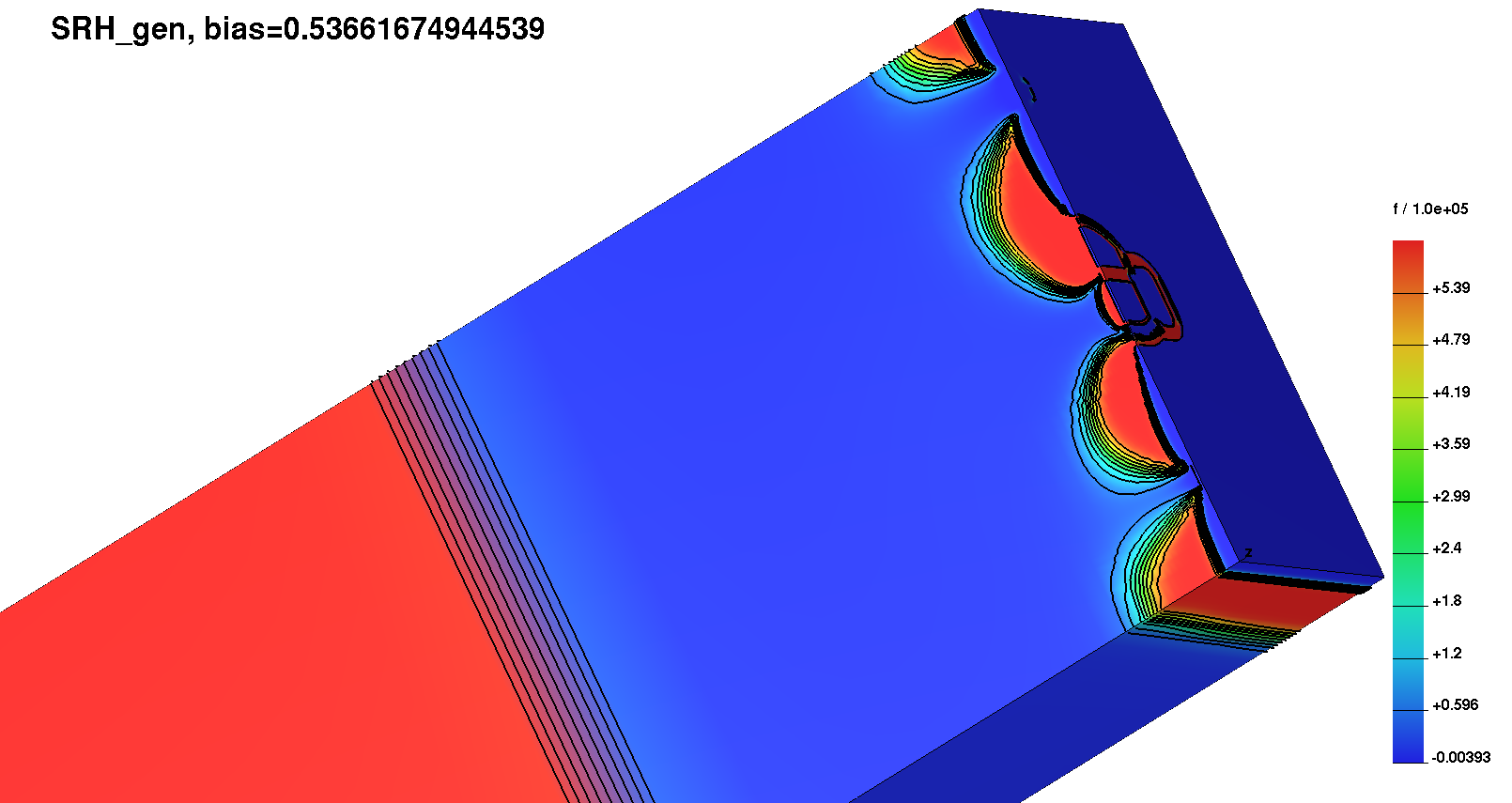}%
   \\[0.25\baselineskip]%
   \includegraphics[width=0.47\textwidth]{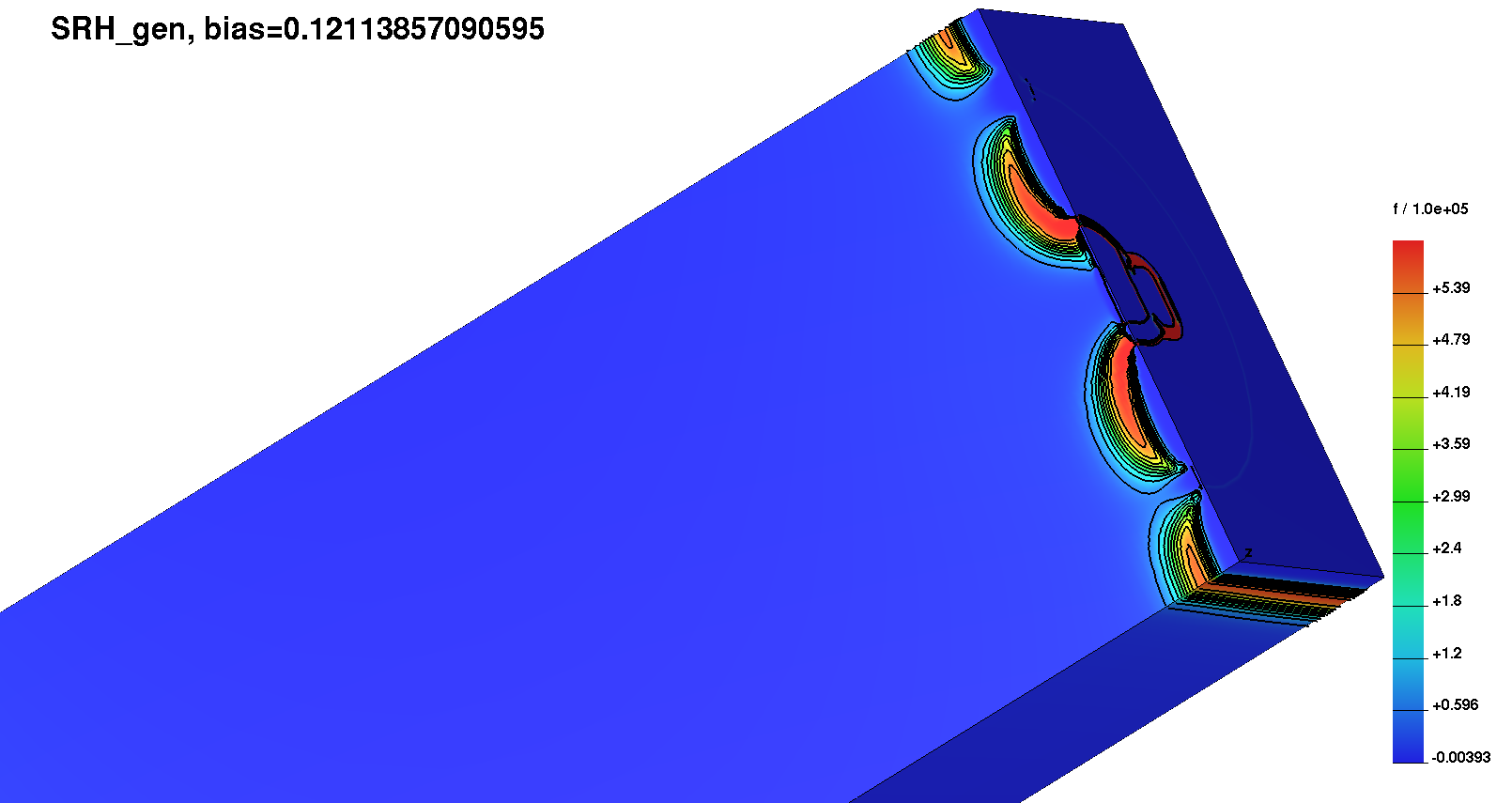}%
   \hfill{}%
   \includegraphics[width=0.47\textwidth]{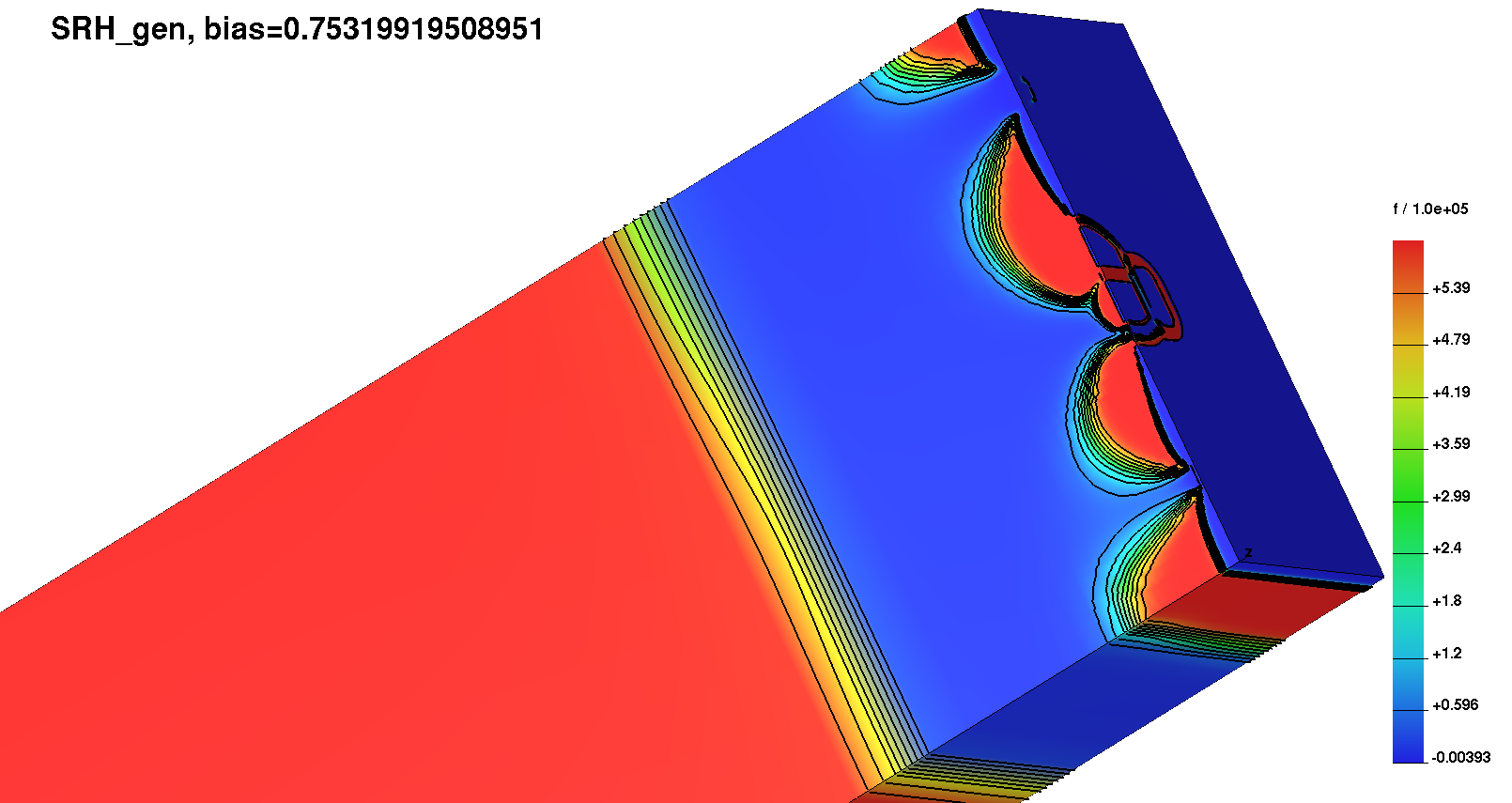}%
   \\[0.25\baselineskip]%
   \includegraphics[width=0.47\textwidth]{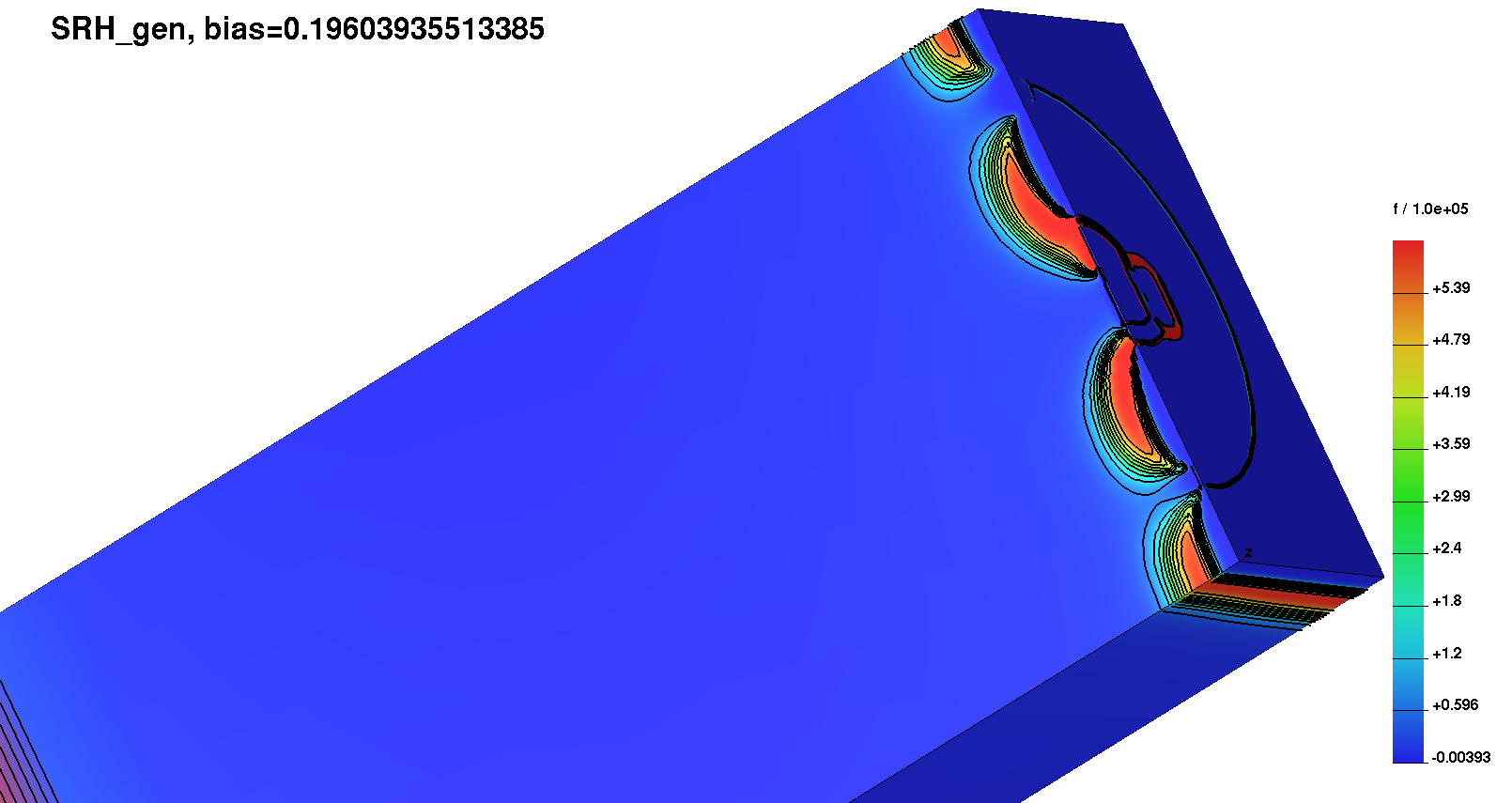}%
   \hfill{}%
   \includegraphics[width=0.47\textwidth]{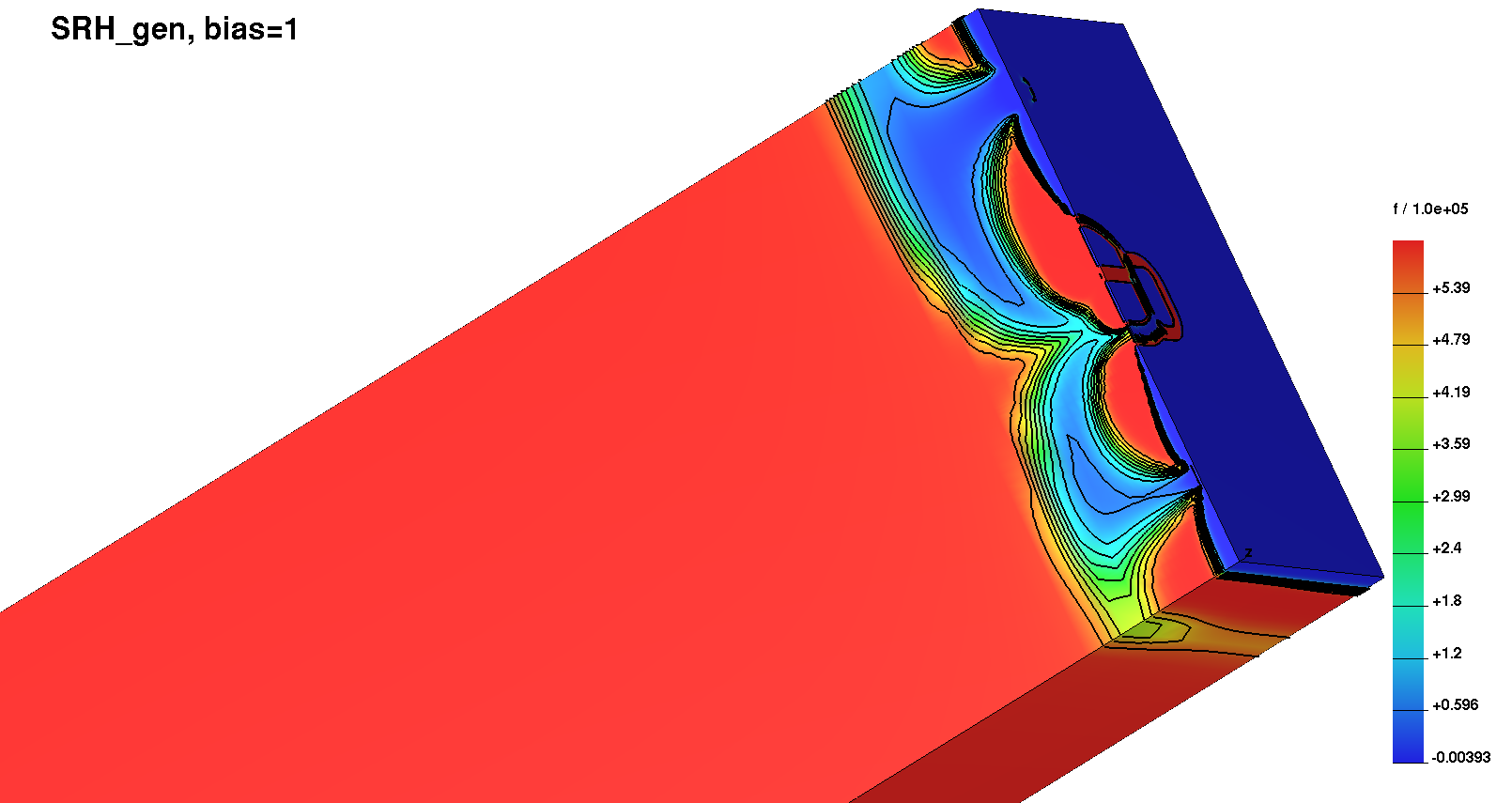}%
   \caption{Детектор Athena WFI: моделирование распространения фронта деплеции%
   }\label{fig:athena}%
\end{figure}

Для численного моделирования важно, чтобы дискретная система сохраняла существенные свойства аналитической системы, такие как принцип максимума, положительную концентрацию носителей и сохранение тока. VFVM позволяет перенести эти свойства в дискретную модель.
Численные расчеты производились по хорошо известной схеме Шарфеттера-Гуммеля (подробное описание находится, например, в~\cite{Gar09}).
\cref{fig:athena} показывает моделирование распространения фронта деплеции в пикселе детектора.
Ввиду симметричности области, вычисления выполнены на половине пикселя.

%--- Summary ----------------------------------------
\section{Заключение}

Сетки Делоне-Вороного  вполне пригодны для корректного описания
геометрии задачи, при этом предоставляя возможность сохранять
основные качественные свойства аналитического решения.  Трудно
представить себе какой-то способ для достижения большего, чем расчет
на короткие времена, численными методами, которые нарушают основные
свойства устойчивости аналитической задачи такие, как ограничение
значений решения, а также диссипативность для любых концентраций и
любого разрешения в пространстве и времени.

Знание не всегда оказывается строго возрастающей функцией времени,
и связь между аналитическими свойствами (нелинейных) систем уравнений
в частных производных, методами конечного объема, сетками Делоне,
и генерацией сеток в прошлом могла быть сильнее  чем в настоящее
время.  В течение последних десятилетий сетки Делоне, согласованные
с границей, стали ``редкими видами'', но есть надежда на их грядущее
возрождение~\cite{GarKudTsv19}.

\newpage{}